\documentclass[12pt]{article}
\usepackage{mathrsfs}
\usepackage{amsmath,amssymb}

\usepackage{amsfonts}
\usepackage{latexsym}
\usepackage{amsthm}
\usepackage{pictex}
\usepackage{color}

\usepackage[T2A]{fontenc}     
\usepackage[cp1251]{inputenc} 
\usepackage{graphicx}

\def\const{{\mathrm{const}}} 
\def\sgn{\mathop{\mathrm{sgn}} \nolimits} 

\def\const{constant}

\def\e{equation}

\def\neigh{neighborhood{}}

\def\pseudor{pseudodifferential operator{}}

\def\wrt{with respect to}

\def\Re{{\rm Re\,}}
\def\Im{{\rm Im\,}}
\def\lb{\label}
\def\dSRN{de Sitter-Reissner-Nordstr{\"o}m }

\newtheorem{theorem}{Theorem}
\newtheorem{definition}{Definition}
\newtheorem{lemma}{Lemma}

\newtheorem{proposition}{Proposition}

\newtheorem{remark}{Remark}

\def\g{\gamma}    
\def\G{\Gamma} \def\cD{{\cal D}}   
\def\d{\delta} \def\cE{{\cal E}}   
\def\D{\Delta} \def\cF{{\cal F}}   
   \def\cG{{\cal G}}   
\def\z{\zeta}  \def\cH{{\cal H}}   
\def\e{\eta}

\def\k{\kappa} \def\cL{{\cal L}}   
\def\l{\lambda}\def\cM{{\cal M}}   
\def\L{\Lambda}   
       
    \def\cP{{\cal P}}   
   \def\cQ{{\cal Q}}   
\def\s{\sigma} \def\cR{{\cal R}}   
\def\S{\Sigma}    
      
   \def\cU{{\cal U}}

\def\o{\omega}    
\def\x{\xi}       
    \def\cZ{{\cal Z}}   
\def\O{\Omega}

\def\vt{\vartheta}
\def\vp{\varphi}

\def\ga{\mathfrak{a}}
\def\gb{\mathfrak{b}}
\def\gq{\mathfrak{q}}

\newcommand{\gD}{\mathfrak{D}}

\newcommand{\gm}{\mathfrak{m}}

\def\Z{{\Bbb Z}}
\def\R{{\Bbb R}}
\def\C{{\Bbb C}}

\def\N{{\Bbb N}}
\def\S{{\Bbb S}}

\def\qq{\quad}
\newcommand{\ma}{\begin{pmatrix}}
\newcommand{\am}{\end{pmatrix}}
\newcommand{\ca}{\begin{cases}}
\newcommand{\ac}{\end{cases}}

\let\geq\geqslant
\let\leq\leqslant
\def\ma{\left(\begin{array}{cc}}
\def\am{\end{array}\right)}

\let\geq\geqslant
\let\leq\leqslant
\def\[{\begin{equation}}
\def\]{\end{equation}}

\def\no{\noindent}

\def\/{\over}

\def\os{\oplus}

\def\Re{\mathop{\rm Re}\nolimits}
\def\Im{\mathop{\rm Im}\nolimits}

\def\sign{\mathop{\rm sign}\nolimits}

\def\const{\mathop{\rm const}\nolimits}

\def\BBox{\hspace{1mm}\vrule height6pt width5.5pt depth0pt \hspace{6pt}}

\begin{document}

\numberwithin{equation}{section}
\numberwithin{theorem}{section}
\numberwithin{proposition}{section}
\numberwithin{definition}{section}
\numberwithin{remark}{section}

\date{\today}

\title{Quasi-normal modes for  Dirac fields in Kerr-Newman-de Sitter black holes}

\author{
Alexei Iantchenko
\begin{footnote}
{Malm{\"o} h{\"o}gskola, Teknik och samh{\"a}lle, 205 06
  Malm{\"o}, Sweden, email: ai@mah.se }
\end{footnote}
}

\maketitle

\begin{abstract} We provide the full asymptotic description of the quasi-normal modes (resonances) in any strip of fixed width for  Dirac fields in slowly rotating Kerr-Newman-de Sitter black holes.  The resonances split in a way similar to the Zeeman effect. The method is based on the extension to Dirac operators of techniques applied by Dyatlov  in \cite{Dyatlov2011}, \cite{Dyatlov2012} to the (uncharged) Kerr-de Sitter black holes.  We show that the mass of the Dirac field does not have effect on the two leading terms in the expansions of resonances.

\end{abstract}


\section{Introduction.}

\subsection{Background.} Kerr--Newman-de Sitter (KN-dS) black hole  is an  exact solution of the Einstein-Maxwell equations which describe electrically charged  rotating black hole  with positive cosmological constant. 

We refer to \cite{DaudeNicoleau2015}  for detailed study in this background, including complete time-dependent scattering theory.

In Boyer-Lindquist coordinates, the exterior of a KN-dS black hole  is described by the four-dimensional manifold $$\R_t\times\cM,\qq \cM= ]r_-,r_+[\times \S_{\vartheta,\varphi}^2,$$ equipped with the Lorentzian metric
 (see Eq. (1.2) in \cite{DaudeNicoleau2015})
\[\lb{(1.2)}g=\frac{\Delta_r}{\rho^2}\left[dt-\frac{a\sin^2\vartheta}{E}d\vp\right]^2-\frac{\rho^2}{\Delta_r}dr^2-\frac{\rho^2}{\Delta_\vt}d\vt^2-\frac{\Delta_\vt\sin^2\vt}{\rho^2}\left[adt-\frac{r^2+a^2}{E}\right]^2,\]
where $$ \rho^2=r^2+a^2\cos^2\vt,\qq E=1+\frac{a^2\L}{3},$$
$$\Delta_r=(r^2+a^2)\left(1-\frac{\L r^2}{3}\right)-2Mr+Q^2,\qq\Delta_\vt=1+\frac{a^2\L\cos^2\vt}{3}.$$ Here  the parameters $M>0,$ $Q\in\R$ and $a\in\R$  are interpreted as the mass, the electric charge and the angular momentum per unit mass of the black hole, and $\L>0$ is the cosmological constant of the universe. In the case of KN-dS black hole the function $\Delta_r$ has three simple positive roots $0<r_c<r_-<r_+$ and one negative root
$r_n=-(r_c+r_-+r_+)<0,$ under the condition that  (see \cite{DaudeNicoleau2015})
\[\lb{2.5-2.6}\frac{a^2}{3}\L\leq 7-4\sqrt{3}\approx 0.072,\qq M_{\rm crit}^- <M<M_{\rm crit}^+,\]
where
\begin{align*}&M_{\rm crit}^\pm=\frac{1}{\sqrt{18\L}}\left( E_-\pm\sqrt{E_-^2-F}\right)^2\left( 2E_-^2\mp\sqrt{E_-^2-F}\right),\\ & E_-=1-\frac{a^2}{3}\L,\qq F=4\L(a^2+Q^2). \end{align*}

In this paper we always assume conditions (\ref{2.5-2.6}) to be fulfilled.
By introduction of the Regge-Wheeler coordinate $x$ via
\[\lb{Regge-W}\frac{d x}{dr}=\frac{r^2+a^2}{\Delta_r}, \] the event and cosmological horizons are  pushed away to $\{x=-\infty\}$ and $\{x=+\infty\}$ respectively. In this article we consider the exterior region of the black hole 
\[\lb{extreg} r_-<r<r_+ \qq\Leftrightarrow\qq -\infty <x<\infty.\]

We consider propagation of the  Dirac fields  with charge  $\gq$ and mass $\gm$ in the exterior region of the black hole. 

Let
\[\lb{3.5} \ga(x)=\frac{\sqrt{\Delta_r}}{r^2+a^2},\qq c(x,D_\varphi)=\frac{aE}{r^2+a^2}D_\varphi+ \frac{\gq Qr}{r^2+a^2},\qq \gb(x)=\gm\frac{r\sqrt{\Delta_r}}{r^2+a^2}.\]  
Note that in the non-rotating case $a=0$ (\dSRN black hole, the massless and chargeless case  was studied in \cite{Iantchenko2015}),  $\Delta_r=r^2F(r),$ where $F(r)=1-\frac{2M}{r}+\frac{Q^2}{r^2}-\frac{\Lambda}{3}r^2$ and $\ga^2(x)=F(r)/r^2= \Delta_r/r^4.$

Let $\G^j,$ $j=0,1,2,3,$ be any $4$ by $4$ Dirac matrices satisfying $\G^i\G^j+\G^j\G^i=2\delta_{ij}I_4,$ and $$\G^5=\G^0\G^1\G^2\G^3.$$ Then $\G^5$ anti-commutes with   $\G^j,$ $j=0,1,2,3.$
Let \[\lb{3.6bis} J=I_4+\alpha (x,\vartheta)\G^3,\qq \alpha(x,\vartheta)=\frac{\sqrt{\Delta_r}}{\sqrt{\Delta_\vartheta}}\frac{a\sin\vartheta}{r^2+a^2}=\ga(x)b(\vt),\qq b(\vt)=\frac{a\sin\vt}{\sqrt{\Delta_\vartheta}},\] be the matrix-valued multiplication operator.
Here $\sup_\vt |\alpha(x,\vt) |$ is exponentially decreasing at both horizons $x\rightarrow\pm\infty.$ By (4.16) in \cite{BelgiornoCacciatori2009} we know that
$\sup_{x,\vartheta} \alpha (x,\vartheta)<1.$ As $\|\G^3\|_\infty=1$ we get that
operator $J$ is invertible and 
\[\lb{3.7} J^{-1}=(1-\alpha^2)(I_4-\alpha\G^3)=\left(1-\frac{\Delta_r}{\Delta_\vt}\frac{a^2\sin^2\vt}{(r^2+a^2)^2}\right)^{-1} \left(I_4-\frac{\sqrt{\Delta_r}}{\sqrt{\Delta_\vt}}\frac{a\sin\vt}{r^2+a^2}\G^3\right).\] 

We consider the  charged Dirac fields with mass $\gm$ represented by 4-spinors belonging to
$$ L^2\left(\R\times\S^2,\frac{\sin\vt}{\sqrt{\Delta_\vt}}dxd\vt d\vp;\,\C^4 \right)$$ and satisfying the evolution equation $i\partial_t\phi=H\phi,\qq H=J^{-1}H_0,$ where
\begin{align*}H_0=&\G^1D_x+\gb(x)\G^0+\ga(x)\left[\sqrt{\Delta_\vt}\left(\G^2\left(D_\vt-i\frac{\cot\vt}{2}\right)+\G^3\frac{E}{\Delta_\vt\sin\vt}D_\vp\right)
-a  \gm\cos\vt \G^5\right]\\&\hspace{3cm}+c(x,D_\vp).\end{align*} Renormalizing spinors $\psi=\left(\frac{\sin\vt}{\sqrt{\Delta_\vt}}\right)^\frac12\phi,$  the new spinor $\psi$ belongs to the Hilbert space independent of parameters of the black hole
$\cH=L^2\left(\R\times\S^2,dxd\vt d\vp;\,\C^4 \right)$ and satisfies the evolution equation
$$i\partial_t\psi=\cD\psi,\qq \cD=J^{-1}\cD_0,\qq \cD_0=\G^1D_x+\gb(x)\G^0+c(x,D_\varphi)+\ga(x)\cD_{\S^2}.$$

Here $\cD_{\S^2}$ is an angular Dirac operator on 2-sphere $\S^2,$
\[\lb{3.10mass}\cD_{\S^2}=\sqrt{\D_\vt}\left[\G^2\left(D_\vt+\frac{i\L a^2\sin(2\vt)}{12\Delta_\vt}\right)+\G^3\frac{E}{\Delta_\vt\sin\vt}D_\vp\right]-a\gm\cos\vt \G^5.\] 
Note that $\cD_0$ is self-adjoint on $\cH,$ while $\cD$ is self-adjoint on slightly modified Hilbert space $\cG$ given by the same space as $\cH$ but equipped with scalar product $\langle.,.\rangle_\cG=(.,J.)_\cH.$
Note that in the expression of the hamiltonian $\cD$  the mass of the Dirac field $\gm$ appears in the coefficient $\gb(x)$ and in the last term in (\ref{3.10mass}). As $\gb(x)$ decays exponentially as $x\rightarrow\pm\infty$ (see (\ref{3.15mass})) it follows that the spectral properties of $\cD$ are independent of the mass $\gm,$ namely  
the  spectrum of $\cD$  is purely absolutely continuous and is given by $\R$ 
(which is proved as in the massless case using Mourre theory, see  \cite{DaudeNicoleau2015}).


We show below that the cut-off resolvent  
$$\cR_\chi(\l)=\chi(\cD-\l)^{-1}\chi,\qq \chi\in C_0^\infty(\R;\C^4)$$ has meromorphic continuation from the upper half-plane to $\C$ with isolated poles of finite rank.
The poles of this meromorphic continuation are called {\em resonances}.

\subsection{Quasi-normal modes.} Resonances for the black holes are  complex characteristic frequencies of the proper solutions of the perturbation equations  which satisfy the boundary conditions appropriate for purely ingoing waves at the event horizon and purely outgoing waves at infinity  (see \cite{ChandrasekharDetweller1975}, \cite{Chandrasekhar1983}).

 These solutions are called 
{\em quasi-normal modes} (QNM). QNMs determine the
late-time evolution of fields in the black hole exterior and
 eventually dominate
the black hole response to any kind of perturbation, therefore containing informations about proper parameters of the black hole: the mass, the electric charge, the angular momentum
per unit mass and De Sitter constant $\L.$ For the physics review we refer to \cite{KokkotasSchmidt1999} and more recent \cite{Bertietal2009}.

The subject has become very popular for the last few decades including the development of   stringent mathematical theory of QNMs (see  \cite{BachelotMotetBachelot1993}, \cite{SaBarretoZworski1997}, \cite{DafermosRodnianski2007}, \cite{Dyatlov2011}, \cite{Dyatlov2012}, \cite{Gannot2014}).

The paper  \cite{SaBarretoZworski1997} provides with  mathematical justification for  localization of QNMs for the wave equation on the
de Sitter-Schwarzschild metric. 
 In Regge-Wheeler coordinates the problem is reduced to the scattering problem for the  Schr{\"o}dinger equation on the line with exponentially decreasing potential.
In the  Schwarzschild  case (zero cosmological constant, which corresponds to asymptotically flat Universe) the Regge-Wheeler potential is only polynomially decreasing and the method does not work due to the possible accumulation of resonances at the origin. A non-zero cosmological constant is needed in order  to define an analytic continuation of the resolvent in a proper space of distributions. Later these works were  complemented by  \cite{BonyHafner2008}, where the authors considered the local energy decay for the wave equation on the
de Sitter-Schwarzschild metric and proved expansion of the solution  in terms of resonances.

In \cite{Dyatlov2011} and \cite{Dyatlov2012} Dyatlov studied the slowly rotating Kerr-de Sitter black holes. Due to cylindrical instead of spherical symmetry the problem can no longer be simply reduced to a scattering problem on the line. The quasi-normal modes split in a way similar to the Zeeman effect. Dyatlov also extended \cite{BonyHafner2008} to the rotating black holes and showed the exponential decay of local energy of linear waves orthogonal to the zero quasi-normal mode. Note also the paper \cite{DaudeNicoleau2015}, where the authors extended their inverse scattering results from \cite{DaudeNicoleau2011} to the scattering for massless Dirac fields by the (rotating) Kerr-Newman-de Sitter black holes.
In the present paper we extend some techniques from \cite{Dyatlov2011}, \cite{Dyatlov2012} to the framework of Dirac operators for the Kerr-Newman-de Sitter black holes. 

Note that some results in \cite{Dyatlov2011} and \cite{Dyatlov2012}  were extended in \cite{Vasy2013}  to the case of perturbations of Kerr-de Sitter metrics when the separation of variables is not possible. For example resonance free regions and polynomials resolvent bounds were establishes by more general and flexible techniques. Recently, there appeared several new papers on applications and extensions of Vasy's method  (see \cite{Dyatlov2015},   \cite{Dyatlov2016} and references therein). 

 As in the present paper we focus on the precise spectral asymptotics where separation of variables is crucial, we do not consider these methods here.

In \cite{Iantchenko2015}, \cite{Iantchenko2016}  we consider scattering of massless uncharged Dirac fields propagating in the outer region of  de Sitter-Reissner-Nordstr{\"o}m (dS-RN) black hole, which is spherically symmetric charged exact solution of the Einstein-Maxwell equations and is special case of non-rotating ($a=0$) KN-dS black hole.

The resonances are approximated by the lattice associated to the trapped set which is a sphere of partially hyperbolic orbits - {\em photon sphere}. Due to radial symmetry, after separation of variables and a Regge-Wheeler transformation the problem is reduced to a family of one-dimensional Schr\"odinger operators on a line with potentials exponentially decaying at infinity and having unique non-degenerate maxima. 

In \cite{Iantchenko2015}, using a special super-symmetric form of the radial Dirac equation, we show that resonances for dS-RN black holes can be obtained as solutions of one-dimensional Schr\"odinger equations and using the method of semiclassical  Birkhoff normal form (as in \cite{Iantchenko2007}, \cite{Iantchenko2008}) we obtain complete asymptotic 
expansions. Moreover, we get similar results for the wave equation in de Sitter-Schwarzschild  metric thus improving  the results in  \cite{SaBarretoZworski1997}. Our results extend to the Dirac operators on spherically symmetric asymptotically hyperbolic manifolds (see \cite{Daudeetal2013}).

Note that the physicist treated the Dirac resonances exactly as solutions of the Schr{\"o}dinger equation  (see also \cite{Chandrasekhar1983}, \cite{Chandrasekhar1980}). In  \cite{Iantchenko2015} Theorem 1 shows a different  point of view and gives exact relation between Schr{\"o}dinger and Dirac resonances.  Indeed, due to the symmetry of the equation, the set of  non-zero Schr{\"o}dinger resonances consists of two sets interposed:
the set of Dirac resonances and its mirror image  with respect to the imaginary axis.

In \cite{Iantchenko2016}
we give an expansion of the massless Dirac fields   in the outer region of  dS-RN black hole in terms of resonances and  describe the decay of local energy for compactly supported data.
The methods extend to the Dirac operators on spherically symmetric asymptotically hyperbolic manifolds.
 
In the case of rotating KN-dS black hole as considered in the present paper the problem in no longer spherically symmetric, but still has cylindrical symmetry. Then following  \cite{DaudeNicoleau2015}, \cite{Dyatlov2011} we show that  it is still possible to decompose operator into the angular and radial parts. But contrary to the  de Sitter-Reissner-Nordstr{\"o}m black hole the radial Dirac operator has no longer super-symmetric form and only leading terms in asymptotic expansion of the resonances can be computed explicitly. 


It is believed that, due to intense gravitation near the event and cosmological horizons of the black hole, even if the Dirac fields are massive, they propagate asymptotically as in the massless case. Below  we support this claim by showing that the two leading terms in the expansion of the black hole resonances are in fact independent of the mass of the Dirac field.

In \cite{Gannot2014} and \cite{Warnick2015} the quasi-normal modes in rather different geometry of Anti-de-Sitter (AdS) black holes are discussed. Such black holes arise in superstring theory via AdS conformal field theory correspondence, that string theory in AdS space is equivalent to conformal field theory in one less dimension (see \cite{Bertietal2009},  \cite{Warnick2015}). The quasi-normal frequencies correspond to the thermalization time scale, which is very hard to compute directly.   Gannot in \cite{Gannot2014} uses a black-box approach to define the quasi-normal modes
after separation of variables and furthermore finds a sequence of quasi-normal frequencies
approaching the real axis exponentially rapidly.   Warnick in \cite{Warnick2015} uses a different approach which applies to asymptotically Anti-de-Sitter black holes and does  not require any separability of the equations
under consideration, nor any real analyticity of the metric. Moreover, the method  can be extended to asymptotically de-Sitter black holes, where it is closely related to approach by Vasy in \cite{Vasy2013}, and  allows consideration
of perturbations which do not vanish on the horizons.

\section{Main results.}

\begin{theorem}\lb{Th1} For $M, Q,\L$ fixed satisfying (\ref{2.5-2.6}), there exists a constant $a_0>0$ such that for $|a|<a_0$ and each $\nu_0,$ there exist constants $C_\l,$ $C_m$ such that the set of quasi-normal modes $\l$ satisfying 
$$\Re\l >C_\l,\qq \Im\l>-\nu_0 $$ coincides modulo ${\mathcal O}(|\l |^{-\infty})$ with the set of pseudopoles 
$$\l=\cF(m,l,k),\qq m\in\Z,\qq l,k\in\Z+\frac12,\qq 0\leq m\leq C_m,\qq |k|\leq l. $$ Here $\cF$ is a complex valued classical symbol of order $1$ in the $(l,k)$ variable, defined and smooth in the cone $\{(m,l,k);\,\,m\in [0,C_m],\,\,|k|\leq l\}\subset\R^3.$ The principal symbol $\cF_0$ of $\cF$ is real-valued and independent of index $m.$  The two leading terms $\cF_0,$ $\cF_1$ in the expansion of  resonances  are  independent of  mass $\gm$ of the Dirac field. Moreover,\\ for $a=0$
\[\lb{leading}\cF=z_0(l+1/2)-i\left(\frac{\alpha}{ z_0}\right)  \left(m+\frac{1}{2}\right)+{\mathcal O}((l+1/2)^{-1}),\]
where \[\lb{z_0} z_0=\left(\frac{M}{r_0^3}-\frac{Q^2}{r_0^4}-\frac{\L}{3}\right)^\frac12,\qq\alpha=\left(\frac{3M}{r_0}-\frac{4Q^2}{r_0^2}\right)^\frac12 z_0^2(x_0);\]
for $|a|<a_0,$
 \[\lb{main-formula}\partial_{\k}\cF_0(m,\pm k,k)=\frac{a}{r_0^2}\left[\frac{4F(r_0)r_0^2}{8Q^2-6Mr_0}\left(1-r_0\frac12 F^{-1}(r_0)F'(r_0)\right)-1\right] -   az_0r_0^{-1}F^{\frac12}(r_0) ,\] where
 $F(r)=1-\frac{2M}{r}+\frac{Q^2}{r^2}-\frac{\Lambda}{3}r^2.$

\end{theorem}

\begin{remark} Formula (\ref{main-formula}) shows the Zeeman-type splitting of the resonances due to rotation. 
\end{remark}

\begin{remark}\lb{R-mass} The fact that the two leading terms in the expansion of  resonances  are  independent of the mass of the Dirac field supports   the well-known  claim that  the massive Dirac fields  propagate asymptotically as in the massless case which is caused by intense gravitation near the event and cosmological horizons of the black hole.

\end{remark}

\begin{remark} Putting $Q=0,$ $r_0=3M,$ we recover the formulas (0.3), (0.4) by Dyatlov in the Kerr-de Sitter case \cite{Dyatlov2012}. Note that in that case $1-r_0\frac12 F^{-1}(r_0)F'(r_0)=0.$ 
\end{remark}

\begin{remark} Note that in the case of chargeless and massless Dirac fields propagating in the exterior of the  static black holes ($a=0$)   formula (\ref{leading})  was obtained in \cite{Iantchenko2015} using a special (super-symmetric) form of the Dirac operator. Moreover, using the method of Birkhoff normal forms, the next order term was obtained
\begin{align*}\cF=&z_0(l+1/2)-i\left(\frac{\alpha}{ z_0}\right)  \left(m+\frac{1}{2}\right)\\
&-\left(\frac{\alpha}{ z_0}\right) \frac{(m+1/2)}{(l+1/2)}\left[-\frac{\alpha}{4z_0^2}(2m+1)+\frac{1}{2} b_{0,2}(2m+1)+ib_{1,2} \right]+{\mathcal O}((l+1/2)^{-2}). \nonumber\end{align*}
\end{remark}

Both massless and massive cases can be treated by similar methods, but as  the massless Dirac fields can be represented by 2-spinors, the massless case is slightly less technical and the methods are more transparent.  Therefore, the strategy of the proof  of Theorem \ref{Th1} will be the following. We explain all the techniques and give the proof of Theorem \ref{Th1} in full details only in the
massless case and  in Section \ref{s-Prelmass} 
we 
indicate how to extend the proofs to the massive Dirac fields.

 In the present paper we do not consider the local energy decay or resonance expansion as studied in the rotationless case in \cite{Iantchenko2015b} and other related properties, leaving these tasks to our forthcoming publications.

{\bf The paper is organized as follows.}\\
We start by  the massless Dirac fields.
 In Section \ref{s-Prel} we explain the decomposition of the Hilbert space using the cylindrical symmetry of the problem and consider the action of  the Dirac operator on the associated subspaces.\\
In Section \ref{s-resolvent} we formulate the main results on meromorphic continuation of the resolvent. We show that the Dirac operator can be written as a tensor product of its angular and radial parts. The resolvent is then represented as a certain contour integral of a tensor product  of angular and radial resolvents.\\
  Properties of angular and radial resolvents are stated in Section  \ref{s-angradres}.  Subsection  \ref{s-red-to-Schr} contains  resolvent identities connecting the Dirac (angular and radial) operators to the diagonal matrix Schr{\"o}dinger-type operators, allowing to transfer some  results from \cite{Dyatlov2011}, \cite{Dyatlov2012} to the Dirac operators.  The proof for the angular resolvent is given in subsection \ref{ss-proofP3.1}.\\  The proof for the radial resolvent      is much longer and is given in Section \ref{s-proofP3.2}.\\ 
In Section \ref{s-Sem-Ref} we reformulate Theorem \ref{Th1} in semi-classical terms.\\ In Sections \ref{s-angular qc}, \ref{s-radial qc} we deduce angular, respectively radial quantization conditions.\\ In Section \ref{s-together} we combine both conditions, deduce formulas (\ref{leading}) and (\ref{main-formula}) and conclude the proof of Theorem (\ref{Th1}) in the massless case.\\
Now, we pass to the general massive case. In  Sections \ref{s-Prelmass} -- \ref{s-normal formmass} we 
indicate how to extend the techniques to the massive Dirac fields and show that the results from previous sections  extend to the massive case.

\newpage
\no{\LARGE {\em Massless Dirac fields.}}

\section{Preliminaries.}\lb{s-Prel}
\subsection{Evolution equation.}
We introduce the matrix-valued multiplication operator \[\lb{3.6} J=I_2-\alpha (x,\theta)\s_2,\qq \alpha(x,\vartheta)=\frac{\sqrt{\Delta_r}}{\sqrt{\Delta_\vartheta}}\frac{a\sin\vartheta}{r^2+a^2}=\ga(x)b(\vt),\qq b(\vt)=\frac{a\sin\vt}{\sqrt{\Delta_\vartheta}},\] with
$\sup_\vt |\alpha(x,\vt) |$ is exponentially decreasing at both horizons $x\rightarrow\pm\infty.$ Here  $\sigma_j,$ $i=1,2,3,$ are the Pauli matrices
\[\lb{Paulim}\sigma_1=\ma 0 &1\\1&0 \am,\qq \sigma_2=\ma 0 &-i\\i&0 \am,\qq\s_3=\ma 1 & 0\\0&-1 \am.\]
Note that in  \cite{DaudeNicoleau2015} the authors used inotation $\G^1=\s_3,$ $\G^2=\s_1,$ $\G^3=-\s_2.$ 

Operator $J$ is invertible and 
\[\lb{3.7} J^{-1}=(1-\alpha^2)(I_2+\alpha\s_2)=\left(1-\frac{\Delta_r}{\Delta_\vt}\frac{a^2\sin^2\vt}{(r^2+a^2)^2}\right)^{-1} \left(I_2+\frac{\sqrt{\Delta_r}}{\sqrt{\Delta_\vt}}\frac{a\sin\vt}{r^2+a^2}\s_2\right).\]

The massless charged Dirac fields are represented by 2-spinors belonging to
$$ L^2\left(\R\times\S^2,\frac{\sin\vt}{\sqrt{\Delta_\vt}}dxd\vt d\vp;\,\C^2 \right)$$ and satisfying the evolution equation $i\partial_t\phi=H\phi,\qq H=J^{-1}H_0,$
$$H_0=\s_3 D_x+\ga(x)\left[\sqrt{\Delta_\vt}\left(\s_1\left(D_\vt-i\frac{\cot\vt}{2}\right)-\s_2\frac{E}{\Delta_\vt\sin\vt}D_\vp\right)\right]+c(x,D_\vp).$$ Renormalizing spinors $\psi=\left(\frac{\sin\vt}{\sqrt{\Delta_\vt}}\right)^\frac12\phi,$ the new spinor $\psi$ belongs to the Hilbert space independent of parameters of the black hole
$\cH=L^2\left(\R\times\S^2,dxd\vt d\vp;\,\C^2 \right)$ and satisfies the evolution equation
$$i\partial_t\psi=\cD\psi,\qq \cD=J^{-1}\cD_0,\qq \cD_0=\s_3D_x+c(x,D_\varphi)+\ga(x)\cD_{\S^2}.$$ Here $\cD_{\S^2}$ is an angular Dirac operator on 2-sphere $\S^2,$
\[\lb{3.10}\cD_{\S^2}=\sqrt{\D_\vt}\left[\s_1\left(D_\vt+\frac{i\L a^2\sin(2\vt)}{12\Delta_\vt}\right)-\s_2\left(\frac{E}{\Delta_\vt\sin\vt} \right)D_\vp\right].\] 
Note that $\cD_0$ is self-adjoint on $\cH,$ while $\cD$ is self-adjoint on slightly modified Hilbert space $\cG$ given by the same space as $\cH$ but equipped with scalar product $\langle.,.\rangle_\cG=(.,J.)_\cH.$
\vspace{2mm}

\subsection{Separation of variables for real $\l$.} In this section we consider the decomposition of the Hilbert space using the cylindrical symmetry 
of  the operator $$\cD=J^{-1}\cD_0,\qq \cD_0=\s_3D_x+c(x,D_\varphi)+\ga(x)\cD_{\S^2},$$ where $\cD_{\S^2}$ is the angular Dirac operator on 2-sphere $\S^2$ given in (\ref{3.10}).

We start with the following identities 

\[\lb{equi-res}(\cD-\l)\psi=\phi\qq\Leftrightarrow\qq\cD(\l)\psi= J\phi,\qq (\cD-\l)^{-1}=[\cD(\l)]^{-1}J,\]
where
\[\lb{3.26}\cD(\l)=(\s_3D_x+c(x,D_\vp)-\l)+\ga(x) (\cD_{\S^2}+\l b(\vt)\s_2).
\] 

Then the stationary Dirac equation $\cD\psi=\l\psi$ can be re-written as $\cD(\l)\psi=0.$

Let \[\lb{3.27} A_{\S^2}(\l)=\cD_{\S^2}+\l b(\vt)\s_2,\qq \cH_{\S^2}=L^2(\S^2,\,d\vt d\vp;\,\C^2).\]

In this section, following  \cite{DaudeNicoleau2015}, we describe the decomposition of the Hilbert space 
in the case  $\l\in\R$.
The construction indicates the way of approaching the case $\l\in\C.$ 

We decompose $\cH_{\S^2}$  onto the angular modes $\{e^{ik\vp}\}_{ k\in \frac12+\Z}$ that are eigenfunctions for $D_\vp$ with anti-periodic boundary conditions (see \cite{BelgiornoCacciatori2009}). Then
\[\lb{3.12bis}
\cH_{\S^2}=\bigoplus_{ k\in \frac12+\Z}\cH_{\S^2}^k,\qq \cH_{\S^2}^k=L^2( (0,\pi),\,d\vt;\,\C^2).
 \] The reduced subspaces $\cH_{\S^2}^k$ remain invariant under the action of $A_{\S^2}(\l)$
and we denote $A_k(\l)=A_{\S^2}(\l)_{|\cH_{\S^2}^k}.$ We have explicitly  
\[\lb{3.30}A_k(\l)= \sqrt{\D_\vt}\left[\s_1\left(D_\vt+\frac{i\L a^2\sin(2\vt)}{12\Delta_\vt}\right)-\s_2\left(\frac{kE}{\Delta_\vt\sin\vt} -\l\frac{a\sin\vt}{\Delta_\vt}\right)\right].\] For each $k\in 1/2 +\Z,$ operator  $A_k(\l)$ is self-adjoint and has discrete simple spectrum $\sigma(A_k(\l))=\{\mu_{kl}(\l)\}_{l\in\Z^*}$ with associated set of eigenfunctions $\{u_{kl}^\l\}_{l\in\Z^*},$ 
$$A_k(\l) u_{kl}^\l(\vt)=\mu_{k,l}(\l)u_{kl}^\l(\vt).$$ Here $\Z^*=\Z\setminus\{0\}.$ Since $\sigma(A_k(\l))$ is discrete, it has no accumulation point and thus
$$\forall k\in 1/2 +\Z,\qq |\mu_{k,l}(\l)| \rightarrow \infty \qq\mbox{as}\qq l\rightarrow \pm\infty.$$
Eigenvalues $\mu_{k,l}(\l)$ of $A_k(\l)$ are also  the  eigenvalues  of $A_{\S^2}(\l)$ with eigenfunctions\\
$Y_{k,l}(\l):=Y_{k,l}(\l,\vt,\vp)=u_{kl}^\l(\vt) e^{ik\vp}.$


By Theorem 3.3 in \cite{DaudeNicoleau2015} we have the following result \begin{lemma}\lb{L-DNangular}Operator $A_{\S^2}(\l),$ $\l\in\R,$ is self-adjoint on $\cH_{\S^2}$ and has pure point spectrum given by a sequence of eigenvalues
$\{\mu_{kl}(\l)\}\in\R$ with associated eigenfunctions $Y_{kl}(\l)\in\cH_{\S^2}.$
Moreover,\\
(i)  ${\displaystyle \cH_{\S^2}=\bigoplus_{(k,l)\in K} {\rm Span}\,(Y_{kl}(\l)),\qq K=\left(\frac12 +\Z\right)\times \Z^*},$\\ \\
(ii) ${\displaystyle A_{\S^2}(\l)Y_{kl}(\l)=\mu_{kl}(\l)Y_{kl}(\l)},$\\ \\
(iii) ${\displaystyle D_\vp Y_{kl}(\l)=kY_{kl}(\l).}$
\end{lemma}
Now, we use the cylindrical symmetry and decompose the Hilbert space $\cH$ onto the angular modes $\{e^{ik\vp}\}_{ k\in \frac12+\Z},$ 
\[\lb{3.12}
\cH=\bigoplus_{ k\in \frac12+\Z}\cH_k,\qq \cH_k=L^2(\R\times (0,\pi),\,dxd\vt;\,\C^2)=L^2(\R;\C^2)\otimes L^2((0,\pi),\,d\vt;\,\C^2).
 \]
 

We choose half-integers $k$  as we want the anti-periodic conditions in variable $\vp:$ the spinors change the sign after a complete rotation (see \cite{BelgiornoCacciatori2009}). Note that $$\mu_{k,-l}(\l)=-\mu_{kl}(\l),\qq Y_{k,-l}(\l)=\s_3Y_{kl}(\l).$$

Using these results we have the decomposition (see  \cite{DaudeNicoleau2015})
$$\cH=\bigoplus_{(k,l)\in I}\cH_{kl}(\l),\qq I=\left(\frac12+\Z\right)\times\N^*,\qq \cH_{kl}(\l)=L^2(\R;\C^2)\otimes Y_{kl}(\l).$$ We choose $I$ instead of $K$  in order to have subspaces $\cH_{kl}$ remain invariant under the action of $\cD(\l)$  (see Section 3.2 in \cite{DaudeNicoleau2015} for details).

Let $\cD_k(\l):=\cD(\l)_{|\cH_k}=\s_3D_x+c(x,k)-\l+\ga(x)A_k(\l)$ be restriction of $\cD(\l)$ to $\cH_k.$

Radial operator $\s_3D_x+c(x,k)-\l$ lets invariant  $\cH_{kl}$ and its action on $\psi=\psi_{kl}\otimes Y_{kl}(\l)\in \cH_{kl}$ is given by

$$ \left[\s_3D_x+c(x,k)-\l\right]\psi= \left(\left[\s_3D_x+c(x,k)-\l\right]\psi_{kl}\right)\otimes Y_{kl}(\l) .$$

Angular operator $A_{\S^2}(\l)$ lets invariant $\cH_{kl}$ and its action on $\psi=\psi_{kl}\otimes Y_{kl}(\l)\in \cH_{kl}$ is given by
\[\lb{3.39} A_{\S^2}(\l)\psi=(\mu_{kl}(\l)\,\s_1\,\psi_{kl})\otimes Y_{kl}(\l) . \]

Then the restriction of $\cD(\l)$ to $\cH_{kl}(\l)$ is given by
\[\lb{3.40}\cD_{kl}(\l)=\s_3D_x+c(x,k)-\l+\mu_{kl}(\l)\ga(x) \s_1 ,\]
with $$c(x,k)=\frac{aEk+\gq Qr}{r^2+a^2},\qq  \ga(x)=\frac{\sqrt{\Delta_r}}{r^2+a^2}$$ satisfying
\[\lb{3.15_0} \ga(x)=a_\pm e^{\k_\pm x}+{\mathcal O}\left(e^{3\k_\pm x}\right),\qq x\rightarrow\pm\infty,\]
\[\lb{3.17_0}c(x,k)=\Omega_\pm(k)+c_\pm e^{2\k_\pm x}+{\mathcal O}\left(e^{4\k_\pm x}\right),\qq x\rightarrow\pm\infty,\]
$$\Omega_\pm=\frac{aEk+\gq Qr_\pm}{r_\pm^2+a^2}. $$ Here $ \k_+<0,\,\,\k_->0$ are fixed constants (surface gravities at cosmological, event horizons respectively) depending on the parameters of the black hole.

\subsection{Properties of eigenvalues $\mu_{kl}$ of angular operator for real $\l$.}\lb{ss-eigang}
In this section we collect some properties of the angular operator
$$A_{\S^2}(\l)=\cD_{\S^2}+\l \frac{a\sin\vt}{\sqrt{\Delta_\vartheta}}\s_2,$$ where $\cD_{\S^2}$ is given in (\ref{3.10}).

 For $\l\in\R$ operator $A_{\S^2}(\l)$ is self-adjoint  on
$\cH_{\S^2}=L^2(\S^2;\C^2)$ and has positive discrete spectrum $\s(A_{\S^2}(\l))=\{\mu_{kl}(\l)\}_{(k,l)\in I},$ 
ordered in such a way that for  each $k\in\frac12+\Z$ and $l\in\N^*$ it follows  $0<\mu_{kl}(\l)< \mu_{k(l+1)}(\l).$ Let $A_k(\l)=A_{\S^2}(\l)_{|\cH_{\S^2}^k}.$ Below we recall some facts
from \cite{DaudeNicoleau2015}  (see proof of Proposition A.1).
\vspace{0.5cm}

  Let $\zeta=\frac{a^2\L}{3},$ $\xi=a\l,$ and consider operator $A_k(\l)$ as operator-valued function $A_k(\z,\x)$ of  complex parameters $\z,\x.$ 
Put
\[\lb{Ak}A_k(\z,\x)=A(\z) {\Bbb D}_{\S^2}^k+B(\z,\x), \]
with $$A(\z)=\sqrt{1+\z\cos^2\vt},\qq B=i\s_1\frac{\z\sin(2 \vt)}{4\sqrt{1+\z\cos^2\vt}}-\s_2\frac{(\z k-\xi)\sin\vt}{\sqrt{1+\z\cos^2\vt}} .$$

 Operator $$ A_k(0,0)\equiv{\Bbb D}_{\S^2}^k=\s_1D_\vt-\s_2\frac{k}{\sin\vt}$$ is the restriction of the standard Dirac operator on $\S^2$ onto the angular mode $\{e^{ik\vp}\},$ $k\in 1/2+\Z.$
 The domain of $A_k(0,0)$ is given by $$\gD=\{u\in \cH_{\S^2}^k,\,\, u\,\,\mbox{is absolutely continuous, }\,\,{\Bbb D}_{\S^2}^ku\in\cH_{\S^2}^k,\,\, u(\pi)=-u(0)\}.$$ 
The spectrum of   $A_k(0,0)$ is simple discrete given by
\[\lb{A.3}\mu_{k,l}(0,0)=\sgn(l)\left(|k|-\frac12+|l|\right),\qq l\in\Z^*\] and
 $A_k(0,0)$ has compact resolvent.

According to ({\ref{2.5-2.6}) we have $\z\in [0,7-4\sqrt{3}]\subset [0,\frac{1}{13.8}]$ and $\x\in\R$ respectively. Now,
we allow parameters $\z,\x$ to be complex  $(\z,\x)\in B(0,\frac{1}{13})\times S,$ where
$B(0,r)=\{z\in\C;\,\,|z|<r\}$   and  $S$ is a narrow strip containing the real axis.  

The operators $A(\z),$ $B(\z,\x)$ are bounded matrix-valued multiplications operators analytic in the variables $(\z,\x)\in B(0,\frac{1}{13})\times S.$ Since the operator $A(\z)$ is also invertible, the operators domain of $A_k(\z,\x)$ is independent on  $(\z,\x)\in B(0,\frac{1}{13})\times S.$ 

Moreover, since for all $u\in\gD,$ $A_k(\z,\x)u$ is a vector-valued analytic function in $(\z,\x),$ and since for $(\z,\x)\in [0,\frac{1}{13}]\times\R$ is self-adjoint on $\cH_{\S^2}^k=L^2( (0,\pi),\,d\vt;\,\C^2),$ then $A_k(\z,\x)$ forms  a self-adjoint holomorphic family of type (A) in variable
$(\z,\x)\in B(0,\frac{1}{13})\times S$ according to Kato's classification.

Using the analytic perturbation theory by Kato \cite{Kato1966}  it was shown in \cite{DaudeNicoleau2015} that  $A_k(\z,\x)$ has compact resolvent for all  $(\z,\x)\in B(0,\frac{1}{13})\times S.$ 
Moreover, 
 for a fixed $k\in1/2+\Z,$ the eigenvalues $$\mu_{kl}(\z,\x),\qq k\in1/2+\Z,\qq l\in\Z^*,$$ of $A_k(\z,\x)$ are simple and depend holomorphically on $(\z,\x)$ in a complex \neigh{} of $[0,\frac{1}{13}]\times\R.$

It was also shown in \cite{DaudeNicoleau2015}, Proposition A.1, that for all $\l\in\R,$ for all $k\in\frac12 +\Z$ and for all $l\in\N^*,$ there exist constants $C_1$ and $C_2$ independent of $k,l$ such that
\[\lb{mubounds}\left(2-e^\frac{1}{26}\right)\left(|k|-\frac12 +l\right)-C_1|k|-C_2-|a\l|\leq\mu_{kl}(\l)\leq e^\frac{1}{26}\left( |k|-\frac12+l\right)+C_1|k| +C_2+|a\l|.\]

\section{Resolvent.}\lb{s-resolvent}
We use the cylindrical symmetry of the operators $\cD,$  $\cD(\l)$ and decompose the Hilbert space $\cH$ onto the angular modes $\{e^{ik\vp}\}_{ k\in \frac12+\Z}$ as in (\ref{3.12}): 
$$
\cH=\bigoplus_{ k\in \frac12+\Z}\cH_k,\qq \cH_k=L^2(\R\times (0,\pi),\,dxd\vt;\,\C^2)=L^2(\R;\C^2)\otimes L^2((0,\pi),\,d\vt;\,\C^2).
 $$
Note that $\cH_k=\cH\cap D'_k,$ where 
\[\lb{((1.2))} D'_k=\{u\in D';\,\, (D_\vp-k)u=0,\qq k\in \frac12+\Z\}.\] This space can be considered  as a subspace of $D'(\R\times\S^2,dxd\vt d\vp;\,\C^2 )$ or of $D'(\S^2,\,d\vt d\vp;\,\C^2)$ along.

In this section we consider the meromorphic continuation of the resolvent  $(\cD-\l)^{-1}$ from the upper half-plane $\C_+$ to the whole complex plane $\C.$ Using  (\ref{equi-res}), (\ref{3.26}), $(\cD-\l)^{-1}=[\cD(\l)]^{-1}J,$ and (\ref{3.12}) this task is reduced to defining an inverse  $[\cD(\l)]^{-1}$  in some subspace of  $\cH_k.$

We will prove the following results.

\begin{theorem}\lb{Th0}\no 1) Let ${\displaystyle \cD_k(\l):=\cD(\l)_{|\cH_k}=\s_3D_x+c(x,k)-\l+\ga(x)A_k(\l)}$ be restriction of $\cD(\l)$ to $\cH_k.$ Here $A_k(\l)$ is given in (\ref{3.30}). The
operator $R(\l,k)=[\cD(\l)]^{-1}$ admits  meromorphic continuation in $\l$ from $\C_+$  into $\C$
as a meromorphic family of operators
$$R(\l,k):\,\,L^2_{\rm comp}\cap D_k'\mapsto L^2_{\rm loc}\cap  D'_k$$
 with poles of finite rank.\\
\no 2) Fix $\delta >0.$  Put $K_{\delta,r}=(r_-+\delta,r_+-\delta) $ and  $ M_{\delta,r}=K_{\delta,r}\times\S^2.$  Denote by the same letters also their   images under the Regge-Wheeler change of variables $r\,\mapsto\, x$ given by
 (\ref{Regge-W}).
 Let $1_{M_{\delta,r}}$ be the operator of multiplication by the characteristic function of $M_{\delta,r}.$ 
Then there exists $a_0>0$  such that if the rotation speed of the black hole satisfies $|a|<a_0,$ we have the following:\\
\no i) Every fixed compact set can only contain $k-$resonances for a finite number of values of $k.$ Therefore, quasi-normal modes form a discrete subset of $\C.$

\no ii) The operators $1_{M_{\delta,r}} R(\l,k)1_{M_{\delta,r}}$ define a family of operators $\cR(\l)=[\cD(\l )]^{-1}$
$$\cR(\l):\,\,L^2(M_{\delta, r})\mapsto L^2(M_{\delta,r})$$ 
meromorphic in $\l\in\C$ with poles of finite rank. 
\end{theorem}


In order to prove the theorem we follow \cite{Dyatlov2011} and represent $R(\l,k)$ as a certain contour integral of a tensor product of two operators acting in different spaces: the  angular and the radial  resolvents. This procedure replaces the separation of variables which can be performed in the non-rotating case of Reissner-Nordstr{\"o}m black holes as in \cite{Iantchenko2015}.  The properties of  the angular and the radial resolvents are stated later in
Propositions \ref{Prop3.1}, \ref{Prop3.2},  in Section \ref{s-angradres} with the proofs given in Sections   \ref{ss-proofP3.1} and \ref{s-proofP3.2}.

 We  start by representing the Hamiltonian as tensor product $H_k(\l)\otimes I_2 +I_2\otimes A_k(\l)$  acting in $L^2(\R,\, dx;\C^2)\otimes \cH_{\S^2}^k,$ where $\cH_{\S^2}^k=L^2( (0,\pi),\,d\vt;\,\C^2).$ Here  $H_k,$ $A_k$ are  operators specified below. 

In the region (\ref{extreg}) $\ga\neq 0$ and we introduce operator
 $$\tilde{\cD}(\l)=[\ga(x)]^{-1}\cD(\l)=[\ga(x)]^{-1}(\s_3D_x+c(x,D_\vp)-\l)+ (\cD_{\S^2}+\l b(\vt)\s_2).$$ Denote $\tilde{\cD}_k(\l)=\tilde{\cD}(\l)_{|\cH_k}$ its restriction to $\cH_k.$
For real $\l,$ its radial $\tilde{\cD}_k^r(\l)$ and angular   $\tilde{\cD}_k^\vt(\l)$    parts  let invariant  $\cH_{kl}$ and the action  $\tilde{\cD}_k(\l)=\tilde{\cD}_k^r(\l)+\tilde{\cD}_k^\vt(\l)$  on $\psi=\psi_{kl}\otimes Y_{kl}(\l)\in \cH_{kl}$ is given by
$$\tilde{\cD}_k(\l)\psi= \left(\tilde{\cD}_k^r(\l)\psi_{kl}\right)\otimes Y_{kl}(\l)+(\s_1\psi_{kl})\otimes (A_{\S^2}(\l)Y_{kl}(\l)) .$$

Let $\phi_{kl}=\s_1\psi_{kl}\in L^2(\R,\, dx;\C^2)$ and 
 $$H_k(\l):=[\ga(x)]^{-1}(\s_3D_x+c(x,k)-\l)\s_1,\qq  A_k(\l)=A_{\S^2}(\l)_{|\cH_{\S^2}^k}$$  acting in $L^2(\R,\, dx;\C^2),$  $\cH_{\S^2}^k={\rm Span}_{l\in Z^*}\,  (Y_{kl}(\l))$ 
respectively.
Then for any $\l\in\C,$ $$\tilde{\cD}_k(\l)=H_k(\l)\otimes I_2 +I_2\otimes A_k(\l)$$  acts in  $L^2(\R,\, dx;\C^2)\otimes\cH_{\S^2}^k.$ 
Let $R_r(\l,\omega,k)=(H_k(\l)+\o)^{-1}
$ and
$R_\vt(\l,\omega,k)=(A_k(\l)-\o
)^{-1}.$ 


We apply the method from \cite{Dyatlov2011}, Section 3.
Take $k\in 1/2+\Z$ and an arbitrary $\delta>0.$ 

Later on in Section \ref{s-angradres} we summarize the properties of  angular $R_\vt$ and radial $R_r$ resolvents. Here we need some of these results.
 From Propositions \ref{Prop3.1}, \ref{Prop3.2} in Section \ref{s-angradres} it follows that the
  resolvents  $$R_r(\l,\omega,k)=(H_k(\l)+\o)^{-1}:\,\,L^2_{\rm comp}(\R,\,dx;\C^2)\,\,\mapsto\,\, H^1_{\rm loc}(\R;\C^2),$$ $$ R_\vt(\l,\omega,k)=(A_k(\l)-\o
)^{-1}:\,\,L^2( (0,\pi),\,d\vt;\,\C^2)\,\,\mapsto\,\, H^1( (0,\pi),\,d\vt;\,\C^2) $$ are {\em meromorphic families} of operators in the sense of Definition 2.1 in \cite{Dyatlov2011}. In particular,  for a fixed values of $\l,$ these families are meromorphic in $\o$ with poles of {\em finite rank} (see Definition 2.2 in \cite{Dyatlov2011}).

A point $\l$ is called regular if the sets $Z_r(\l)=\{\o\in\C;\,\,(\l,\o)\in Z_r\},$ $Z_\vt(\l)=\{\o\in\C;\,\,(\l,\o)\in Z_\vt\}$ do not intersect. 
Here $Z_r,$ $Z_\vt$ are devisors of $R_r,$ $R_\vt$ respectively (see Definition 2.1 in \cite{Dyatlov2011}).
From part 3), Proposition \ref{Prop3.2}, it follows that the set of all {\em regular points}  is non-empty. Indeed, every $\l\in\R$ is regular.

We fix an angle $\psi\in(0,\pi)$ and a regular point $\l.$ Let  $\gamma$  be an {\em admissible contour} on $\C$ at $\l,$ which mean that (see Definition 2.3 in \cite{Dyatlov2011})\\  i) $\g$ is smooth simple contour  given by the rays $\arg\o=\pm\psi$ outside of some compact subset of $\C,$ \\
\no (ii) $\gamma$ separates $\C$ into two regions, $\G_r$ and $\G_\vt,$ such that sufficiently large positive real numbers lie in $\G_\vt$ and $Z_\alpha(\l)\subset\G_\alpha,$ $\alpha\in\{ r,\vt\}.$

Now, the angle $\psi$ of  admissible contour $\g$ at infinity is  chosen as in part 2), Proposition \ref{Prop3.2}.

 From Propositions \ref{Prop3.1}, \ref{Prop3.2} it follows  the following property which assures that admissible contour exists at every regular point.\\
{\em For any compact $K_\l\subset \Omega\subset\C$ there exist constants $C$ and $R$ such that for $\l\in K_\l$ and $|\o|\geq R$,\\
\no (i) for $|\arg\o|\leq\psi$ and $|\pi-\arg\o|\leq\psi$ we have $(\l,\o)\not\in Z_r$ and $\| R_r(\l,\o,k)\|\leq 1/|\o|;$\\
\no (ii) for $\psi\leq |\arg\o|\leq\pi-\psi,$ we have $(\l,\o)\not\in Z_\vt$ and $\|R_\vt(\l,\o,k)\|\leq 1/|\o|.$} 

 Then we can construct restriction of the resolvent $$\cR(\l)=(\cD-\l)^{-1}=[\cD(\l)]^{-1}J
=[\tilde{\cD}(\l)]^{-1}[\ga(x)]^{-1}J$$ to $\cH_k$ at any regular point $\l$ as a contour integral
\[\lb{(2.1)}
R(\l,k)=\tilde{R}(\l,k)[\ga(x)]^{-1},\qq \tilde{R}(\l,k)=\frac{1}{2\pi i}\int_\gamma R_r(\l,\omega,k)\otimes R_\vt(\l,\omega,k) d\o\] for some admissible $\g.$ The orientation of $\gamma$ is chosen so that $\G_r$ always stays on the left.

\begin{figure}[htbp]
\includegraphics[width =16 cm]{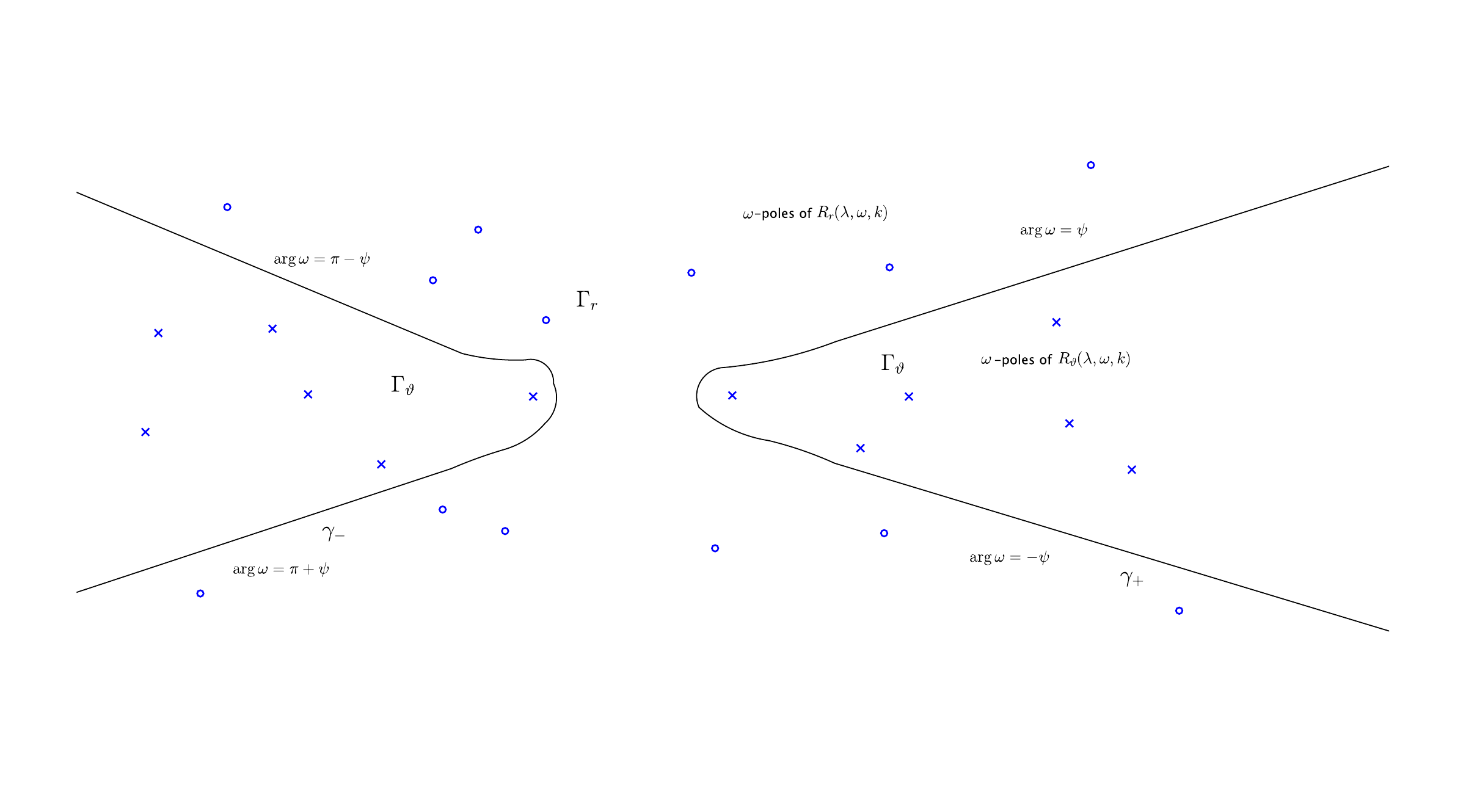}
 \caption{Admissible contour.\label{admcontour}}
\end{figure}

 Integral (\ref{(2.1)}) converges and is independent of the choice of an admissible contour $\g.$ The set of regular points is open and $\tilde{R}(\cdot,k)$ is holomorphic on this set.

From Proposition 2.3 of Dyatlov in \cite{Dyatlov2011} it follows that 
the set of all non-regular points is discrete and the operator $\tilde{R}(\cdot,k)$ given by integral in (\ref{(2.1)}) is defined as an operator on $L^2_{\rm comp}\cap D'_k$ with poles of finite rank. As this can be done for any $\delta >0$ ($\delta$ is the constant from the definition of an admissible contour in part 2), Proposition \ref{Prop3.2}) then $R(\l,k)$ is defined as an operator 
$L^2_{\rm comp}\cap D'_k\mapsto L^2_{\rm loc}\cap  D'_k$ and is meromorphic in $\l$ with poles of finite rank.

Now, we show that $\tilde{\cD}_k(\l)\tilde{R}(\l,k) f=f$ in the sense of distributions for each $f\in L^2_{\rm comp}(\R,\, dx;\C^2) \otimes L^2( (0,\pi),\,d\vt;\,\C^2).$
Let $\l$ be a regular point. Then $\tilde{R}(\l,k)$  is well-defined. By analyticity, we can assume that $\l$ is real which allows us to use the results from Section \ref{s-Prel}. Then $L^2( (0,\pi),\,d\vt;\,\C^2)$ has an orthonormal basis of eigenfunctions $\{u_{kl}^\l\}$ of $A_k(\l) =A_{\S^2}(\l)_{|\cH^k_{\S^2}},$ 
$A_k(\l) u_{kl}^\l(\vt)=\mu_{k,l}(\l)u_{kl}^\l(\vt).$

Let $Y_{kl}(\l)= Y_{kl}(\l, \vt,\vp)=u_{kl}^\l(\vt) e^{ik\vp}.$ 
Put $\o_0=\mu_{k,l}(\l).$ Then \[\lb{res}R_\vt(\l,\omega,k) Y_{kl}(\l)=\frac{ Y_{kl}(\l)}{\o_0-\o}.\]

Let $f,h\in C_0^\infty (\R,\,dx;\C^2),$  $\phi=Y_{kl}(\l)$ as above and $\chi\in C^\infty(\S^2,\,d\vt d\vp;\,\C^2).$

We need to prove
$$I:=\langle \tilde{R}(\l,k)(f(x)\otimes\phi(\vt,\vp), \tilde{\cD}_k(\l) (h(x)\otimes \chi(\vt,\vp)\rangle=\langle f,h\rangle\cdot\langle\phi,\chi\rangle.$$ If   $\g$ is admissible contour,
 then  \begin{align*}I=&\frac{1}{2\pi i}\int_\gamma \langle R_r(\l,\omega,k) f, H_k(\l) h\rangle\cdot \langle R_\vt(\l,\omega,k)\phi,\chi\rangle\\
&+  \langle R_r(\l,\omega,k) f, h\rangle\cdot \langle R_\vt(\l,\omega,k)\phi,A_{\S^2}(\l)\chi\rangle d\o.\end{align*}
Using (\ref{res}) we replace $\gamma$ by a closed bounded contour $\g'$ which contains $\o_0,$ but no poles of $R_r.$ Then
 \begin{align*}I=&\frac{1}{2\pi i}\int_{\gamma'} \langle \left( 1-\o R_r(\l,\omega,k)\right) f,  h\rangle\cdot \langle R_\vt(\l,\omega,k)\phi,\chi\rangle\\
&+  \langle R_r(\l,\omega,k) f, h\rangle\cdot \langle \left( 1+\o R_\vt(\l,\omega,k)\right) \phi,\chi\rangle d\o\\=&\frac{1}{2\pi i}\int_{\gamma'} \langle f,  h\rangle\cdot \langle R_\vt(\l,\omega,k)\phi,\chi\rangle+  \langle R_r(\l,\omega,k) f, h\rangle\cdot \langle\phi,\chi\rangle d\o\\
=&\frac{1}{2\pi i}\int_{\gamma'} \frac{ \langle f,  h\rangle\cdot \langle \phi,\chi\rangle}{\o_0-\o}d\o= \langle f,  h\rangle\cdot \langle \phi,\chi\rangle,\end{align*}
which shows $\tilde{\cD}_k(\l)\tilde{R}(\l,k) f=f$ and finishes the proof of the first part of Theorem \ref{Th0}.

Part 2) of  Theorem \ref{Th0} follows from an analogue of 
Proposition 3.3 in \cite{Dyatlov2011}.\begin{proposition}\lb{Prop_3.3Dy}
Let $\chi\in C_0^\infty(\R;\C^2).$ Then there exists $a_0>0$ (depending on the support of $\chi$) and a constant $C_k$ such that for $|a|<a_0$ and $|k|\geq C_k (1+|\l |),$ $\l$ is not a pole of $R(\cdot,k)$ and we have
\[\lb{(3.12)} \|\chi R(\l,k)\chi\|_{L^2\cap\cD_k'\mapsto L^2}\leq\frac{C_k}{|k|}.\]
\end{proposition}
{\bf Proof.}
This fact  follows from Propositions \ref{Prop3.1}, \ref{Prop3.2}. We choose the constants $\psi,$ $C_r$  from part 2) of Proposition \ref{Prop3.2} and the constants  $C_\vt,$ $C_\psi$  from Proposition \ref{Prop3.1}. Put $\o_0=k/3.$ If $C_k$ is large enough, then
$$|k| >1+C_\vt |a\l |,\qq \o_0 >C_\psi |a\l | +C_r(1+|\l |).$$ We choose the contour $\g$ consisting of the rays $\{\arg\o=\pm\psi,\,\,\pi-\arg\o=\pm\psi,\,\,|\o|\geq\o_0\}$ and the arc $\{ |\o |=\o_0,\,\,|\arg\o|\geq \psi\}.$

By (\ref{(3.2)}) and (\ref{3.4}), all poles of $R_\vt$ lie inside $\g$ (namely, in the region $\{|\o|\geq \o_0,\,\,\psi\leq|\arg\o|\leq \pi-\psi\}$) and \[\lb{3.13}\|R_\vt(\l,\o)\|_{L^2(\S^2)\cap D_k'\mapsto L^2(\S^2)}\leq \frac{C}{|\o|}\] for each $\o$ on $\g.$ Now, suppose that $|a| <a_0=(3C_r)^{-1},$ then (\ref{3.7-3.7}) is satisfied inside $\g.$ 
Then (\ref{(3.12)}) follows from (\ref{(2.1)}): $R(\l,k)=\frac{1}{2\pi i}\int_{\gamma} R_r(\l,\omega,k)\otimes R_\vt(\l,\omega,k) d\o,$ (\ref{3.8-3.8}) and (\ref{3.13}),
as $$\left\|\frac{1}{2\pi i}\int_{\gamma} R_r(\l,\omega,k)\otimes R_\vt(\l,\omega,k) d\o\right\|\leq \frac{C}{|\o_0|}. $$\phantom{1}\hfill\BBox

\begin{remark}\lb{Remark-poles-together} It follows  that a number $\l\in\C$ is a pole of $[\cD(\l)]^{-1}$ if and only if there exist $k\in1/2+\Z,$ $\o\in\C$ such that $(\l,\o,k)$ is a pole for both $ R_r(\l,\o,k)$ and $R_\vt(\l,\omega,k).$\end{remark}

The proof of Theorem \ref{Th0} is now finished.

The following important result comes from the fact that the only trapping in our problem is normally hyperbolic.
\begin{theorem}[Resonance free strip]\lb{Th5-Dy} Fix $\delta >0$ and $s>0.$ Then there exist $a_0>0,$ $\nu_0>0,$ and $C$ such that for $|a| <a_0,$
$$\|\cR(\l)\|_{L^2(M_{\delta, r})\mapsto L^2(M_{\delta,r})}\leq C|\l|^s,\qq |\Re\l|\geq C,\qq |\Im\l|\leq \nu_0.$$
\end{theorem}
The proof of this theorem follows the lines of the proof of Theorem 5 in \cite{Dyatlov2011} and is explained in Section \ref{ss-ResFree}.

\section{Properties of  angular and radial resolvents.}\lb{s-angradres}
In this section we state the main technical results used in the proof of  Theorem \ref{Th0}: Propositions  \ref{Prop3.1} and \ref{Prop3.2}. 
 They are analogue of Propositions 3.1 and 3.2 by Dyatlov in \cite{Dyatlov2011}. Proposition  \ref{Prop3.1} is proved in Subsection \ref{ss-proofP3.1}. Proof of Proposition \ref{Prop3.2} is much longer and will be given later in Section \ref{s-proofP3.2}.

\subsection{Statement of results.}

\begin{proposition}[Angular resolvent estimates]\lb{Prop3.1} There exists a two-sided inverse
$$R_\vt(\l,\o)=(A_{\S^2}(\l)-\o)^{-1}:\,\,L^2(\S^2)\,\, \mapsto\,\,H^1(\S^2),\qq (\l,\o)\in \C^2,$$ with the following properties:\\
\no 1) $R_\vt(\l,\o)$ is meromorphic with poles of finite rank and it has the following meromorphic decomposition near $\Im\l =0,$ $\o=\mu_{k,l}(\l):$
\[\lb{(3.1)}R_\vt(\l,\o)=\frac{S_\vt(\l,\o)}{\o-\mu_{k,l}(\l)}, \] where $S_\vt(\l,\o)$ and $\mu_{k,l}(\l)$ are holomorphic in some $a-$independent  \neigh s of    $\Im\l =0,$ $\o=\mu_{k,l}(\l).$ Moreover, if $\Im\l=0,$ then  $\mu_{k,l}(\l)$ satisfies (\ref{mubounds}).

\no 2) There exists a constant $C_\vt$ such that 
\[\lb{(3.2)}\|R_\vt(\l,\o)\|_{L^2(\S^2)\cap D_k'\mapsto L^2(\S^2)}\leq \frac{C_\vt}{|k|}\qq\mbox{for}\,\ |\o |\leq 0.4 |k|, \,\,|k|\geq C_\vt |a\l|;\] and
\[\lb{(3.3)}\|R_\vt(\l,\o)\|_{L^2(\S^2)\cap D_k'\mapsto L^2(\S^2)}\leq \frac{1}{|\Im\o |}\qq\mbox{for}\,\,|\Im\o| > 2a|\Im\l |.\]

\no 3) For every $\psi >0,$ there exists a constant $C_\psi$ such that
\[\lb{3.4} \|R_\vt(\l,\o)\|_{L^2(\S^2)\cap D_k'\mapsto L^2(\S^2)}\leq \frac{C_\psi}{|\o |}\qq\mbox{for}\,\,\psi\leq |\arg\o|\leq\pi-\psi,\,\,| \o |\geq C_\psi |a\l |.
\]

\end{proposition}

\begin{proposition}[Radial resolvent estimates]\lb{Prop3.2} There exists a family of operators
$$R_r(\l,\o,k)=(H_k(\l)+\o
)^{-1}:\,\,L^2_{\rm comp}(\R,\,dx;\C^2)\,\,\mapsto\,\, H^1_{\rm loc}(\R;\C^2),\qq (\l,\o)\in \C^2$$ with the following properties:

\no 1) For each $k\in 1/2+\Z,$ $R_r(\l,\o,k)$ is meromorphic with poles of finite rank.

\no 2) Take $\delta >0. $ There exists $\psi >0$ and $C_r$ such that \[\lb{3.7-3.7} |\o| \geq C_r,\qq |\arg \o|\leq \psi\qq\mbox{or}\qq |\pi-\arg \o|\leq \psi,\qq |ak|\leq|\o|/C_r,\qq |\l|\leq |\o|/C_r, \] $(\l,\o, k)$ is not a pole of $R_r$ and we have 
\[\lb{3.8-3.8} \|1_{K_{\delta,r}} R_r(\l,\o,k) 1_{K_{\delta,r}}\|_{L^2\mapsto L^2}\leq\frac{C_r}{|\o |}.
\] Here  $1_{K_{\delta,r}}$ is the operator of multiplication by the characteristic function of the Regge-Wheeler image of $(r_-+\delta,r_+-\delta)$ (see  Theorem \ref{Th0}).

\no 3) Resolvent $R_r(\l,\o,k)$ does not have any poles for real $\o$ and real $\l.$ 

\no 4) Assume that $R_r(\l,\o,k)$ has pole at $(\l,\o,k).$ Then there exists a nonzero solution $f\in C^\infty(r_-,r_+)$ to the equation $$(H_k(\l)+\o)f=[\ga(x)]^{-1}\left((\s_3D_x+c(x,k)-\l) \s_1+\o\right)f=0$$ such that the functions
$$|r-r_\pm |^{\mp i\frac{(\l-\Omega_\pm (k))}{2\k_\pm}}f(r)$$ are real analytic at $r_\pm$ respectively.

 \end{proposition}

\subsection{Reduction to the Schr{\"o}dinger-type operators.}\lb{s-red-to-Schr} 
 In this section we derive useful resolvent identities connecting the  angular and radial  Dirac operators to the diagonal matrix Schr{\"o}dinger-type operators similar to those studied in \cite{Dyatlov2011}, \cite{Dyatlov2012}. This allows us in many cases to apply the results obtained in these papers. The results of this section will be used throughout the whole paper and we will formulate them for the semi-classically scaled version of Dirac operators.
\vspace{5mm} 

{\large{\bf Angular operator.}}
We consider an $h$-dependent version of the operator $A_{\S^2}(\l)$ from Section \ref{ss-eigang} 
\begin{align*}& A_{\S^2,h}(\l)=\D^{\frac12}_\vt \left(\s_1\left[hD_\vt+hq_1(\vt)\right]- \s_2\left[q_2(\vt)hD_\vp-\l q_3(\vt)\right]\right)=\ma 0 &D_+\\ D_-&0\am\\ \\
&=\D^{\frac12}_\vt\ma 0 & hD_\vt +hq_1(\vt) +i\left( q_2(\vt)hD_\vp -\l q_3(\vt)\right)\\ hD_\vt +hq_1(\vt) -i\left( q_2(\vt)hD_\vp -\l q_3(\vt)\right)& 0\am
,\end{align*} where
$$q_1(\vt)=\frac{i\L a^2\sin(2\vt)}{12\Delta_\vt},\qq q_2(\vt)=\frac{E}{\Delta_\vt\sin\vt},\qq q_3(\vt)=\frac{a\sin\vt}{\Delta_\vartheta}.$$
Then
\[\lb{sqrang} \left[A_{\S^2,h}(\l)\right]^2= \ma \cP_+&0\\ 0&\cP_-\am ,\] where 
$$\cP_+=D_+D_-=\D_\vt\left\{\left[hD_\vt+hq_1(\vt)\right]^2+\left[q_2(\vt)hD_\vp-\l q_3(\vt)\right]^2-hq_2'(\vt)hD_\vp +\l hq_3'(\vt)\right\}, $$
$$\cP_-=D_-D_+=\D_\vt\left\{\left[hD_\vt+hq_1(\vt)\right]^2+\left[q_2(\vt)hD_\vp-\l q_3(\vt)\right]^2+hq_2'(\vt)hD_\vp -\l hq_3'(\vt)\right\}. $$ 

It is well-known (see \cite{HardtKonstantinovMennicken2000}, \cite{HardtMennicken2001}) that for any $\l\in\C$ the non-zero spectrum of operators $\cP_\pm(\l)$ coincides.  Moreover, 
\[\lb{spectrumA}\s(A_{\S^2,h}(\l))\setminus\{0\}=\{\o\in\C;\,\o^2\in\s(\cP_+(\l))\}\setminus\{0\}=\{\o\in\C;\,\o^2\in\s(\cP_-(\l))\}\setminus\{0\}.
\]
Note that 
\begin{align*}\left[q_2(\vt)hD_\vp-\l q_3(\vt)\right]^2=&\left(\frac{1+\frac{\L a^2}{3}}{\sin\vt \D_\vt}hD_\vp -\l \frac{a\sin\vt}{\Delta_\vartheta}\right)^2\\
=&\frac{\left(1+\frac{\L a^2}{3}\right)^2}{\sin^2\vt \D_\vt^2}\left(hD_\vp -\l a\sin^2\vt\frac{1}{1+\frac{\L a^2}{3}}\right)^2.\end{align*}

The leading term of $\cP_\pm$ is given by
\begin{align*}& \cP_{+,0}(\l)= \cP_{-,0}(\l)=\D_\vt h^2D_\vt^2+\left[q_2(\vt)hD_\vp-\l q_3(\vt)\right]^2\\
&=\Delta_\vt h^2D_\vt^2+\frac{\left(1+\frac{\L a^2}{3}\right)^2}{\sin^2\vt \D_\vt}\left(hD_\vp -\l a\sin^2\vt\frac{1}{1+\frac{\L a^2}{3}}\right)^2 . \end{align*}
We compare it with the leading term of the angular operator in Dyatlov \cite{Dyatlov2012}, Formula (1.11) (after introducing parameter $h$):
$$P_{\vt,0}(\l)= \Delta_\vt h^2D_\vt^2+\frac{\left(1+\frac{\L a^2}{3}\right)^2}{\sin^2\vt \D_\vt}\left(hD_\vp -\l a\sin^2\vt\right)^2,$$
which  is  the same as ours (with  $\breve{\l}=\l/ (1+\frac{\L a^2}{3})=\l/E$ instead of $\l$.)
Let $\o\in\rho(A_{\S^2,h}(\l)).$ Then 
 we have the following  resolvent identity
\[\lb{resDiSch}
(A_{\S^2,h}(\l)-\o)^{-1}=(A_{\S^2}(\l)+\o)\ma (\cP_+-\o^2)^{-1}&0\\
0&(\cP_--\o^2)^{-1}\am 
\] 
or
\[\lb{resDiSch0}(A_{\S^2,h}(\l)-\o)^{-1}=\ma \o(\cP_+-\o^2)^{-1}&D_+(\cP_--\o^2)^{-1}\\
D_-(\cP_+-\o^2)^{-1}&\o(\cP_--\o^2)^{-1}\am, \] where
$ D_\pm= \D^{\frac12}_\vt \left[hD_\vt +hq_1(\vt) \pm i\left( q_2(\vt)hD_\vp -\l q_3(\vt)\right)\right].$

\vspace{5mm}

{\large{\bf Radial operator.}}
 In $\cH=L^2(\R )\os L^2(\R )$ we consider  two Dirac operators $$\cD_{h,\pm}(\l):=-\cD_{h,0}\pm (c(x)-\l)I_2,\qq \cD_{h,0}:=-ih\s_3\partial_x+q(x)\s_1=\ma -ih\partial_x  & q(x)\\q(x) & ih\partial_x\am.$$ Here we only assume  the natural conditions on the real-valued functions $q,c$ so that $\cD_{h,\pm}(\l)$ and all operators below are well-defined and for $\l\in\R$ self-adjoint in $\cH.$  

The product of the operators is given by
\[\lb{produkt}\cD_{h,+}\cD_{h,-}=  \cD_{h,0}^2-(c-\l)^2 I_2-\left[\cD_{h,0},(c-\l) I_2\right].\]Here
$$ \cD_{h,0}^2=-I_2h^2\partial^2_x +\ma q^2&-ihq'(x)\\ihq'(x)&q^2\am,$$  is the  matrix Schr{\"o}dinger operator. The commutator is given by
$$\left[\cD_{h,0},(c-\l) I_2\right]= -ih\s_3 c'(x).$$
The operator $\cD_{h,0}^2$ is self-adjoint in  $\cH=L^2(\R )\os L^2(\R )$ and unitary equivalent to
\[\lb{diagon}U\cD_{h,0}^2U^{-1}=\ma \cP_{h,0}^-&0\\
0&\cP_{h,0}^+\am,\qq \cP_{h,0}^\pm=-h^2\partial_x^2 +q^2\pm h q'.\]
Here, $$U=\frac{i}{\sqrt{2}} \ma 1 & i\\1 & -i\am,\qq U^{-1}=-\frac{i}{\sqrt{2}} \ma 1 & 1\\-i & i\am .$$ We get also
\[\lb{diagon_bis}U\left(\cD_{h,0}^2-(c-\l)^2 I_2\right)U^{-1}=\ma \cP_h^-(\l)&0\\
0&\cP_h^+(\l)\am,\qq \cP_h^\pm(\l)=-h^2\partial_x^2 +q^2-(c-\l)^2\pm h q'.\]
Now, using (\ref{produkt}) we get

\[\lb{produkt-bis} \cD_{h,+}\cD_{h,-}= U^{-1}\ma \cP_{h}^-(\l)&0\\
0&\cP_{h}^+(\l)\am U +ih c'(x)\s_3.
\] If $\cP_h^\pm(\l)$ are invertible, denote $$\cR(\l)=U^{-1}\ma [\cP_{h}^-(\l)]^{-1}&0\\
0&[\cP_{h}^+(\l)]^{-1}\am U$$ and write
$$\cD_{h,+}\cD_{h,-}\cR(\l)=  I +ih c'(x)\s_3\cR(\l).$$
Then we get \begin{align*}&[\cD_{h,+}(\l)]^{-1}\left(I_2+ih\s_3 c'(x)\left(\cD_{h,0}^2-(c-\l)^2 I_2\right)^{-1}\right)=\cD_{h,-}(\l)\left(\cD_{h,0}^2-(c-\l)^2 I_2\right)^{-1}\end{align*}
which leads to
\begin{align}&[\cD_{h,+}(\l)]^{-1}=\cD_{h,-}(\l)U^{-1}\ma [\cP_h^-(\l)]^{-1}&0\\
0&[\cP_h^+(\l)]^{-1}\am U\lb{rad-res-id}\\&\cdot\left(I_2+ih\s_3 c'(x)U^{-1}\ma [\cP_h^-(\l)]^{-1}&0\nonumber\\
0&[\cP_h^+(\l)]^{-1}\am U\right)^{-1}.\end{align}

As $U$ is constant matrix, identity (\ref{rad-res-id}) can be extended to $q=\o\ga$ complex with $\o\in\C.$ This is used in Section \ref{ss-proof2}.

\subsection{Angular resolvent estimates, proof of Proposition \ref{Prop3.1}.}\lb{ss-proofP3.1} 
Here, we prove Proposition \ref{Prop3.1}.

\no 1)  In Section \ref{ss-eigang} it was shown that for $a^2\L,$ $|a\Im  \l |$ small enough $A_k(\l)=A_{\S^2}(\l)_{|\cH_{\S^2}^k}$  has compact resolvent.  Then for each $\o,$ operator $A_{\S^2}(\l)-\o:$ $H^1(\S^2)\mapsto L^2(\S^2)$ is Fredholm and by Proposition 2.2 in \cite{Dyatlov2011}, $R_\vt(\l,\o)=(A_{\S^2}(\l)-\o)^{-1}$ is meropmorphic family of operators $ L^2(\S^2)\mapsto H^1(\S^2).$ 
Decomposition (\ref{(3.1)}) and the meromorphic property of the resolvent for
 $|a\Im  \l |$ small  follow from  Section \ref{ss-eigang}. General case follows from  representation (\ref{resDiSch}) in Section \ref{s-red-to-Schr} (putting $h=1$ there):
$$
R_\vt(\l,\o)=(A_{\S^2}(\l)-\o)^{-1}=(A_{\S^2}(\l)+\o)\ma (\cP_+(\l)-\o^2)^{-1}&0\\
0&(\cP_-(\l)-\o^2)^{-1}\am ,
$$
where $$\cP_\pm(\l)=\D_\vt\left\{\left[D_\vt+q_1(\vt)\right]^2+\left[q_2(\vt)D_\vp-\l q_3(\vt)\right]^2\mp (q_2'(\vt)D_\vp +\l q_3'(\vt))\right\}, $$
 where
$$q_1(\vt)=\frac{i\L a^2\sin(2\vt)}{12\Delta_\vt},\qq q_2(\vt)=\frac{E}{\Delta_\vt\sin\vt},\qq q_3(\vt)=\frac{a\sin\vt}{\Delta_\vartheta}.$$ 
We observe that operators $\cP_\pm(\l)$ are holomorphic families of elliptic second order differential operators on the sphere of the same type as operator $P_\vt(\l)$ in the proof of Proposition  3.1 in \cite{Dyatlov2011}. Therefore, for each $\o,$ we know that the operators $\cP_\pm(\l)-\o^2:\,\,H^2(\S^2)\,\, \mapsto\,\,L^2(\S^2)$ are Fredholm and by Proposition 2.2 in \cite{Dyatlov2011}, $R_\vt(\l,\o)$ is a meromorphic family of operators $L^2(\S^2)\,\, \mapsto\,\,H^1(\S^2).$

\vspace{5mm}

In order to prove  (\ref{(3.2)}) and (\ref{(3.3)}) we need some preliminary inequalities following from Section \ref{ss-eigang}. 
Recall that for each $\l\in\R$ and $k\in1/2+\Z,$ $\mu_{kl}$ are the eigenvalues of the self-adjoint operator $A_k(\l)$ on $\cH_{\S^2}^k=L^2( (0,\pi),\,d\vt;\,\C^2).$

 Let $\zeta=\frac{a^2\L}{3},$ $\xi=a\l,$ and consider operator $A_k(\l)=A_{\S^2}(\l)_{|\cH_{\S^2}^k}$ as operator-valued function $A_k(\z,\x)$ of  complex parameters $\z,\x$ as in (\ref{Ak}):
$ A_k(\z,\x)=A(\z) {\Bbb D}_{\S^2}^k+B(\z,\x)$ 
 with $  A(\z)=\sqrt{1+\z\cos^2\vt},$
\begin{align*}B(\z,\x)=&\frac{1}{4\sqrt{1+\z\cos^2\vt}}\\
&\ma 0 &i(\z\sin(2 \vt)+4(\z k-\xi)\sin\vt)\\i(\z\sin(2 \vt)-4(\z k-\xi)\sin\vt) & 0\am . \end{align*}

 Operator ${\displaystyle A_k(0,0)\equiv{\Bbb D}_{\S^2}^k=\s_1D_\vt-\s_2\frac{k}{\sin\vt}}$ is  restriction of the standard Dirac operator on $\S^2$ onto the angular mode $\{e^{ik\vp}\},$ $k\in 1/2+\Z,$ with domain
 $$\gD=\{u\in \cH_{\S^2}^k,\,\, u\,\,\mbox{is absolutely continuous}\,\,{\Bbb D}_{\S^2}^ku\in\cH_{\S^2}^k,\,\, u(\pi)=-u(0)\},$$ and 
 simple discrete spectrum 
$\mu_{k,l}(0,0)=\sgn(k)\left(|k|-\frac12+|l|\right),$  $l\in\Z^*.$

Now, using \cite{Kato1966}  (Chap. VII, Sect. 3, Theorem 3.6) as in \cite{DaudeNicoleau2015} we get for all $\z\in [0,\frac{1}{13}],$ $k\in\frac12+\Z$ and $l\in\N^*,$
$$|\mu_{kl}(\z,0)-(|k|-\frac12 +l)|\leq\left(e^\frac{1}{26}-1\right)(|k|-\frac12 +l)+2\left(e^\frac{1}{26}-1\right)\left(1+\frac{1}{26}\right)(|k|+\frac14 ),$$
\begin{align*}\mu_{kl}(\z,0)&\geq \left(2-e^\frac{1}{26}\right)(|k|-\frac12 +l)-2\left(e^\frac{1}{26}-1\right)\left(1+\frac{1}{26}\right)(|k|+\frac14 ) \\
&\geq \left(4+\frac{1}{13}-(3+\frac{1}{13})e^\frac{1}{26}\right)\left(|k|+\frac14 \right) >0.8 \left(|k|+\frac14 \right). \end{align*}

For $l<0$ we use $\mu_{k,-l}=-\mu_{kl}$ and get $\mu_{kl}(\z,0) <-0.8 \left(|k|+\frac14 \right).$


Therefore, if $u\in H^1(\S^2)\cap D_k',$ then
$$\| u\|_{L^2}\leq \frac{\| (A(\z){\Bbb D}_{\S^2}^k +B(\z,0) -\o)u\|_{L^2}}{d(\o, \R\setminus (0.8\left(|k|+\frac14\right))(-1,1) )} .$$ 
Now, $$ B(\z,\x)-B(\z,0)=\frac{1}{\sqrt{1+\z\cos^2\vt}}
\ma 0 &-i\xi\sin\vt\\i\xi\sin\vt & 0\am  ,$$ 
 $$(B(\z,\x)-B(\z,0))^2=\frac{\xi^2\sin^2\vt}{1+\z\cos^2\vt} I_2$$ and we get
$$\left\| B(\z,\x)-B(\z,0)\right\|_{L^2(\S^2)\cap D_k'}\leq |\xi|=a|\l |.$$ Therefore,
\[\lb{3.5-3.5} \| u\|_{L^2}\leq \frac{\| A_k(\z,\x) -\o)u\|_{L^2}}{d(\o, \R\setminus (0.8\left(|k|+\frac14\right))(-1,1) )-|a \l|}, \] provided that the denominator is positive. 

\vspace{5mm}

\no 2) Using (\ref{3.5-3.5}) we prove (\ref{(3.2)}). Let $|\o |\leq 0.8\left(|k|+\frac14\right))/2. $ Then $d(\o, \R\setminus (0.8\left(|k|+\frac14\right))(-1,1) )\geq 0.8\left(|k|+\frac14\right))/2$ and
$$d(\o, \R\setminus (0.8\left(|k|+\frac14\right))(-1,1) )-|a \l|\geq 0.4\left(|k|+\frac14\right)- |a \l|\geq 0.2 |k| $$ if
${\displaystyle  |k|\geq |a \l|/0.2.}$

In order to prove (\ref{(3.3)}) we calculate $\Im A_k(\z,\x)=\frac12(A_k(\z,\x)-A_k^*(\z,\x)) $  and get 
\begin{align*}\Im A_k(\z,\x)&=\Im (B(\z,\x)-B(\z,0))=\frac{1}{4\sqrt{1+\z\cos^2\vt}}\\
&\frac12\left\{\ma 0 &-i4\xi \sin\vt\\i 4\xi\sin\vt & 0\am - \ma 0 &-i4\overline{\xi} \sin\vt\\i 4\overline{\xi}\sin\vt) & 0\am\right\}\\
&= \frac{1}{4\sqrt{1+\z\cos^2\vt}} \ma 0 &-i4\Im\xi \sin\vt\\ i 4\Im \xi\sin\vt & 0\am.\end{align*} 
As $$(\Im A_k(\z,\x))^2=  \frac{(\Im \xi)^2\sin^2\vt}{1+\z\cos^2\vt} I_2$$ we get
$$\left\| \Im A_k(\z,\x)\right\|_{L^2(\S^2)\cap D_k'}\leq |\Im\xi|=a|\Im\l |.$$
However, for $u\in H^1(\S^2)\cap D_k',$
\begin{align*}&\|(A_k(\z,\x)-\o)u\|\cdot \|u\|\geq |\Im((A_k(\z,\x)-\o)u,u)|\geq |\Im\o|\cdot \|u\|^2-|(\Im A_k(\z,\x),u,u)|\\
&\geq (|\Im\o|-a|\Im\l |)\| u\|.
\end{align*} If $|\Im\o| >2a |\Im\l |$ we get  $\|(A_k(\z,\x)-\o)u\|\geq a|\Im\l |$ and therefore (\ref{(3.3)}).

\vspace{5mm}

\no 3) We use (\ref{3.5-3.5}).
  If $ \psi\leq |\arg\o|\leq\pi-\psi$  (equivalently $ \psi\leq \arg \o \leq \pi- \psi$ or  $-\pi+\psi\leq \arg \o \leq -\psi$) then $d(\o, \R\setminus (0.8\left(k+\frac14\right))(-1,1) )\geq (|k|+|\o|)/C_2$ with some constant $C_2$ depending on $\psi.$ We get 
$$d(\o, \R\setminus (0.8\left(k+\frac14\right))(-1,1) )-|a \l|\geq  (|k|+|\o|)/C_2 -C_1|a \l|\geq |\o|/C_2 - C_1 |a \l|,$$ which implies (\ref{3.4}).

\section{Radial resolvent estimates, proof of Proposition \ref{Prop3.2}.}\lb{s-proofP3.2} 
Here we prove Proposition \ref{Prop3.2}. 

\subsection{Preliminaries.}
In this section we construct the outgoing (Jost) solutions to the radial Dirac equation.
We consider the radial resolvent
$R_r(\l,\o,k)=(H_k(\l)+\omega)^{-1},$
where
$$H_k(\l):=[\ga(x)]^{-1}(\s_3D_x+c(x,k)-\l)\s_1.$$
Note that $R_r(\l,\o,k)$ satisfies
\begin{align*}&\left([\ga(x)]^{-1}(\s_3D_x+c(x,k)-\l) \s_1+\o\right)R_r(\l,\o,k)f=f,\qq f\in C_0^\infty (\R,\,dx:\C^2),\\
&\Leftrightarrow\qq \left(\underline{\s_3D_x+c(x,k)-\l+\o \ga(x) \s_1}\right) \s_1R_r(\l,\o,k)f=\ga(x)f\\
&\Leftrightarrow\qq R_r(\l,\o,k)=\s_1\left[\s_3D_x+c(x,k)-\l+\o \ga(x) \s_1\right]^{-1} \ga(x).
\end{align*}

Let  $R(\l,\o,k)=[ \cD_k^r(\l,\o)]^{-1}$ be resolvent  of  ${\displaystyle \cD_k^r(\l,\o)=\s_3D_x+c(x,k)-\l+\o \ga(x) \s_1}.$ 

Here $$c(x,k)=\frac{aEk+\gq Qr}{r^2+a^2},\qq  \ga(x)=\frac{\sqrt{\Delta_r}}{r^2+a^2}$$ satisfying
(see  (\ref{3.15_0}), (\ref{3.17_0}))
\[\lb{3.15} \ga(x)=a_\pm e^{\k_\pm x}+{\mathcal O}\left(e^{3\k_\pm x}\right),\qq x\rightarrow\pm\infty,\]
\[\lb{3.17}c(x,k)=\Omega_\pm(k)+c_\pm e^{2\k_\pm x}+{\mathcal O}\left(e^{4\k_\pm x}\right),\qq x\rightarrow\pm\infty,\]
$$\Omega_\pm=\frac{aEk+\gq Qr_\pm}{r_\pm^2+a^2},\qq \O_- >\O_+.$$


Write $$V(x;\l,\omega,k)=c(x,k)-\l+\o \ga(x)\s_1,\qq \cD_k^r(\l,\o):=\s_3D_x+V(x;\l,\omega,k).$$

We consider the Dirac equation  (\ref{DS}) 
\[\lb{DS} \cD_k^r(\l,\o)f=\left(\s_3D_x+V(x;\l,\omega,k)\right)f=0\]
for a vector valued
function $f(x)= f_1(x)e_++f_2(x)e_-,$ where $f_1, f_2$ are  the functions of $x\in\R, $
$$f(x)= f_1(x)e_++f_2(x)e_-,  \ \ \ e_+= \ma 1\\0\am,
   \   e_-= \ma 0\\1\am$$
\[\lb{1.2}
 \left\{\begin{array}{c}
           -if_1'+c(x,k)f_1+\o\ga(x)f_2=\l f_1 \\
           if_2'+c(x,k)f_2+\o\ga(x) f_1=\l f_2
         \end{array}\right. ,\qq (\l,\o) \in \C^2.
\] Note that in \cite{IantchenkoKorotyaev2013} we studied the properties of similar equation with $c(x,k)=0$ and $\ga(x)$ with compact support. Note also  the following property which will be used in  the proof of part 3) of  Proposition \ref{Prop3.2}.
\begin{remark}\lb{Rem-conj}Note that if $f=(f_1(\l,\o),f_2(\l,\o))^T$ is solution of (\ref{DS}) with $(\l,\o)\in\C^2,$ then $\tilde{f}:=(\overline{f}_2(\overline{\l},\overline{\o}),\overline{f}_1(\overline{\l},\overline{\o}))^T$ is also the solution of (\ref{DS})  with the same $(\l,\o).$\end{remark}

 Due to (\ref{3.15}) and (\ref{3.17}) we have
\[\lb{3.17bis}V(x,\l,\omega,k)=\Omega_\pm(k)-\l+\o a_\pm\s_1e^{\k_\pm x}+c_\pm e^{2\k_\pm x}+{\mathcal O}\left(e^{3\k_\pm x}\right),\qq x\rightarrow\pm\infty\]
and the system (\ref{1.2}) in the limit $x\rightarrow\pm\infty$ is given by  $$
 \left\{\begin{array}{c}
           -if_1'+\left(\Omega_\pm(k)+{\mathcal O}\left(e^{2\k_\pm x}\right)\right)f_1+\o\left(a_\pm e^{\k_\pm x}+{\mathcal O}\left(e^{3\k_\pm x}\right)\right)f_2=\l f_1 \\
           if_2'+\left(\Omega_\pm(k)+{\mathcal O}\left(e^{2\k_\pm x}\right)\right)f_2+\o\left(a_\pm e^{\k_\pm x}+{\mathcal O}\left(e^{3\k_\pm x}\right)\right) f_1=\l f_2
         \end{array}\right. ,\qq (\l,\o) \in \C^2.
$$

 The Regge-Wheeler coordinate $x$ is given by $$x=\int_{r_0}^r\frac{s^2+a^2}{\Delta_r(s)}ds,$$ where $r_0\in (r_-,r_+)$ is a fixed number or explicitly
$$x=\frac{1}{2\k_-}\ln(r-r_-)+\frac{1}{2\k_+}\ln(r_--r)+\frac{1}{2\k_c}\ln(r-r_c)+\frac{1}{2\k_n}\ln(r-r_n)+C\in \R $$ with an integration constant $C.$ Here 
$$\k_\sigma=\frac{\Delta_r'(r_\sigma)}{2(r_\sigma^2+a^2)},\qq \sigma=-,+,c,n, $$ and
$\k_->0,$ $\k_+<0.$

We recall the following result by Dyatlov  (\cite{Dyatlov2012}, Proposition 4.1) 
\begin{proposition}[Dyatlov]\lb{Prop4.1Dy} There exists a constant $X_0$ such that for $\pm x>X_0,$ it follows $r=r_\pm\mp F_\pm (e^{ \k_\pm x}),$ where $F_\pm(w)$ are real analytic on $[0, e^{\pm\k_\pm X_0})$ and holomorphic in the discs $\{|w| <e^{\pm \k_\pm X_0}\}\subset \C.$
\end{proposition}
{\bf Proof.} For $x(r)$ near $r=r_+$ we have $2\kappa_+ x(r)=\ln (r_+-r)+G(r),$ where $G$ is holomorphic near $r=r_+.$ Then $$w^2:=e^{2\kappa_+ x}=(r_+-r) e^{G(r)}. $$ We apply  the inverse function theorem to solve for $r$ as a function of $w$ near zero.  

Together with the similar analysis near
 $r=r_-$ it implies that there exists a constant $X_0>0$ such that for $\pm x >X_0,$ we have $r=r_\pm \mp F_\pm (e^{ \kappa_\pm x}),$ where $F_\pm (w)$ are real analytic on $[0, e^{\pm \kappa_\pm X_0})$
and holomorphic in the discs $\{ |w| < e^{\pm \kappa_\pm X_0}\}\in\C.$ \hfill\BBox

We get by this proposition that \[\lb{4.4}V(x,\l,\omega,k)=V_\pm(e^{\k_\pm x}),\qq \pm x >X_0,\] where $V_\pm(w)$ are functions holomorphic in the discs $\{|w| <e^{\pm \k_\pm X_0}\},$ and $V_\pm(e^{\k_\pm x})\rightarrow\Omega_\pm(k)-\l$ as $x\rightarrow\pm\infty.$

Using  (\ref{3.17bis}) we can consider Hamiltonian $\s_3D_x+V(x,\l,\omega,k)$ as asymptotic perturbation for $x\rightarrow\pm\infty$ of the ``free'' operators $\s_3D_x+\Omega_\pm(k) -\l.$ Then the ``free'' Jost solutions $ f^{0,\pm}$ such that $[\s_3D_x+\Omega_\pm(k) -\l]f^{0,\pm}=0$ are given by (compare with \cite{IantchenkoKorotyaev2013})
 $$f^{0,-}(x;\l,k)= \ma 0\\e^{-i(\l-\Omega_-(k)) x}\am,\qq  f^{0,+}(x;\l,k)=  \ma e^{i(\l-\Omega_+(k)) x}\\ 0\am .$$

\begin{definition}[Outgoing solution]\lb{DefOutg} Let $\l$  and $(k,l)\in I$ be fixed constants.  A function $f(x):=f(x;\l,\o,k)$ is called outgoing at $-\infty$ if and only if  $$f(x)=e^{-i(\l-\Omega_-(k))  x} \ma v_1^-(e^{\k_- x};\l,\o,k)\\v_2^-(e^{\k_- x};\l,\o,k)\am
$$
and $v_1^-(\bullet;\l,\o,k),$ $v_2^-(\bullet;\l,\o,k)$ are holomorphic in a \neigh{} of zero. 

 A function $f(x):=f(x;\l,\o,k)$ is called outgoing at $+\infty$ if and only if $$ f(x)=  e^{i(\l-\Omega_+(k)) x} \ma  v_1^+(e^{\k_+ x};\l,\o,k)\\  v_2^+(e^{\k_+ x};\l,\o,k)\am
$$
and $v_1^+(\bullet;\l,\o,k),$ $v_2^+(\bullet;\l,\o,k)$ are holomorphic in a \neigh{} of zero.  

If $f$ is outgoing at both infinities, we call it outgoing.
\end{definition}
In variable $w=e^{\k_+ x}$ we get $e^{i(\l-\Omega_+(k)) x} =w^{i\frac{(\l-\Omega_+(k))}{\k_+}}.$ Moreover, as $w^2=(r_+-r)e^{G(r)},$ with function $G$ holomorphic near $r_+,$ we get $$e^{i(\l-\Omega_+(k)) x}=(r_+-r)^{i\frac{(\l-\Omega_+(k))}{2\k_+}}e^{G(r)i\frac{(\l-\Omega_+(k))}{2\k_+}}.$$ Therefore,
for the outgoing at $+\infty$ function we have
$$ f(x)=  (r_+-r)^{i\frac{(\l-\Omega_+(k))}{2\k_+}}u_+(r;\l,\o,k),$$ 
for the outgoing at $-\infty$ function we have
$$ f(x)=  (r-r_-)^{-i\frac{(\l-\Omega_-(k))}{2\k_-}}u_-(r;\l,\o,k),$$ with $\C^2$-valued functions $u_\pm(r;\l,\o,k)$ real analytic at $r=r_\pm$ respectively.

We get
 \begin{proposition}\lb{PropSol} Let $ \cD_k^r(\l,\o)=\s_3D_x+V(x;\l,\o,k),$ $V(x;\l,\o,k)=c(x,k)-\l+\o \ga(x)\s_1$. 

\no 1)  There exist  outgoing at $\pm\infty$ solutions $f^\pm$ to the equation $ \cD_k^r(\l,\o) f^\pm=0$ of the form 
$$f^\pm(x)=e^{\pm i(\l-\Omega_\pm(k))  x} \ma v_1^\pm (e^{\k_\pm x};\l,\o,k)\\v_2^\pm(e^{\k_\pm x};\l,\o,k)\am,$$ where  $v_\pm( w ;\l,\o,k)$ are holomorphic in $\{ |w| < W_\pm(e^{\pm \k_\pm X_0})\}$ and 
\begin{align}& v_1^-( 0;\l,\o,k)=0,\qq v_2^-( 0;\l,\o,k)=\frac{1}{\G\left(1+2\frac{\l-\O_-}{i\k_-}\right)},\label{4.6-}\\
&  v_1^+( 0;\l,\o,k)=\frac{1}{\G\left(1-2\frac{\l-\O_+}{i\k_+}\right)},\qq v_2^+( 0;\l,\o,k)= 0.\label{4.6}\end{align}

These solutions are holomorphic in $(\l,\omega)$ and are unique unless \[\lb{condition}
-2\frac{\l-\O_-}{i\k_-}\qq\mbox{or respectibely}\qq 2\frac{\l-\O_+}{i\k_+}\qq\mbox{is positive integer}
.\]
Equivalently
\[\lb{condition_bis}\l\in \O_++i\frac12 \k_+ \Z_+ \in\C_- \,\,(\mbox{for}\,\,v^+),\qq \l\in \O_- -i\frac12 \k_- \Z_+ \in\C_-\,\,(\mbox{for}\,\,v^-).\]

\no 2) For $\l$ as in (\ref{condition}) the solutions $f^\pm$ can be identically zero.  However, assume that one of the solutions $f_\pm$ is identically zero. Then every solution $f$ to the equation $ \cD_k^r(\l,\o) f=0$ is outgoing at the corresponding infinity.
\end{proposition}
{\bf Proof.} 

\no 1)  We consider only function $f^+,$ as construction of $f^-$ is similar. Let $$f^+(x)=   e^{i(\l-\Omega_+) x} \ma  v_1^+(e^{\k_+ x};\l,\o,k)\\  v_2^+(e^{\k_+ x};\l,\o,k)\am .$$  Recall that $V(x;\l,\o,k)=c(x,k)-\l+\o \ga(x)\s_1.$ We write  $v^+(e^{\k_+ x}):=v^+(e^{\k_+ x};\l,\o,k),$ $v^+=( v^+_1,v^+_2)^T.$ Now, 
 \begin{align*}&(D_x+c(x,k)-\l)e^{i(\l-\Omega_+) x}v_1^+(e^{\k_+ x})+\o \ga(x) v_2^+(e^{\k_+ x})
\\&=e^{i(\l-\Omega_+) x}\left(\l-\Omega_++c(x,k)-\l+D_x\right)v_1^+(e^{\k_+ x})+\o \ga(x) v_2^+(e^{\k_+ x})e^{i(\l-\Omega_+) x},
\end{align*}
and we get the first equation with $w=e^{\k_+ x}$
$$\left(w \k_+ D_w+c(x,k)-\Omega_+\right)v_1^+(w)+\o \ga(x) v_2^+(w)=0. $$

Now, from
 \begin{align*}&(-D_x+c(x,k)-\l)e^{i(\l-\Omega_+) x}v_2^+(e^{\k_+ x})+\o \ga(x) e^{i(\l-\Omega_+) x}v_1^+(e^{\k_+ x})\\&=e^{i(\l-\Omega_+) x}\left(-\l+\Omega_++c(x,k)-\l-D_x\right)v_2^+(e^{\k_+ x})+\o \ga(x) v_1^+(e^{\k_+ x})e^{i(\l-\Omega_+) x}, 
\end{align*} we get
$(-w \k_+ D_w+\Omega_++c(x,k)-2\l)v_2^+(w)+\o \ga(x) v_1^+(w)=0
.$

Then equation $D_{k}^r(\l,\o)f^+=0$ is equivalent to the system

$$ \left\{\begin{array}{c}
\left(w \k_+ D_w+c(x,k)-\Omega_+\right)v_1^+(w)+\o \ga(x) v_2^+(w)=0           
 \\
        (-w \k_+ D_w+\Omega_++c(x,k)-2\l)v_2^+(w)+\o \ga(x) v_1^+(w)=0 .
         \end{array}\right.$$

Let $x >X_0.$ Then there are holomorphic functions $c^+(w),a+(w)$ such that
\[\lb{behatinf} c(x,k)=c^+(e^{\k_+ x}),\qq   \ga(x)=a^+(e^{\k_+ x}),\qq x >X_0.\]
Now, we consider
$$ \left\{\begin{array}{c}
\left(w \k_+ D_w+c^+(w)-\Omega_+\right)v_1^+(w)+\o a^+(w) v_2^+(w)=0           
 \\
        (-w \k_+ D_w+\Omega_++c^+(w)-2\l)v_2^+(w)+\o a^+(w) v_1^+(w)=0 .
         \end{array}\right.$$

We construct the Taylor series of $v^+=(v_1^+,v_2^+)^T:$
$$v_i^+(w)=\sum_{j\geq 0}v_{i,j}w^j,\qq i=1,2.
$$

Now, we omit $+$ in all indexes. Let $$ c=\sum_{j\geq 0} c_jw^j,\qq c_0=\O,\qq a=\sum_{j\geq 1}a_jw^j.$$

At order $w^0,$ $j=0,$ the system 
\[\lb{analsyst0} \left\{\begin{array}{l}
{\displaystyle c_0-\O=0 }          
 \\
{\displaystyle    c_0 v_{2,0}+(\O-2\l) v_{2,0}=0  }.
         \end{array}\right.\]
is  satisfied as $c_0-\O=0$ and if  $(\O_++c_{0,+}-2\l)v_{2,0}^+=0$ which  force $v_{2,0}^+=0$ unless $\l=\O_+.$ We get the system of equations of order $w^j,$ $j\geq 1:$ 

\[\lb{analsyst} \left\{\begin{array}{l}
{\displaystyle -i\k j v_{1,j}+\sum_{1\leq l\leq  j}c_l v_{1,j-l} +\o\sum_{1\leq l< j}  a_l v_{2,j-l}=0 }          
 \\
{\displaystyle      i \k jv_{2,j}+\sum_{0\leq l < j}c_l v_{2,j-l}+(\O-2\l) v_{2,j} +\o\sum_{1\leq l\leq j}  a_l v_{1,j-l}=0  }.
         \end{array}\right.\]

We solve system (\ref{analsyst}) by induction.
For $j=1$ we get

$$ \left\{\begin{array}{c}
{\displaystyle -i\k v_{1,1}+c_1 v_{1,0}=0\qq\Rightarrow\,\, v_{1,1}=\frac{c_1}{i\k} v_{1,0} }          
 \\ \\
{\displaystyle        i \k v_{2,1}+c_0 v_{2,1}+(\O-2\l) v_{2,1} +\o a_1v_{1,0}=0 \qq\Rightarrow\,\,  v_{2,1}=-\frac{\o a_1}{ i \k  +2(\O-\l) }v_{1,0} }.
         \end{array}\right. ,$$ as $c_0=\O.$
{\rm We can choose $v_{1,0}=1$ if we want a Jost solution $\psi^+.$ But then solution will not be holomorphic! We will chose Gamma function with simple zeros precisely at poles of the solution, see below. } In order to get $v_{2,1}$ we need that the denominator is non-zero: $$ i \k +2(\O-\l) \neq 0\qq\Leftrightarrow\qq  \l\neq \O_++i\frac12\k_+\in\C_-  $$ as $\k=\k_+ <0.$
If $\l =\O_++ i\frac12\k_+$ and $\o\neq 0$ then there is no solution satisfying $v_{1,0}=1$ as necessarily $v_{1,0}=0$ and $v_{1,1}=0,$ $v_{2,1}$ is arbitrary.

Now, suppose we know $\{v_{1,l}, v_{2,l}\}_{l\leq j}.$ Then the system of equations of order $w^j$ 
\[\lb{analsyst_j+1} \left\{\begin{array}{c}
{\displaystyle -i (j+1)\k v_{1,j+1}+\sum_{1\leq l\leq  j+1}c_l v_{1,j+1-l} +\o\sum_{1\leq l< j+1}  a_l v_{2,j+1-l}=0 }          
 \\
{\displaystyle      i \k (j+1)v_{2,j+1}+\sum_{0\leq l < j+1}c_l v_{2,j+1-l}+(\O-2\l) v_{2,j+1} +\o\sum_{1\leq l\leq j+1}  a_l v_{1,j+1-l}=0  }.
         \end{array}\right.\] has solutions
$$  v_{1,j+1}=\frac{1}{i (j+1)\k}\left(\sum_{1\leq l\leq j+1}c_l v_{1,j+1-l} +\o\sum_{1\leq l< j+1}  a_l v_{2,j+1-l}\right), $$
$$ v_{2,j+1}=-\frac{1}{i \k (j+1)+2(\O-\l)}\left(\sum_{1\leq l < j+1}c_l v_{2,j+1-l}+\o\sum_{1\leq l\leq j+1}  a_l v_{1,j+1-l}\right). $$ The last denominator is non-zero if 
$\l\neq \O_++i \frac12\k_+ (j+1)\in\C_-,$  i.e. if $2\frac{\l-\O_+}{i\k_+} $ is not positiv integer.

Now, if $\l= \O+i \frac12\k (j+1)$ for some $j=j_0,$ then solution $v_{2,j+1}$ exists and is arbitrary if and only if
$$\sum_{1\leq l < j+1}c_l v_{2,j+1-l}+\o\sum_{1\leq l\leq j+1}  a_l v_{1,j+1-l}=0$$ - which is only possible for discrete values of $\o$ - or if solution is identically zero in these discrete values. Thus if we want a holomorphic solution we chose $$v_{1,0}=\frac{1}{\G\left(1-2\frac{\l-\O_+}{i\k_+}\right)} .$$

The convergence of the series is proved by induction.

\no 2) The same argument as in the proof of Proposition 4.3 in \cite{Dyatlov2011} applies. Assume that $f^+(x;\l_0,\o_0,k_0)\equiv 0$ (similar for $f^-$). Then $\nu=2(\l_0-\O_+)/i\k_+$ has to be a positive integer. As in part 1), we can construct a nonzero solution $f^1$ to the equation $\cD_k^r(\l,\o)g=0$ with  $$f^1(x)=   e^{-i(\l_0-\Omega_+) x} \ma  \hat{v}_1(e^{\k_+ x})\\  \hat{v}_2(e^{\k_+ x})\am $$ and $\hat{v}$ is holomorphic at zero. We can see that  
\begin{align*}f^1(x)= &  e^{i(\l_0-\Omega_+) x}  e^{-2i(\l_0-\Omega_+) x}\hat{v}(e^{\k_+ x})= e^{i(\l_0-\Omega_+) x}  e^{\frac{-2i(\l_0-\Omega_+) x}{2\k_+x}2\k_+x}\hat{v}(e^{\k_+ x})\\=&e^{i(\l_0-\Omega_+) x}w^\nu \hat{v}(e^{\k_+ x})=:e^{i(\l_0-\Omega_+) x}v(e^{\k_+ x}),\end{align*} where $v(w)=w^\nu\hat{v}(w)$ is holomorphic at zero. Therefore, $ f^1$ is outgoing at $+\infty.$ Note that $f^1(x)=o\left(e^{i(\l_0-\Omega_+) x}\right) $ as $x\rightarrow+\infty.$

Now, since $f^+(x;\l_0,\o_0,k_0)\equiv 0,$ we can define
$$f^2(x)=\lim_{\l\rightarrow\l_0}\G\left(1-2\frac{\l-\O_+}{i\k_+}\right)f^+(x;\l,\o_0,k_0) ,$$
which is an outgoing solution to $\cD_k^r(\l,\o)f^2=0$ and satisfies
$$f^2(x)=e^{ i(\l-\Omega_\pm(k))  x} \left(\ma 1\\0\am+o(1)\right),$$ where  $v_\pm( w ;\l,\o,k)$ are holomorphic in $\{ |w| < W_\pm(e^{\pm \k_\pm X_0})\}$ and 
\[\lb{4.6_next} v_2^-( 0;\l,\o,k)=\frac{1}{\G\left(1+2\frac{\l-\O_-}{i\k_-}\right)},\qq  v_1^+( 0;\l,\o,k)=\frac{1}{\G\left(1-2\frac{\l-\O_+}{i\k_+}\right)}.\]    

\hfill \BBox



\subsection{Proof of parts 1), 3) and 4) in Proposition \ref{Prop3.2}. } 

{\bf Proof of 1).} The integral kernel of the  resolvent $$R(\l,\o,k):=[ \cD_k^r(\l,\o)]^{-1},\qq  \cD_k^r(\l,\o)=\s_3D_x+c(x,k)-\l+\o \ga(x) \s_1,$$ is given by
$$ R(x,y;\l,\o,k)=\left\{\begin{array}{lr}
                 \frac{1}{\det(f^+,f^-)}f^+(x,\l)(f^-(y,\l))^T & \mbox{if}\,\, y<x, \\
                 \frac{1}{\det(f^+,f^-)}f^-(x,\l)(f^+(y,\l))^T & \mbox{if}\,\, x<y,
               \end{array}\right. 
$$
where $f^+=\psi^+,$ $f^-=\vp^-$ are the outgoing solutions.

The Wronskian  $W(\l,\o,k):=\det(f^+,f^-)$ is independent of $x.$ We define  the integral operator $S$ by its kernel $S(x,y;\l,\omega,k)=W(\l,\o,k)\cdot R(x,y;\l,\o,k).$

Using Definition \ref{DefOutg} of the outgoing solutions  we have
$$\det(f^+,f^-)=\ma e^{i(\l-\Omega_+(k)) x}  v_1^+(e^{\k_+ x};\l,\o,k) &e^{-i(\l-\Omega_-(k)) x}  v_1^-(e^{\k_- x};\l,\o,k)\\  e^{i(\l-\Omega_+(k)) x}v_2^+(e^{\k_+ x};\l,\o,k) &e^{-i(\l-\Omega_-(k)) x}v_2^-(e^{\k_- x};\l,\o,k)\am $$
$$=e^{i(\O_-(k)-\Omega_+(k)) x} \left(v_1^+(e^{\k_+ x};\l,\o,k)v_2^-(e^{\k_- x};\l,\o,k)-v_1^-(e^{\k_- x};\l,\o,k)v_2^+(e^{\k_+ x};\l,\o,k)\right).$$
Then $ \cD_k^r(\l,\o)R(\l,\o,k)f=f.$ 

Now, let $ R_r(\l,\o,k)=\s_1R(\l,\o,k) \ga(x).$ Then $(H_k(\l)+\o)R_r(\l,\o,k)=I.$
Resolvent $R_r(\l,\o,k)$ is a meromorphic family of operators $L^2_{\rm comp}\mapsto H_{\rm loc}^1.$ The proof of the fact that $R_r$ has poles of finite rank is done by induction as
in \cite{Dyatlov2011}, Section 4. It
follows  by differentiating $l$ times the identity $(H_k(\l)+\o)R_r(\l,\o,k) =I_2$ in $\l.$  

{\bf Proof of 3) } {\bf Resolvent $R(\l,\o,k)$ does not have any poles for real $\o$ and real $\l.$}  Assume that $\l$ and $\o$ are both real and $R$ has a pole at $(\l,\o, k).$ Then by Proposition \ref{PropSol} there exists the corresponding resonant state $f,$  i.e. a  nonzero
solution $f\in C^\infty(\R,\C^2)$  to the equation $ \cD_k^r(\l,\o)f=0$ such that $f$  is outgoing in the sense of Definition \ref{DefOutg}. We know that it has the asymptotics
$$ f^+(x)=e^{ i (\l-\O_+)x}v^+,\qq v^+=\ma C_+\\0\am+{\mathcal O}(e^{\k_+ x}),\qq x\rightarrow +\infty, $$ 
$$ f^-(x)-=e^{ -i (\l-\O_-)x}v^-,\qq v^-=\ma 0\\C_-\am+{\mathcal O}(e^{\k_- x}),\qq x\rightarrow -\infty, $$ for some nonzero constants $C_\pm.$ Since entries of $V$ are real-valued and due to Remark \ref{Rem-conj} both $f(\l,\o, k)$ and $\tilde{f}(\l,\o, k)$ solve Dirac equation (\ref{DS}). Then the Wronkian $W=\det(f,\tilde{f})$ must be constant. However,
$$\det(f,\tilde{f})\rightarrow \pm|C_\pm|^2\qq\mbox{as}\,\, x\rightarrow\pm\infty$$ which leads to the contradiction.\hfill \BBox


{\bf Proof of 4)} follows the lines of the proof of 4) of Proposition 3.2 in \cite{Dyatlov2011}.
If neither $f^\pm$ is identically zero, then the resolvent $R$ has a pole if and only if the functions $f^\pm$ are linearly dependent, or, equivalently, if there exists a nonzero outgoing (at both ends) solution $f$ to the equation $ \cD_k^r(\l,\o)f=0.$ Now, if one of $f^\pm,$ say $f^+,$ is identically zero, then by part 2) of  Proposition  \ref{PropSol}, $u^-$ will be an outgoing solution at both infinities.

\subsection{Proof of 2) in Proposition \ref{Prop3.2}.}\lb{ss-proof2}

Now, we are able  to prove part 2) from Proposition \ref{Prop3.2}.

Take $h>0$ such that $|\Re\o|=h^{-1}$ and let $\mathfrak{s}=\sign\Re\o.$ Put $$ \tilde{\l}=h\l,\qq \tilde{k}=hk,\qq \tilde{\Omega}_\pm=h\Omega_\pm,\qq\tilde{\mu}=h \Im\o . $$
Then (\ref{3.7-3.7}) implies
\[\lb{3.7-3.7bis} |\tilde{\mu}| \leq \epsilon_r,\qq |a\tilde{k}|\leq\epsilon_r,\qq |\tilde{\l}|\leq \epsilon_r, \] where $\epsilon_r>$ and $h$ can be made arbitrary small by choice of $C_r$ and $\psi.$  
Let
$$ \tilde{\cD}_{h,+}(\l):=h\cD_k^r(\l,\o)=\s_3 h D_x+\tilde{V}_h,\qq \tilde{V}_h(x,\tilde{\l},\tilde{\mu},\tilde{k})=(\mathfrak{s}+i\tilde{\mu} )\ga(x)\s_1+\tilde{c}_h(x,\tilde{k})-\tilde{\l},$$
$$\tilde{c}_h(x,\tilde{k})=\frac{aE\tilde{k}}{r^2+a^2}+h\frac{\gq Qr}{r^2+a^2} .$$
Now, we apply the  resolvent identity (\ref{rad-res-id})
\begin{align*}&[\cD_{h,+}(\l)]^{-1}=\cD_{h,-}(\l)U^{-1}\ma [\cP_h^-(\l)]^{-1}&0\lb{rad-res-id}\\
0&[\cP_h^+(\l)]^{-1}\am U\\&\cdot\left(I_2+ih\s_3 c'(x)U^{-1}\ma [\cP_h^-(\l)]^{-1}&0\nonumber\\
0&[\cP_h^+(\l)]^{-1}\am U\right)^{-1}.\end{align*}

Here ${\displaystyle \cP_h^\pm(\l)=(hD_x)^2 +\tilde{W}_h^\pm(x),}$
where
\[\lb{Wh}\tilde{W}_h^\pm(x)=\tilde{W}_h^\pm(x,\tilde{\l},\tilde{\mu},\tilde{k})=(\mathfrak{s}+i\tilde{\mu} )^2\ga^2(x)-(\tilde{c}_h(x,\tilde{k})-\tilde{\l})^2\pm h (\mathfrak{s}+i\tilde{\mu} )\ga'(x).\]
Note that similarly to (4.4) in \cite{Dyatlov2011} and (\ref{4.4}) we get  \[\lb{4.4Dyatlov}\tilde{W}_h^\pm(x)=W_\pm^\pm(e^{\k_\pm x}),\qq \pm x >X_0,\] where $W_\pm^\pm(w)$ are functions holomorphic in the discs $\{|w| <e^{\pm \k_\pm X_0}\},$ and $W_\pm^\pm(e^{\k_\pm x})\rightarrow-(\Omega_\pm(k)-\l)^2$ as $x\rightarrow\pm\infty.$

The operators $ \cP_h^\pm(\l)$ are of the same type as considered by Dyatlov (his operator $P_x$  given by (4.2) in 
\cite{Dyatlov2011} in unscaled version $h=1$).
Note that the notion of the outgoing vector-function from Definition \ref{DefOutg} is consistent with Definition 4.1 of outgoing function in \cite{Dyatlov2011} in the sense that   each component $f_i,$ $i=1,2,$ of the outgoing at $+\infty$ ($-\infty$) vector-function $f=(f_1,f_2)^T$  is outgoing at the same $\infty$ in the sense of  Definition 4.1  in \cite{Dyatlov2011}.

We apply Proposition 4.4 from \cite{Dyatlov2011} to $\cP=(hD_x)^2 +W(x),$ where
$W(x)$ is either $\tilde{W}_h^+(x)$ or $\tilde{W}_h^-(x),$ given in (\ref{Wh}).

\vspace{5mm}

We get that, if $u\in H^2_{\rm loc}(\R;\C)$ is any outgoing function in the sense of  Definition 4.1 in \cite{Dyatlov2011} and if  $v=\cP u$ is supported in $K_r,$ then\\
\no 1) $u$ can be extended holomorphically to the two half-planes $\{\pm\Re z >X_0\} $ and  $\cP_z u=0$ there. Here  $\cP_z=(hD_z)^2+W(z),$  and $W(z)$ is well defined by 
(\ref{4.4Dyatlov}). 

\no 2) If $\gamma$ is a contour in $\C$ given by $\Im z= F(\Re z),$ $x_-\leq \Re z \leq x_+$ with $F(x)=0$ for $|x| \leq X_0,$ then  there is a restriction to $\g$ of the holomorphic extension of $u$ by \[\lb{scu} u_\g(x)=u(x+iF(x)), \qq \cP_\g u_\g=v,\qq \cP_\g=\left(\frac{1}{1+iF'(x)}hD_x\right)^2+W(x+iF(x)).\] 

\no 3) If $\g$ is as above, $x_\pm=\pm\infty,$ and $F'(x)=C=\const$ for large $|x|.$ Then $$u_\g(x)={\mathcal O}\left(e^{\mp\Im ((1+iC)(\l-\O_\pm)x}\right),\qq\mbox{as}\,\,x\rightarrow\pm\infty. $$ Moreover, if  $\Im (1+iC)(\l-\O_\pm) >0,$ then $u_\g\in H^2(\R;\C^2).$


\vspace{5mm}

Let $$\tilde{P}_{h,\gamma}^\pm=\left(\frac{1}{1+iF'(x)}hD_x\right)^2+\tilde{W}_h^\pm(x+iF(x)), $$ where
$$W:=\tilde{W}_h^\pm(x)=\tilde{W}_h^\pm(x,\tilde{\l},\tilde{\mu},\tilde{k})=(\mathfrak{s}+i\tilde{\mu} )^2\ga^2(x)-(\tilde{c}_h(x,\tilde{k})-\tilde{\l})^2\pm h (\mathfrak{s}+i\tilde{\mu} )\ga'(x).$$ 
Let $u$  be an outgoing function in the sense of Definition 4.1 in \cite{Dyatlov2011} and define $u_\gamma(x),$ $x_-\leq x\leq x_+,$ by (\ref{scu}).

 Then $\tilde{P}_{h,\gamma}^\pm u_\gamma= v,$ where $v$ is supported in $K_r,$ Regge-Wheeler   image of $(r_-+\delta_r,r_+-\delta_r).$ Then as in  \cite{Dyatlov2011}, Section 6, it follows that
$$\|u_\gamma\|_{L^2}\leq C\|v\|_{L^2}.$$ Note that though Dyatlov proved it for $\tilde{P}_{0,\gamma}$ his proof also extends to $\tilde{P}_{h,\gamma}^\pm.$  Then we get  
\[\lb{boundScr} \| 1_{K_r} [\tilde{\cP}_{h}^\pm(\l) ]^{-1} 1_{K_r} \|_{L^2\mapsto L^2}\leq C.\]

In order to transfer bound (\ref{boundScr}) over to the Dirac operator we use  the resolvent identity (\ref{rad-res-id}) and two lemmas which we formulate below. We omit the proofs as they repeat the arguments from the proofs of Lemmas 1, 2 in \cite{Iantchenko2015b}.

Let \[\lb{semiclSobolev}H^s=\{ u\in L^2(\R,\C),\,\,\|u\|_{H^s} <\infty\},\,\, \|u\|_{H^s}^2:=\sum_{k=0}^s\int_\R|(h\partial_x)^k u(x)|^2 dx \] be the standard semi-classical Sobolev spaces, and we define
$\cH^s=H^s\os H^s. $

\begin{lemma}\lb{l-interpol} Suppose (\ref{3.7-3.7bis}). Let $\cP_h$ be either $ \tilde{\cP}_{h}^+(\l)$ or $ \tilde{\cP}_{h}^-(\l)$  and   $\chi\in C_0^\infty (\R;\C).$ Then for $j=1,2$
\[\lb{interpol} \|\chi\cP_h^{-1}\chi\|_{\cL\left(H^0,H^{j} \right)}\lesssim \|\chi\cP_h^{-1}\chi\|_{\cL\left(H^0,H^0 \right)}\
\]
\end{lemma}

From the radial resolvent identity (\ref{rad-res-id}) it follows

\begin{lemma}\lb{p-3.1} Suppose (\ref{3.7-3.7bis}). Let $ \tilde{\cD}_{h,+}(\l),$ $ \tilde{\cP}_{h}^\pm$ be as before.   Let $\chi\in C_0^\infty(\R;\C^2)$ and $\chi_i\in C_0^\infty (\R;\C),$ $i=1,2.$ Then
\[\lb{3.1} \|\chi [\tilde{\cD}_{h,+}(\l)]^{-1}\chi\|_{\cL\left(\cH^0,\cH^{0} \right)}\lesssim\left(\|\chi_1[\tilde{\cP}_h^-]^{-1}\chi_1\|_{\cL\left(H^0,H^{1} \right)}+\|\chi_2[\tilde{\cP}^+_h]^{-1}\chi_2\|_{\cL\left(H^0,H^{1} \right)}\right).\]
\end{lemma}

These two lemmas imply 
$$ \| 1_{K_r}[ \tilde{\cD}_{h,+}(\l)]^{-1} 1_{K_r} \|_{L^2\mapsto L^2}\leq C$$ or (\ref{3.8-3.8})
$$ \|1_{K_r} R_r(\l,\o,k) 1_{K_r}\|_{L^2\mapsto L^2}\leq\frac{C_r}{|\o |}.
$$

\newpage

\section{Semi-classical reformulation of Theorem \ref{Th1}.}\lb{s-Sem-Ref}
Now, we pass to the proof of Theorem \ref{Th1}. We need to reformulate the statement in the semi-classical terms.

We start by recalling Definition 1.1 from \cite{Dyatlov2012}.
\begin{definition}\lb{Def1.1} Let $h>0$  and $\cR(\l;h):\,\,\cH_1\mapsto\cH_2$  be a meromorphic family of operators for $\l\in\cU(h)\subset\C,$ with $\cH_j$ Hilbert spaces. Assume $\O(h)\subset \cU(h)$ be open and $\cZ(h)\subset\C$ be a finite subset (the elements of $\cZ(h)$ may have multiplicities). {\bf The poles of $\cR$ in $\O(h)$ are  simple with polynomial resolvent estimate, given modulo ${\mathcal O}(h^\infty)$ by $\cZ(h),$} if for sufficiently small $h$  there exist maps $\cQ$ and $\Pi$ from  $\cZ(h)$ to $\C$ and the algebra of bounded operators $\cH_1\mapsto\cH_2,$  such that:\\ \\
\no (i) for each $\hat{\l}'\in\cZ(h),$ $\hat{\l}=\cQ(\hat{\l}')$ is a pole of $\cR,$ $|\hat{\l}-\hat{\l}'|={\mathcal O}(h^\infty),$ and $\Pi(\hat{\l}')$ is a rank one operator; \\ 
\no(ii) there is a constant $N$ such that $\|\Pi(\hat{\l}')\|_{\cH_1\mapsto\cH_2}={\mathcal O}(h^{-N})$ for each $\hat{\l}'\in\cZ(h)$ and
$$\cR(\l;h)=\sum_{\hat{\l}'\in\cZ(h)}\frac{\Pi(\hat{\l}')}{\l-\cQ(\hat{\l}')} +{\mathcal O}_{\cH_1\mapsto\cH_2}(h^{-N}),\qq\l\in\O(h).$$
\\ \\
So every pole of $\cR$ in $\O(h)$ lie in the image of $\cQ.$
\end{definition}

Theorem \ref{Th1} results from the  following $h-$dependent version (similarly to Proposition 1.2 by Dyatlov \cite{Dyatlov2012}).
\begin{proposition}\lb{Pro1.2Dy} Let $\nu_0 >0$ and  $h>0.$  Then, for sufficiently small $a$ (independently of $\nu_0$), the poles of $R(\l)$ in the region
\[\lb{1.7} |\Im\l|\ <\nu_0,\qq h^{-1}<|\Re\l |<2h^{-1},\]
are simple with  a polynomial resolvent estimate, given modulo ${\mathcal O}(h^\infty)$ by
\[\lb{1.8}\l=h^{-1}\cF^\l(m, hl,hk;h),\qq m\in\Z,\qq l,k\in\Z+\frac12,\qq 0\leq m\leq C_m,\qq C_l^{-1}\leq hl\leq C_l,\qq |k|\leq l. \] Here $C_m$ and $C_l$ are some constants and $$\cF^\l(m,\tilde{l},\tilde{k};h)\sim\sum_{j\geq 0}h^j\cF_j^\l(m,\tilde{l},\tilde{k})$$ is a classical symbol. The principal symbol $\cF_0$ of $\cF$ is real-valued and independent of $m.$ Moreover, {\bf for} $\bf{a=0}$ 
\[\lb{leading_bis}\cF^\l=z_0(\tilde{l}+h/2)-ih\left(\frac{\alpha}{ z_0}\right)  \left(m+\frac{1}{2}\right)+{\mathcal O}(h^2),\]
where $z_0,$ $\alpha$ are given in (\ref{z_0});
{\bf for} $\bf{|a|<a_0},$
 \[\lb{main-formula_bis} \partial_{\tilde{k}}\cF_0(m,\pm \tilde{k},\tilde{k})=\frac{a}{r_0^2}\left[\frac{4F(r_0)r_0^2}{8Q^2-6Mr_0}\left(1-r_0\frac12 F^{-1}(r_0)F'(r_0)\right)-1\right]-az_0r_0^{-1}F^{\frac12}(r_0) ,\] where
 $F(r)=1-\frac{2M}{r}+\frac{Q^2}{r^2}-\frac{\Lambda}{3}r^2.$


\end{proposition}

Proposition \ref{Pro1.2Dy}  is consequence  of the following results.
\begin{proposition}\lb{Prop1.3Dy} Take $\delta>0$ and let
 $M_{\delta,r}$ be the image of $ (r_-+\delta,r_+-\delta)\times \S^2$  under the Regge-Wheeler transformation $r\,\,\mapsto\,\, x$ (see (\ref{Regge-W})).
 Let $1_{M_{\delta,r}}$ be the operator of multiplication by the characteristic function of $M_{\delta,r}.$ Then, for   $a$ small enough and fixed $\nu_0,$ the poles of the cut-off resolvent
$$R_M(\l)=1_{M_{\delta,r}} R(\l)1_{M_{\delta,r}}:\,\,L^2(M_{\delta,r})\mapsto L^2(M_{\delta,r})$$   in the region (\ref{1.7}) are simple with polynomial resolvent estimate $L^2\mapsto L^2,$   given modulo ${\mathcal O}(h^\infty)$ by (\ref{1.8}). 
\end{proposition}

\begin{proposition}\lb{Prop1.4Dy} Let $R_M(\l)$ be as in Proposition \ref{Prop1.3Dy} and denote its restriction $R_M(\l,k)=R_M(\l)_{|\cH_k},$ where 
 $\cH_k=\cH\cap D'_k,$ is the subspace of angular momentum $k\in \frac12+\Z$ as in  (\ref{((1.2))}).

There exists a constant $C_k$ such that for each $k\in\frac12+\Z,$\\
\no 1) if $h|k|>C_k,$ then $R_M(\l,k)$ has no poles in the region (\ref{1.7})  and its $L^2\mapsto L^2$ norm is ${\mathcal O}(|k|^{-1}).$\\ 
\no 2) if  $h|k|\leq C_k,$ then the poles of $R_M(\l,k)$ in the region (\ref{1.7}) are simple with a polynomial resolvent estimate $L^2\mapsto L^2,$   given modulo  ${\mathcal O}(h^\infty)$ by (\ref{1.8})
\end{proposition}

Propositions \ref{Prop1.3Dy}, \ref{Prop1.4Dy}  are analogue for the Dirac case of 
 Propositions 1.3, 1.4 in \cite{Dyatlov2012}.

The construction of the cut-off resolvent in Propositions \ref{Prop1.3Dy}, \ref{Prop1.4Dy} follows from Part 2) of Theorem \ref{Th0}.


Part 1) of Proposition \ref{Prop1.4Dy} is a reformulation of Proposition \ref{Prop_3.3Dy}.


In Sections \ref{s-angular qc}, \ref{s-radial qc}, \ref{s-together} we prove the second part of 
Proposition \ref{Prop1.4Dy}, namely expansion (\ref{1.8}).

\section{Angular quantization condition.} \lb{s-angular qc}

\subsection{Main result.}\lb{ss-angular-scaling} 
For small $h>0$ we put
$$ \tilde{\l}=h\Re\l,\qq\tilde{\nu}=\Im\l,\qq \tilde{k}=hk,\qq \tilde{\Omega}_\pm=h\Omega_\pm,\qq  \tilde{\o}=h\Re\o,\qq\tilde{\mu}=\Im\o.$$ 

The following results are   analogue for the Dirac case of Proposition 1.6 in \cite{Dyatlov2012}.
\begin{proposition} (Angular)\lb{Prop1.6} The poles $\tilde{\o}+ih\tilde{\mu}$ of $R_\vt(\l,\o,k)=(A_k(\l)-\o)^{-1}$ as a function of $\o$ in the region \[\lb{1.15}
1<\tilde{\l}<2,\qq |\tilde{\nu}| <\nu_0,\qq |\tilde{k}| <C_k,\qq C_\vt^{-1}<\tilde{\o}<C_\vt,\qq |\tilde{\mu}| <C_\vt
\] are simple  with polynomial resolvent estimate $L^2\mapsto L^2,$  given modulo ${\mathcal O}(h^\infty)$ by  \[\lb{1.16 }(\tilde{\o}+ih\tilde{\mu})^2=\cF^\vt(hl,\tilde{\l},\tilde{\nu},\tilde{k};h),\qq l+1/2\in\Z,\qq 0\leq m\leq C_m,\qq \max(|\tilde{k}|,C_l^{-1})\leq h| l | \leq C_l, \] (in a sense of Definition \ref{Def1.1}) 
for some constant $C_l.$ The principal part $\cF_0^\vt$ of the classical symbol $\cF^\vt$ is real-valued, independent of $\tilde{\nu},$ and (in the non-rotating case)
$$\cF^\vt=\tilde{l}(\tilde{l}+h)+{\mathcal O}(h^\infty)\qq\mbox{for}\,\, a=0.$$
Moreover, \[\lb{viktig-ang}\cF_0^{\vt}(\pm \tilde{k},\tilde{\l},\tilde{k})=(E\tilde{k}-a\tilde{\l})^2,\qq \partial_{\tilde{l}}\cF_0^{\vt}(\pm \tilde{k},\tilde{\l},\tilde{k})=\pm2\tilde{k}+{\mathcal O}(a^2),\qq E=1+\frac{\L a^2}{3},\] and consequently, $\partial_{\tilde{k}}\cF^\vt_0(\pm \tilde{k},\tilde{\l},\tilde{k})=-2a\tilde{\l}+{\mathcal O}(a^2).$
\end{proposition}

In the following sections \ref{ss-Joint-s}--\ref{ss-Angu} we prove Proposition \ref{Prop1.6}.

\subsection{Joint spectrum and diagonalization.}\lb{ss-Joint-s}
We start by recalling Definition A.1 from \cite{Dyatlov2012}.
 We say that $\o=(\o_1,\ldots,\o_n)\in\C^n$ belongs to the  {\em joint spectrum} of (matrix-valued) \pseudor s $P_1,\ldots,P_n\in \Psi^{k_j},$ $k_j\geq 0,$ on a compact manifold $M,$ if the joint eigenspace $$\{u\in C^\infty(M);\,\,P_ju=\o_j u,\,\,j=1,\ldots n\} $$ is nontrivial. 
The classes $ \Psi^k$ are the natural generalization to the matrix-valued case of the usual classes of \pseudor s  defined in \cite{Dyatlov2012}, Section 2.

Put
 \begin{align*}&P_1(\tilde{\l},\tilde{\nu};h):=hA_{\S^2}(\l)=A_{\S^2,h}(\tilde{\l} +ih\tilde{\nu})\\
&=\sqrt{\D_\vt}\left[\s_1\left(hD_\vt+\frac{ih\L a^2\sin(2\vt)}{12\Delta_\vt}\right)-\s_2\left(\frac{E}{\Delta_\vt\sin\vt} \right)hD_\vp\right]+(\tilde{\l} +ih\tilde{\nu} )\frac{a\sin\vt}{\sqrt{\Delta_\vartheta}}\s_2,\\
&P_2(h):=I_2hD_\vp.\end{align*}
Then, $(h\l,h\o,hk)$ is a pole of  $R_\vt^h(\o)=(hA_k(\l)-h\o
)^{-1}$ 
 if and only if $(\tilde{\o}+ih\tilde{\mu},\tilde{k})$ lies in the joint spectrum of the operators $(P_1,P_2).$ If  $h\Im\l=\tilde{\nu}=0$ then from 
Section \ref{ss-eigang} we know that  the joint spectrum is given by
$$\left( h\mu_{kl}(\l),hk\right),\qq (k,l)\in(1/2+\Z)\times Z^*.$$ 
 In the rotationless case $a=0$
$$ P_1(\tilde{\l},\tilde{\nu};h)_{|a=0}=\sqrt{\D_\vt}\left[\s_1hD_\vt-\s_2\frac{1}{\sin\vt}hD_\vp\right] =h\sqrt{\D_\vt}{\mathbb D}_{\S^2},$$ where ${\mathbb D}_{\S^2}=\s_1D_\vt-\s_2\frac{1}{\sin\vt}D_\vp$ is the usual Dirac operator on $\S^2$ (with our choice of weight on the spinor) and the joint spectrum of $(P_1,P_2)$ is given by the spherical harmonics
$$\left(\pm\,(\tilde{l}+h/2),\tilde{k}\right),\qq \tilde{k}\in h(1/2+\Z),\qq \tilde{l}\in h(1/2+\N_0), |\tilde{k}|\leq |\tilde{l}|. $$

Now, we use the results obtained in Section \ref{s-red-to-Schr}, Formula (\ref{sqrang}):    $$ \left[P_1(\tilde{\l},\tilde{\nu};h)\right]^2=\left[A_{\S^2,h}(\tilde{\l} +ih\tilde{\nu})\right]^2= \ma \cP_+(\tilde{\l} +ih\tilde{\nu})&0\\ 0&\cP_-(\tilde{\l} +ih\tilde{\nu})\am ,$$ where
\begin{align}&\cP_\pm(\l)=\D_\vt\left\{\left[hD_\vt+hq_1(\vt)\right]^2+\left[q_2(\vt)hD_\vp-\l q_3(\vt)\right]^2\mp h(q_2'(\vt)hD_\vp -\l q_3'(\vt))\right\}\lb{orifoperators}\\
&=\Delta_\vt \left[hD_\vt+hq_1(\vt)\right]^2+\frac{ E^2}{\sin^2\vt \D_\vt}\left(hD_\vp -\l a\sin^2\vt\frac{1}{ E}\right)^2\nonumber\\
&\hspace{4cm}\mp h\left(\Delta_\vt  q_2'(\vt)hD_\vp -\l \Delta_\vt q_3'(\vt)\right)\nonumber
, \end{align}  where
$$q_1(\vt)=\frac{i\L a^2\sin(2\vt)}{12\Delta_\vt},\qq q_2(\vt)=\frac{E}{\Delta_\vt\sin\vt},\qq q_3(\vt)=\frac{a\sin\vt}{\Delta_\vartheta}.$$
Note that for any $\l\in\C$ (see \ref{spectrumA}) $$\s(A_{\S^2,h}(\l))\setminus\{0\}=\{\o\in\C;\,\o^2\in\s(\cP_+(\l))\}\setminus\{0\}=\{\o\in\C;\,\o^2\in\s(\cP_-(\l))\}\setminus\{0\}.$$

 The principal symbol of both $ \cP_\pm(\tilde{\l} +ih\tilde{\nu})$ is 
$$p=\D_\vt\xi_\vt^2+\frac{ E^2}{\sin^2\vt\Delta_\vt}\left(\xi_\vp-\tilde{\l} a\sin^2\vt \frac{1}{E}\right)^2,\qq E=1+\frac{\L a^2}{3}. $$

Note that in \cite{Dyatlov2012}, page 1130, Dyatlov considered symbol $$p_{10}(\vt,\xi_\vt,\xi_\vp;\tilde{\l})=\D_\vt\xi_\vt^2+\frac{ E^2}{\sin^2\vt\Delta_\vt}\left(\xi_\vp-\tilde{\l} a\sin^2\vt \right)^2$$ which is analogue to our $p$ and   
 coincides with ours after rescaling $\tilde{\l}=E\breve{\l}.$
Then we follow \cite{Dyatlov2012}. Let $$ \left\{\begin{array}{l}
\mu_{10\pm}(\vt,\xi_\vt,\xi_\vp;\tilde{\l})=\pm\D_\vt^{-\frac12}\sqrt{p_{10}}=\pm\D_\vt^{-\frac12}\left(\D_\vt\xi_\vt^2+\frac{E^2}{\sin^2\vt\Delta_\vt}\left(\frac{1}{\sqrt{\D_\vt}}\xi_\vp-\breve{\l} a\sin^2\vt \right)^2  \right)^\frac12
 \\
   p_{20}=\xi_\vp 
         \end{array}\right.$$ and
 $ {\mathbf p}(\vt,\xi_\vt,\xi_\vp):=(\mu_{10\pm}, p_{20}).$ 
In the non-rotating case $a=0$ we have $$ \left\{\begin{array}{l}
\mu_{10\pm}(\vt,\xi_\vt,\xi_\vp;\tilde{\l})=\pm\left(\xi_\vt^2+\frac{\xi_\vp^2}{\sin^2\vt}  \right)^\frac12
 \\
   p_{20}=\xi_\vp .
         \end{array}\right.$$ So in $a=0$ case we get $|\mu_{10\pm}|\geq \xi_\vp.$ 
 Then the set   $ {\mathbf p}^{-1}(\tilde{\o},\tilde{k})$  is non-empty  if $|\tilde{\o}|\geq\tilde{k}.$ If $a\neq 0$ small, we get $|\mu_{10\pm}|\geq \xi_\vp$ and 
 the set   $ {\mathbf p}^{-1}(\tilde{\o},\tilde{k})$  is non-empty  if $|\tilde{\o}|\geq E(\tilde{k}-\breve{\l}a).$ 

Fix $\xi_\vp\geq \epsilon_k.$ Then
$\mu_{10\pm}(\vt,\xi_\vt,\xi_\vp;\tilde{\l})$ as a function of $(\vt,\xi_\vt)$ for small $a$ 
has unique critical point at $(\vt,\xi_\vt)=(\pi/2,0)$ with critical value $\pm E(\xi_\vp-\breve{\l}a)$ with signature of the hessian $(\pm 1,\pm 1).$ It is enough to verify for $a=0:$
$$
 \left\{\begin{array}{lc}
(\mu_{10\pm})'_\vt=\mp\left(\xi_\vt^2+\frac{\xi_\vp^2}{\sin^2\vt}  \right)^{-\frac12}\sin^{-3}(\vt)\cos\vt=0&\Leftrightarrow\vt=\pi/2
 \\
(\mu_{10\pm})'_{\xi_\vt}=\pm\left(\xi_\vt^2+\frac{\xi_\vp^2}{\sin^2\vt}  \right)^{-\frac12}\xi_\vt=0&\Leftrightarrow\xi_\vt=0.
         \end{array}\right.$$
$$
 \left\{\begin{array}{l}
(\mu_{10\pm})''_{\vt\vt}=\left[\mp\left(\xi_\vt^2+\frac{\xi_\vp^2}{\sin^2\vt}  \right)^{-\frac12}\sin^{-3}(\vt)\right]'_\vt\cos\vt\pm\left(\xi_\vt^2+\frac{\xi_\vp^2}{\sin^2\vt}  \right)^{-\frac12}\sin^{-3}(\vt)\sin\vt|_{(\vt,\xi_\vt)=(\pi/2,0)}\\
\hfill=\pm(\xi_\vp)^{-1}
 \\
(\mu_{10\pm})''_{\xi_\vt\xi_\vt}=\left[\pm\left(\xi_\vt^2+\frac{\xi_\vp^2}{\sin^2\vt}  \right)^{-\frac12}\right]'_{\xi_\vt}\xi_\vt\pm\left(\xi_\vt^2+\frac{\xi_\vp^2}{\sin^2\vt}  \right)^{-\frac12}|_{(\vt,\xi_\vt)=(\pi/2,0)}=\pm(\xi_\vp)^{-1}.
         \end{array}\right.$$

 For $\epsilon>0$ (will be chosen small enough) let  $\tilde{K}_\epsilon=\{ (\tilde{\o},\tilde{k});\,\, \tilde{k}\geq \epsilon,\,\,E^2(\tilde{k}-aE^{-1}\tilde{\l})^2\leq\tilde{\o}^2\leq C_\vt^2 \}\subset\R^2.$ 
Now, we may apply Proposition 2.10 of Dyatlov in \cite{Dyatlov2012} to the function $\mu_{10\pm}(\cdot,\cdot,\xi_\vp)$ and obtain a function $F_+(\tilde{\o};\tilde{k})$ on  $\tilde{K}_\epsilon$ such that 
${F_+}|_{\tilde{K}_\epsilon}=0$ and $\partial_{\tilde{\o}}F_+>0.$ 
If $a=0$ then $F_\pm=\tilde{\o}\mp\tilde k.$

Moreover, by the same Proposition 2.10 of Dyatlov in \cite{Dyatlov2012}  $$F'_+(\mu_{10\pm}(\pi/2,0,\xi_\vp))=(\det\nabla^2\mu_{10\pm}(\pi/2,0,\xi_\vp))^{-1/2}. $$

We have  the following version of  Proposition 3.1 by Dyatlov from \cite{Dyatlov2012} applied to either $ \cP_+(\tilde{\l} +ih\tilde{\nu})$ or $ \cP_-(\tilde{\l} +ih\tilde{\nu})$ defined in (\ref{orifoperators}):  \begin{proposition}\lb{Prop3.1_2012} 
Let $\tilde{K}=\{ (\tilde{\o},\tilde{k});\,\,C_\vt^{-1}\leq\tilde{\o}\leq C_\vt,\,\,\tilde{\o}^2\geq E^2(\tilde{k}-aE^{-1}\tilde{\l})^2\}\subset\R^2$  and
$\tilde{K}_\pm =\{(\tilde{\o},\tilde{k}\in\tilde{K};\,\,\tilde{\o}=\pm E(\tilde{k}-aE^{-1}\tilde{\l})\}.$

There are functions $G_\pm(\tilde{\o},\tilde{k};h)$ such that:\\
\no 1) $G_\pm$ is a complex-valued classical symbol in $h,$ smooth in a fixed \neigh{} of $\tilde{K}.$ For $(\tilde{\o},\tilde{k})$ near $\tilde{K}$ and $|\tilde{\mu}|\leq C_\vt,$ symbol $G_\pm(\tilde{\o}+ih\tilde{\mu},\tilde{k})$ can be defined  as an asymptotic analytic Taylor series for $G_\pm$ at $(\tilde{\o},\tilde{k}).$

\no 2) For $a=0$, $G_\pm(\tilde{\o},\tilde{k};h)=-h/2+\sqrt{\tilde{\o}+h^2/4}\mp\tilde{k}.$

\no 3) $G_-(\tilde{\o},\tilde{k};h)-G_+(\tilde{\o},\tilde{k};h)=2\tilde{k}.$

\no 4) The  principal symbol $F_\pm$ of $G_\pm$  is real-valued, $\partial_{\tilde{\o}}F_\pm >0,$   $\mp\partial_{\tilde{k}} F_\pm >0$ on $\tilde{K}$ and  ${F_\pm}_{|\tilde{k}_\pm}=0.$

\no 5) For sufficiently small $h$, the set of elements $((\tilde{\o}+ih\tilde{\mu})^2,\tilde{k})$ of the joint spectrum of $([P_1]^2,P_2)$ lies within ${\mathcal O}(h)$ of
$\tilde{K}$ and coincides modulo ${\mathcal O}(h^\infty)$ with the set of solutions to the quantization condition $$\tilde{k}\in h(\Z+\frac12),\qq G_\pm(\tilde{\o}+ih\tilde{\mu},\tilde{k})\in h\N^*.$$ It is required  $G_\pm\geq 0,$ the condition $G_+\in\Z$ and $G_-\in\Z$ are equivalent, the corresponding joint eigenspaces are one dimensional.
\end{proposition}

The function $\cF^{\vt}(\tilde{l},\tilde{\l},\tilde{\nu},\tilde{k};h)$ in Proposition \ref{Prop1.6}  can be  defined  as   the solution $\tilde{\o}+ih\tilde{\mu}$ to  the equation
\[\lb{fG}G_+(\tilde{\o}+ih\tilde{\mu},\tilde{k},\tilde{\l}+ih\tilde{\nu};h)=\tilde{l}-\tilde{k}. \]  Function $\cF^{\vt}$ is uniquely defined as $ \cP_+(\tilde{\l} +ih\tilde{\nu})$ and $ \cP_-(\tilde{\l} +ih\tilde{\nu})$ have the same non-zero spectrum.

\subsection{The bottom of the well asymptotics for the angular eigenvalues.}\lb{ss-Angu}
In this sections we show (\ref{viktig-ang}).

We consider $$P_1(\tilde{\l},\tilde{\nu};h)_{\cH_{\S^2}^k} =hA_{\S^2}(\l)|_{\cH_{\S^2}^k}=U\D^{\frac12}_\vt \tilde{P}_1(\tilde{\l},\tilde{\nu};h)U^*|_{\cH_{\S^2}^k} ,$$ where
$$\D^{\frac12}_\vt\tilde{P}_1= \D^{\frac12}_\vt\left(\begin{array}{cc}
   \mu_+(\vt,hD_\vt,hD_\vp;\tilde{\l})&0 \\
   0  &   \mu_-(\vt,hD_\vt,hD_\vp;\tilde{\l})\\
         \end{array}
       \right)+\D^{\frac12}_\vt W(h). $$

Here $W(h)$ is admissible \pseudor{} of order $0$ (see \cite{HelfferRobert1983}) 
and $$\D^{\frac12}_\vt\mu_{10\pm}(\vt,\xi_\vt,\xi_\vp;\tilde{\l})=\pm\sqrt{p_{10}}=\pm\left(\D_\vt\xi_\vt^2+\frac{E^2}{\sin^2\vt\Delta_\vt}\left(\xi_\vp-\breve{\l} a\sin^2\vt \right)^2  \right)^\frac12,\qq \breve{\l}=E^{-1}\tilde{\l}.$$ So in the leading order  the quantization condition for $ P_1(\tilde{\l},\tilde{\nu};h)_{\cH_{\S^2}^k}$is the same as for
$$\left(\begin{array}{cc}
  +p_{10}^\frac12(\vt,hD_\vt,\tilde{k};\tilde{\l})&0 \\
   0  &  - p_{10}^\frac12(\vt,hD_\vt,\tilde{k};\tilde{\l})\\
         \end{array}
       \right). $$ Following \cite{Dyatlov2012}, Section B.3, we introduce new variable $y=\cos\vt$  and consider (for $h=1$)
\[\lb{angreDy}\left(\begin{array}{cc}
  +P_y^\frac12&0 \\
   0  &  - P_y^\frac12\
         \end{array}
       \right),\qq P_y=D_y(1-y^2)(1+ (E-1) y^2)D_y +\frac{ E^2(a\breve{\l}(1-y^2)-k)^2}{(1-y^2)(1+ (E-1) y^2)}.\] We
 study the bottom of the well asymptotics for the eigenvalues of $P_y^\frac12.$ The critical point for the principal symbol of $P_y$   is $(0,0).$ In order to pass from bottom of the well to the barrier-top problem we rescale the operator, introducing the papameter $y'=e^{i\pi/4}y.$ Eigenvalue $-i\o^2$ of the rescaled operator $-iP_y$ is given by
$$-i\o^2=-i E^2(k-a\Re\breve{\l} )^2-i(2m+1)\sqrt{U_0(0)}+\ldots, $$
$$U_0(0)=-V''(0)/2=k^2-(a\Re\breve{\l})^2- (E-1)(a\Re\breve{\l}-k)^2=2a\Re\breve{\l}+{\mathcal O}(a^2)>0.$$
Now, multiplying by $i$ and rescaling with $h$ we get 
$$(\tilde{\o}+ih\tilde{\mu})^2=E^2(\tilde{k}-a\Re\breve{\l} )^2+(2m+1)h\sqrt{U_0(0)}+\ldots, $$
which shows (\ref{viktig-ang}).

Note that (\ref{angreDy}) gives the leading part of $P_1.$ In \cite{Dyatlov2012} the full angular Hamiltonian is given by $P_y$ with complex parameter $\l$ instead of real $\breve{\l}.$
In order to get the quantization condition to any order, to calculate function $G_+$ in (\ref{fG}) and function  $\cF^{\vt}$ in  (\ref{1.16 }), we apply the above method to the original operators  $\cP_\pm(\l)$ defined in (\ref{orifoperators}).

\section{Radial quantization condition.}\lb{s-radial qc}

\subsection{Formulation of the radial quantization condition.} 

As in previous section, for small $h>0$ we put
\[\lb{scaling_par} \tilde{\l}=h\Re\l,\qq\tilde{\nu}=\Im\l,\qq \tilde{k}=hk,\qq \tilde{\Omega}_\pm=h\Omega_\pm,\qq  \tilde{\o}=h\Re\o,\qq\tilde{\mu}=\Im\o.\]


\begin{proposition}[Radial]\lb{Prop1.5}  Let $C_\o$ be a fixed constant  and put $K_r=(r_-+\d,r_+-\d).$ Recall that $R_r(\l,\o,k)=(H_k(\l)+\o)^{-1},$ where $H_k(\l):=[\ga(x)]^{-1}(\s_3D_x+c(x,k)-\l)\s_1.$  Then, the poles of $1_{K_r}R_r(\l,\o,k)1_{K_r}$   as a function of $\o$ in the region \[\lb{1.13} 
1 <\tilde{\l} <2,\qq |\tilde{\nu}| <\nu_0,\qq |\tilde{k}| <C_k,\qq |\tilde{\o}|,\,\,|\tilde{\mu}| <C_\o,
\] are simple with polynomial resolvent estimate $L^2\mapsto L^2,$  given modulo ${\mathcal O}(h^\infty)$ by
\[\lb{1.14}\tilde{\o}+ih\tilde{\mu}=\cF^{r,+}(m,\tilde{\l},\tilde{\nu},\tilde{k};h)\,\,\mbox{or}\,\, \tilde{\o}+ih\tilde{\mu}=\cF^{r,-}(m,\tilde{\l},\tilde{\nu},\tilde{k};h),\qq m\in\Z,\qq 0\leq m\leq C_m, \] for some constant $C_m.$ The principal part $\cF_0^{r,\pm}$ of the classical symbol $\cF^{r,\pm}$ is real-valued, independent of $m$ and $\tilde{\nu}.$ Moreover,
$$\cF^{r,+}_0(\tilde{\l},\tilde{k})=\tilde{\l}r_0F^{-\frac12}(r_0)+\frac{a}{r_0}\left[H\left(1-r_0\frac12 F^{-1}(r_0)F'(r_0)\right)-F^{-\frac12}(r_0)\tilde{k} \right]$$ and
  $$\cF^{r,-}_0(\tilde{\l},\tilde{k})=-\cF^{r,+}_0(\tilde{\l},\tilde{k})
,\qq H=\frac{\tilde{k} 4F^\frac12(r_0)r_0^2}{8Q^2-6Mr_0},\qq F(r)=1-\frac{2M}{r}+\frac{Q^2}{r^2}-\frac{\Lambda}{3}r^2.$$ 
 For $\l,k$ satisfying (\ref{1.13}), every pole $\o$ satisfies $|\tilde{\o}| >\epsilon$ for some $\epsilon >0.$
 \end{proposition} 

Note that this proposition is an analogue for the Dirac case of Proposition 1.5,  \cite{Dyatlov2012},  for  the Kerr-de Sitter black holes. In that case instead of two symbols $\cF^{r,\pm}$ there is one symbol $\cF^r$ and
$$\cF^r=ih(m+1/2)+\frac{3\sqrt 3 M}{\sqrt{1-9\L M^2}}(\tilde{\l}+ih\tilde{\nu})+{\mathcal O}(h)\qq\mbox{for}\,\, a=0,$$
$$\cF_0^r(\tilde{\l},\tilde{k})=\frac{3\sqrt 3 M}{\sqrt{1-9\L M^2}}\tilde{\l}-\frac{\tilde{k}}{M\sqrt{3}\sqrt{1- 9M^2\Lambda}} a+{\mathcal O}(a^2).$$
In the next sections we prove Proposition \ref{Prop1.5}.

\subsection{Resonance free strip} \lb{ss-ResFree}

As  in \cite{Dyatlov2011}, the definition of an outgoing solution implies the following proposition.
\begin{proposition}\lb{Prop4.4_Dirac} Let $\delta_r >0$ and  $K_r$ be the image of the set $(r_-+\delta_r,r_+-\delta_r)$ via the Regge-Wheeler change of variables $r\,\mapsto\, x.$ Suppose  that $X_0$ is  sufficiently large so that Proposition \ref{Prop4.1Dy} holds and $K_r\subset (-X_0,X_0).$ Let $f\in H^1_{\rm loc}(\R;\C^2)$ be any outgoing function in the sense of  Definition \ref{DefOutg} and suppose that $g=\cD_k^r(\l,\o)f$ is supported in $K_r.$ Then\\
\no 1) $f$ has an analytic extension to the two half-planes $\{\pm\Re z >X_0\} $ and satisfies the equation $\cD_k^r(z,\l,\o)f=0$ there. Here  $\cD_k^r(z,\l,\o)=\s_3D_z+V(z;\l,\omega,k),$ $V(z;\l,\omega,k)=c(z,k)-\l+\o \ga(z)\s_1,$ an $V(z;\l,\omega,k)$ is well defined by
(\ref{4.4}). 

\no 2) Suppose $\gamma$ be a contour in $\C$ given by $\Im z= F(\Re z),$ $x_-\leq \Re z \leq x_+$ and $F(x)=0$ for $|x| \leq X_0.$ Then  there is  a restriction to $\g$ of the holomorphic extension of $f$ by $$f_\g(x)=f(x+iF(x))$$ and $f_\g$ satisfies the equation $\cD_\g f_\g=g,$ where $$\cD_\g=\frac{1}{1+iF'(x)}\s_3D_x+V(x+iF(x);\l,\o,k). $$

\no 3) If $\g$ is as above, with $x_\pm=\pm\infty,$ and $F'(x)=c=\const$ for large $|x|.$ Then $$f_\g(x)={\mathcal O}\left(e^{\mp\Im ((1+ic)(\l-\O_\pm)x}\right),\qq\mbox{as}\,\,x\rightarrow\pm\infty. $$ Therefore, if  $\Im (1+ic)(\l-\O_\pm) >0,$ then $f_\g\in H^1(\R;\C^2).$
 \end{proposition}
 Put $h=|\Re\l|^{-1},$  $\mathfrak{s}=\sign \Re\l$ and
$$ \tilde{\nu}=\Im\l,\qq \tilde{k}=hk,\qq \tilde{\Omega}_\pm=h\Omega_\pm,\qq  \tilde{\o}=h\Re\o,\qq\tilde{\mu}=\Im\o.$$

Consider the rescaled operator
\begin{align*}&  \tilde{\cD}_k^r(\l,\o):=h\cD_k^r(\l,\o):=\s_3 h D_x+\tilde{V}_h(x,\tilde{\l},\tilde{\mu},\tilde{k}),\\  &\tilde{V}_h(x,\tilde{\l},\tilde{\mu},\tilde{k})=(\tilde{\o}+ih\tilde{\mu})\ga(x)\s_1+\tilde{c}_h(x,\tilde{k})-(\mathfrak{s}+ih\tilde{\nu}),\qq \tilde{c}_h(x,\tilde{k})=\frac{aE\tilde{k}}{r^2+a^2}+h\frac{\gq Qr}{r^2+a^2} .\end{align*}
As in  \cite{Dyatlov2011} one can show that in order to prove Theorem \ref{Th5-Dy}
it suffices to prove the following:

 {\em  Let $h$ be small enough and suppose the conditions
 \[\lb{(7.1)} |\tilde{\o}|\leq C',\qq |\tilde{k}|\leq C',\qq |\tilde{\mu}|\leq 1/C'',\qq |\nu|\leq 1/C''.\]
\vspace{2mm}

For each $g(x)\in L^2\cap\cE'(K_r)$ let $f(x)$ be solution to the equation $\tilde{\cD}_k^r(\l,\o) f=g$ which is outgoing in the sense of Definition (\ref{DefOutg}). Then
\[\lb{(7.2)} \|f\|_{L^2(K_r)}\leq Ch^{-1-\epsilon}\|g\|_{L^2}\] for some $\epsilon >0.$}

Bound (\ref{(7.2)}) will follow from (\ref{(7.4)}) which itself is proved  by the construction of an escape function and conjugation by exponential weight as in  \cite{Dyatlov2011}. We will explain below why the method also works for the Dirac operator.

We consider the leading part of $\tilde{\cD}_k^r$
$$\cD_h-1=\left(\begin{array}{cc}
         hD_x+c_0(x,\tilde{k})-\mathfrak{s}&     \tilde{\o}\ga(x)  \\
     \tilde{\o}\ga(x)  &   -hD_x+c_0(x,\tilde{k}) -\mathfrak{s}\\
         \end{array}
       \right) $$
with  the principal symbol (which is also the principal symbol of $\tilde{\cD}_k^r$) \begin{align*}&p(x,\xi)- I_2 =\left(\begin{array}{cc}
         \xi+c_0(x,\tilde{k})-\mathfrak{s}&\tilde{\o}\ga(x) \\
     \tilde{\o}\ga(x)  &   -\xi+c_0(x,\tilde{k})-\mathfrak{s} \\
         \end{array}
       \right) .
\end{align*}

 The
eigenvalues $\mu_\pm$ of  $p(x,\xi)$ are zeros of the determinant $$\det (p(x,\xi)-(\mathfrak{s}+\mu) I_2)=-\xi^2+(c_0(x,\tilde{k})-(\mathfrak{s}+\mu))^2-\tilde{\o}^2\ga^2(x) .$$ We get   $$\mu_\pm=c_0(x,\tilde{k})-\mathfrak{s}\pm\sqrt{\xi^2+\tilde{\o}^2\ga^2(x)}. $$

We  apply the method of complex scaling  (see for example \cite{AguilarCombes1971}). Consider the contour $\gamma$ in the complex plane given by $\Im x=F(\Re x),$ with $F$ defined by
\[\lb{7.3}F(x)=\left\{\begin{array}{ll}0,\qq &|x|\leq R;\\
F_0(x-R),\qq &x\geq R;\\
-F_0(-x-R),\qq &x\leq -R.\end{array}\right. \] Here $R>X_0$ is large and $F_0\in C_0^\infty (0,\infty)$ is fixed function such that $F_0' \geq 0$ and $F_0''\geq 0$ for all $x$ and $F_0'(x)=1$ for $x\geq 1.$ 

Now, let $f(x)$  be an outgoing in the sense of Definition \ref{DefOutg}
 solution to the equation $\tilde{\cD}_k^r(\l,\o) f=g\in L^2\cap\cE'(K_r).$ By Proposition \ref{Prop4.4_Dirac} we can define the restriction $f_\g$ of $f$ to $\g$ and $\tilde{\cD}_\g(\l,\o) f_\g=g,$ where $$\cD_\g=\frac{1}{1+iF'(x)}\s_3D_x+V(x+iF(x);\l,\o,k). $$ For $a$ and $h$ small $f_\gamma\in H^1(\R)$, and in order to prove (\ref{(7.2)}), it is enough to show that for each $f_\gamma\in H^1(\R)$, we have \[\lb{(7.4)} \|f_\g\|_{L^2(K_r)}\leq Ch^{-1-\epsilon}\|\tilde{\cD}_\g(\l,\o) f_\g\|_{L^2}.\]
The eigenvalues of the semi-classical principal symbols of $\tilde{\cD}_k^r(\l,\o)$ and $\cD_\g$ are given by
$\mu_\pm$, $\mu_{\g,\pm},$ respectively. Here
 \begin{align*}&\mu_\pm=c_0(x,\tilde{k})-\mathfrak{s}\pm\sqrt{\xi^2+\tilde{\o}^2\ga^2(x)},\\ &\mu_{\g,\pm}=c_0(x+iF(x),\tilde{k})-\mathfrak{s}\pm\sqrt{\frac{\xi^2}{(1+iF'(x))^2}+\tilde{\o}^2\ga^2(x+iF(x),\tilde{k})}.  \end{align*}

Now, recall  (\ref{3.17}):$$c_0(x,\tilde{k})=\Omega_{0,\pm}(\tilde{k})+c_{0,\pm} e^{2\k_\pm x}+{\mathcal O}\left(e^{4\k_\pm x}\right),\qq x\rightarrow\pm\infty,\qq \Omega_{0,\pm}=\frac{aE\tilde{k}}{r_\pm^2+a^2}, $$ and (\ref{behatinf}) $$ c_0(x,k)=c_0^\pm(e^{\k_\pm x}),\qq   \ga(x)=a^\pm(e^{\k_\pm x}),\qq \pm x >X_0,\qq c_0^\pm(0)=\O_{0,\pm},\,\, a^\pm(0)=0$$ for some holomorphic functions $c_0^\pm(w),$ $a^\pm(w).$ 

Bound (\ref{(7.4)}) follows by the construction of an escape function and conjugation by exponential weight as in  \cite{Dyatlov2011} using 
an analogue of  Proposition 7.2 in that paper:\\
{\em Suppose $x>X_0>0.$ There exists a constant $C$ such that for $R$ large enough and $\delta >0$ small enough,
\[\lb{(7.5)}\mbox{if}\,\, x\geq R+1,\,\,\mbox{then}\,\,|\mu_{\g,\pm}|\geq 1/C>0,\]
\[\lb{(7.7)}\mbox{if}\,\,|\mu_{\g,\pm}(x,\xi)|\leq \delta,\,\, \mbox{then}\,\, |\mu_{\pm}(x,\xi)|\leq C\delta,\,\,|\nabla(\Re\mu_{\pm}(x,\xi)|^2\leq C\delta.\]
\[\lb{(7.6)}\mbox{if}\,\,|\mu_{\g,\pm}(x,\xi)|\leq e^{2k_+ R},\,\, \mbox{then}\,\, \Im\mu_{\pm}(x,\xi)\leq 0.\]
Similar facts hold if $x<-X_0<0$ with $\k_-$ instead of $\k_+.$
}

It is enough  to prove this statement for $\mathfrak{s}=1.$
\\ \no Formula (\ref{(7.5)}). Suppose $x\geq R+1.$ Then 
\begin{align*}\mu_{\g,\pm}&=c_0(x+iF(x),\tilde{k})-1\pm\sqrt{-i\xi^2/2+\tilde{\o}^2\ga^2(x+iF(x),\tilde{k})}\\&=\O_{0,+}-1\pm\sqrt{-i\xi^2/2}+{\mathcal O}(e^{-R})=\O_{0,+}-1\pm(1-i)|\xi |/2+{\mathcal O}(e^{-R})\end{align*} and
$$|\mu_{\g,\pm}|^2=(\O_{0,+}-1\pm |\xi |/2)^2+|\xi|^2/4+{\mathcal O}(e^{-R}).$$ Minimizing the function $f_\pm(x)=(\O_{0,+}-1\pm x)^2+x^2$ we get $f_+(x)\geq f(\mp\frac23 (\O_{0,+}-1))=C_\pm>0,$ where $C_+=\frac59(\O_{0,+}-1))^2,$  $C_-=\frac{29}{9}(\O_{0,+}-1))^2.$

Taking $R$ large enough we get (\ref{(7.5)}). 

For the rest of the proof we assume that $R\leq x\leq R+1.$ 

\no Formulas (\ref{(7.7)}), (\ref{(7.6)}).  Suppose that $|\mu_{\g,\pm}(x,\xi)|\leq \delta,$ then as in \cite{Dyatlov2011} we get $$ \frac{\xi^2}{(1+iF'(x))^2}=(\O_{0,+}-1)^2+{\mathcal O}(\delta+e^{-R}F(x))$$ which implies $F'(x)\leq c\delta$ and leads to (\ref{(7.5)}) and (\ref{(7.6)}) (by repeating the proof in \cite{Dyatlov2011}).\BBox

Then we can apply the method of \cite{Dyatlov2011} involving construction of escape function and conjugation by exponential weights.

This achieves the proof of  Theorem \ref{Th5-Dy} and shows a resonance free strip. The only trapping in our situation is normally hyperbolic and generate the radial poles  which are studied in the next section.

\subsection{The trapping point asymptotics for the radial resonances.}

In this section we study the  radial poles in the region (\ref{1.13}) appearing in the  trapping case.  In the non-trapping case
 there is an arbitrary large strip free of radial poles. 

We rescale operator $\cD_k^r(\l,\o)$ and use notation (\ref{scaling_par}).
Let \begin{align*}&  \tilde{\cD}_k^r(\l,\o):=h\cD_k^r(\l,\o):=\s_3 h D_x+\tilde{V}_h(x,\tilde{\l},\tilde{\mu},\tilde{k}),\\  &\tilde{V}_h(x,\tilde{\l},\tilde{\mu},\tilde{k})=(\tilde{\o}+ih\tilde{\mu})\ga(x)\s_1+\tilde{c}_h(x,\tilde{k})-(\tilde{\l}+ih\tilde{\nu}),\qq \tilde{c}_h(x,\tilde{k})=\frac{aE\tilde{k}}{r^2+a^2}+h\frac{\gq Qr}{r^2+a^2} .\end{align*}
We consider  the leading part of  $\tilde{\cD}_k^r(\l,\o)$
\[\lb{PrinsPart}\cD_h(\tilde{\l},\tilde{\o},\tilde{k})=\cD_h(\tilde{\o},\tilde{k})-\tilde{\l}=\left(\begin{array}{cc}
         hD_x+c_0(x,\tilde{k})-\tilde{\l}&\tilde{\o}\ga(x) \\
     \tilde{\o}\ga(x)  &   -hD_x+c_0(x,\tilde{k}) -\tilde{\l}\\
         \end{array}
       \right). \]
The principal symbol of $\cD_h(\tilde{\l},\tilde{\o},\tilde{k})$ is given by \[\lb{prinss}p(x,\xi; \tilde{\l},\tilde{\o},\tilde{k})=p(x,\xi; \tilde{\o},\tilde{k})-\tilde{\l} I_2 =\left(\begin{array}{cc}
         \xi+c_0(x,\tilde{k})-\tilde{\l}&\tilde{\o}\ga(x) \\
     \tilde{\o}\ga(x)  &   -\xi+c_0(x,\tilde{k})-\tilde{\l} \\
         \end{array}
       \right) , 
\]
with the
eigenvalues   $$\mu_\pm=\mu_\pm(x,\xi; \tilde{\l},\tilde{\o},\tilde{k})=c_0(x,\tilde{k})-\tilde{\l}\pm\sqrt{\xi^2+\tilde{\o}^2\ga^2(x)}. $$

Note the oddness property of the principal symbol:
$\mu_-(x,\xi; \tilde{\l},\tilde{\o},\tilde{k})=-\mu_+(x,\xi; -\tilde{\l},\tilde{\o},-\tilde{k})$ and that it does not depend on the sign of $\tilde{\o}.$

We study the trapping properties of the  
eigenvalues  of  $\mu_\pm.$   Bicharacteristics are the solutions of the Hamilton system
\begin{align*}\left\{\begin{array}{l}
        \dot{x}_\pm=\partial_\xi\mu_\pm(x_\pm(t),\xi_\pm(t)),\\
 \dot{\xi}_\pm=-\partial_x\mu_\pm(x_\pm(t),\xi_\pm(t))\\
         \end{array}
       \right.\qq\Leftrightarrow\qq\left\{\begin{array}{l}
        \dot{x}_\pm=\pm\frac12\frac{1}{\sqrt{\xi^2+\tilde{\o}^2a^2(x)}}2\xi,\\
 \dot{\xi}_\pm=-c_0'(x)\mp \frac12\frac{1}{\sqrt{\xi^2+\tilde{\o}^2a^2(x)}}\tilde{\o}^2(\ga^2(x))'\\
         \end{array}
       \right.\end{align*}
with initial condition $(x_\pm(0),\xi_\pm(0))=(x_\pm^0,\xi_\pm^0).$

 The only critical point of the Hamiltonian $\mu_\pm$ is 
 given by $(x,\xi)=(x_\pm,0),$ where $x_\pm$ is solution to $ \left[c_0\pm |\tilde{\o}|\ga\right]'(x)=0.$

 Note that  by symmetry of the symbol we get that passing from $x_+$ to $x_-$ is equivalent to the change of sign $\tilde{k}\,\,\mapsto\,\,-\tilde{k}.$

At $(x_\pm,0)$ we get $$ (\mu_\pm)''_{xx}=(c\pm |q|)''(x_\pm),\qq (\mu_\pm)''_{\xi\xi}=\pm |q(x_\pm)|^{-1},\qq(\mu_\pm)''_{x\xi}=0.$$
 $$ (\mu_\pm)''_{xx}=(c_0\pm |\tilde{\o}|\ga)''(x_\pm),\qq (\mu_\pm)''_{\xi\xi}=\pm (|\tilde{\o}|\ga(x_\pm))^{-1},\qq(\mu_\pm)''_{x\xi}=0.$$

Note that for $a$ small enough $(\mu_+)''_{xx}(x_+,0) <0, $ $(\mu_-)''_{xx}(x_-,0) >0,$ and the critical points are of hyperbolic type.

Now, using that there is an arbitrary large strip free of radial poles in the nontrapping cases we know that the only radial poles in the region (\ref{1.13}) appear in the trapping case. Then we may assume  that $| \tilde{\o}^{2}- \tilde{\o}_0^{2}(\tilde{\l},\tilde{k})| <\varepsilon_r,$
where $\tilde{\o}_0^{-2}$ is the value of the function
\[\lb{valf}F_V(r;\tilde{\l},\tilde{k})=\frac{\ga^2}{\tilde{\l}-c_0(x)}=\frac{\Delta_r}{\left[\tilde{\l}(r^2+a^2)-aE\tilde{k}\right]^2} \]
at its only maximum point. 
Under assumption (\ref{1.13}) , $1/C\leq \tilde{\o}_0^{2}\leq C$ for some constant $C.$

We see that the only critical points of $\mu_\pm$  in the set $\{ \mu_\pm \geq -\varepsilon_r^3\}$  are non-degenerate hyperbolic critical points at $(x_\pm,0).$

 Now, we establish a microlocal form for $\tilde{\cD}_k^r(\l,\o)$ near  $(x,\xi)=(x_\pm,0).$ Firstly, we will diagonalize the 2 by 2 matrix $\tilde{\cD}_k^r(\l,\o)$ modulo $h^\infty,$ using the method of M. Taylor explained in \cite{HelfferSjostrand1990} and \cite{BruneauRobert1999}. 

The principal symbol $p(x,\xi; \tilde{\l},\tilde{\o},\tilde{k})
$ (\ref{prinss}) can be diagonalized exactly via unitary Foldy-Wouthuysen transform $U_0$ (see \cite{Thaller1992})
\[\lb{diagS}U_0(x,\xi)^*p(x,\xi; \tilde{\l},\tilde{\o},\tilde{k}) U_0(x,\xi)=\ma \mu_+(x,\xi) &0\\ 0 &\mu_-(x,\xi)\am,\] with $U_0=U_0(x,\xi)=$ $$\frac{1}{\sqrt{2}}\frac{1}{\sqrt{\D^2+( \xi +c_0(x))\D}}\left(\begin{array}{cc}
        \xi+c_0(x) +\D&-(\tilde{\o}\ga(x))\\
   \tilde{\o}\ga(x) & \xi +c_0(x)+\D  \\
         \end{array}
       \right),\,\, \Delta=\sqrt{(\xi+c_0(x))^2+(\tilde{\o}\ga(x))^2}.$$

By Weyl quantization of the symbol identity (\ref{diagS}) we diagonalize the principle part of the operator $\cD_h=\cD_h(\tilde{\l},\tilde{\o},\tilde{k})$ in (\ref{PrinsPart}) 
 $$\tilde{\cD}_h(\l)=U(x,hD_x)^*\cD_hU(x,hD_x)= \left(\begin{array}{cc}
   \mu_+(hD,x)&0 \\
   0  &   \mu_-(hD,x)\\
         \end{array}
       \right)+hR(h), $$ where \pseudor{} $R(h)$ has classical symbol and $\|R(h)\|={\mathcal O}(1).$

Now, by iteration as in Section 3.1 of \cite{HelfferSjostrand1990} and Section 4.1 of \cite{BruneauRobert1999} we can diagonalise operator $\tilde{\cD}_k^r$ to the infinite order in $h.$

\begin{proposition}[Decoupling]\lb{Prop-decoupling} There exist unitary $U=U(x,hD_x,h)$  such that 
\[\lb{decoupling} U^*\tilde{\cD}_k^r(\l,\o) U= \left(\begin{array}{cc}
   \mu_{+}(x,hD,h)&0 \\
   0  &   \mu_{-}(x, hD,h)\\
         \end{array}
       \right)+{\mathcal O}(h^\infty)\] in $\cL(L^2,L^2).$ Here
 $\mu_{\pm}(x,hD,h)$ is \pseudor{} with classical symbol $\mu_{\pm}(x,\xi,h)=\sum_{k=0}^{\infty}h^k\mu_{k,\pm}(x,\xi)$ and $\mu_{0,\pm}(x,\xi)=\mu_{\pm}(x,\xi).$
\end{proposition}

We pass now to reduction of $ \mu_\pm(x,hD,h)$ to the microlocal normal form.
It is enough to consider $ \mu_+(x,hD,h)$ as  the other case is similar. 

The proof of Proposition 4.3 in \cite{Dyatlov2012} can  also be applied here using that $\mu_\pm(x,\xi)$ is of the same type as symbol $p_0(x,\xi)$ in Section 4.3 in  \cite{Dyatlov2012}.
 The idea is to apply an analogue of  Theorem 12 in \cite{CdV_Parisse1994} to the symbol
  $\mu_+(x,\xi)-\mu_+(x_+,0).$ Here $ \mu_+(x,\xi)$ is the principal symbol of $ \mu_+(x,hD,h).$ Then 
 operator $$\mu_+(x,hD_x)-\mu_+(x_+,0)$$  can be transformed microlocally near $(x_+,0)$ into $f(hxD_x)$ by constructing a simplectomorphism $\Phi$ from a \neigh{} of $$K_0=\{|x-x_+|\leq \epsilon_0,\,\,|\xi|\leq\epsilon_0\}\subset T^*\R$$ onto a \neigh{} of the origin in $T^*\R,$ $\Phi(x_+,0)=(0,0),$ and operators $B_1,$ $B_2$ quantizing $\Phi$ near $K_0\times\Phi(K_0).$

Let $f_0$ be the principal part of $f.$ Then, $\mu_+(x,\xi)\circ\Phi^{-1}=f_0(x\xi).$  The level set $ \{ \mu_+(x,\xi)=\mu_+(x_+,0)\}$ at the trapped energy 
containes in particular the outgoing trajectory $$\{ x>x_+,\,\, 
\xi=\sqrt{\left(c_0(x_+)-c_0(x)+\tilde{\o}\ga(x_+)\right)^2-(\tilde{\o}\ga(x))^2}\}.$$ Then we can choose $\Phi$ mapping this trajectory  into $\{x>x_+,\,\,\xi=0\}.$

The function $f(s;h)$ is not uniquely defined; however, its Taylor decomposition at $s=0,$ $h=0$ is and we can compute
\begin{align*}
\mu_+(x,\xi)=&\mu_+(x_+,0)+\frac12|\tilde{\o}\ga(x_+)|^{-1}\xi^2-\frac12\left|(c_0+ |\tilde{\o}|\ga)''(x_+)\right|(x-x_+)^2+\ldots \\
=&\mu_+(x_+,0)+\frac12|\tilde{\o}\ga(x_+)|^{-1}\left[\xi^2-|\tilde{\o}|\ga(x_+)\left|(c_0+ |\tilde{\o}|\ga)''(x_+)\right|(x-x_+)^2\right]+\ldots,
\end{align*} where we used that $(\mu_+)''_{xx}(x_+,0)= (c_0+ |\tilde{\o}|\ga)''(x_+)<0.$ Denote $$\sqrt{ |\tilde{\o}|\ga(x_+)\left|(c_0+ |\tilde{\o}|\ga)''(x_+)\right|}=\sigma. $$
Performing the first symplectic change of variables 
$x\mapsto\tilde{x}:=\s^\frac12 (x-x_+),$ $\xi\mapsto\tilde{\xi}=\s^{-\frac12}\xi  $
we get  \begin{align*}
\mu_+(x,\xi)=&\mu_+(x_+,0)+\frac12|\tilde{\o}\ga(x_+)|^{-1}\s\left[\tilde{\xi}^2-\tilde{x}^2\right]+\ldots
\end{align*} 
Then the second linear symplectic change of variables
$\tilde{\xi}-\tilde{x}=\sqrt{2}\eta,$ $ \tilde{\xi}+\tilde{x}=\sqrt{2}y$ leads to
 \begin{align*}
\mu_+(x,\xi)=&\mu_+(x_+,0)+|\tilde{\o}\ga(x_+)|^{-1}\s\left[y\eta \right]+\ldots
\end{align*} The composition of these two linear maps is the linear part of the symplectomorphism $\Phi.$ Recall $\mu_+(x,\xi)\circ\Phi^{-1}=f_0(x\xi)$ and
we get for $f_0(s),$ $s=y\eta$
\begin{align*} f_0(s)=&\mu_+(x_+,0)+|\tilde{\o}\ga(x_+)|^{-1}\s s+{\mathcal O}(s^2)\\
=& c_0(x_+)-\tilde{\l}+ \tilde{\o}\ga(x_+)+|\tilde{\o}\ga(x_+)|^{-\frac12}\sqrt{ \left|(c_0+ |\tilde{\o}|\ga)''(x_+)\right|} s+{\mathcal O}(s^2).
\end{align*}

Similarly, we can consider the case of critical point $(x,\xi)=(x_-,0).$

As result, by applying Proposition \ref{Prop-decoupling}, we get an analogue for  Dirac operators of Proposition 4.3 in \cite{Dyatlov2012}. \begin{proposition}\lb{Prop4.3} The microlocal normal form of $\tilde{\cD}_k^r(\l,\o)$ is given by
 $$\left(\begin{array}{cc}
 hxD_x-\beta_+&0 \\
   0  &  -hxD_x+\beta_-\\
         \end{array}
       \right) ,$$ where $\beta_\pm=\beta_\pm(\tilde{\l},\tilde{\nu},\tilde{\o},\tilde{\mu},\tilde{k};h)$ is classical symbol. Moreover, the principal part $\beta_{0,\pm}$ of $\beta_\pm$ is real-valued, independent of $\tilde{\nu},$ $\tilde{\mu},$ and vanishes if and only if $\tilde{\o}=\tilde{\o}_0(\tilde{\l},\tilde{k}),$  where $\tilde{\o}_0^{-2}$ is defined by   (\ref{valf}). Moreover,
$$\beta_{0,+}=-\frac{(c_0(x_+)-\tilde{\l}+ \tilde{\o}\ga (x_+))\sqrt{|\tilde{\o}|\ga(x_+)}}{\sqrt{ \left|(c_0+|\tilde{\o}|\ga)''(x_+)\right|}}+{\mathcal O}\left( c_0(x_+)-\tilde{\l}+ \tilde{\o}\ga(x_+)\right)^2 $$
and  $\beta_{0,-}(\tilde{\l},\tilde{\nu},\tilde{\o},\tilde{k})=\beta_{0,+}(-\tilde{\l},\tilde{\nu},\tilde{\o},-\tilde{k}).$
\end{proposition}

The radial quantization symbol $\cF^{r,+}(m,\tilde{\l},\tilde{\nu},\tilde{k};h)$ can be obtained as the solution $\tilde{\o}+ih\tilde{\mu}$ to the equation
$$\beta_+(\tilde{\l},\tilde{\nu},\tilde{\o},\tilde{\mu},\tilde{k};h)=-ihm,\qq m\in\Z,\qq 0\leq m\leq C_m.$$
So we get in the leading order   $h=0$  near $x=x_+$
\[\lb{star} -(c_0(x_+)-\tilde{\l}+ \tilde{\o}\ga(x_+))(|\tilde{\o}|a)^{\frac12}(x_+)=-ihm\sqrt{ \left|(c_0+\tilde{\o}\ga)''(x_+)\right|},\qq c_0=\frac{a E\tilde{k}}{r^2+a^2}.\] We need to solve it \wrt{} $\tilde{\o}.$ 
Let $r_+$ be defined via Regge-Wheeler transform $x_+=x(r_+)$ and let $x_0=x(r_0),$ where 
  $r_0=\frac{3M}{2}+\sqrt{\left(\frac{3M}{2}\right)^2-2Q^2},$ is the critical point for $a=0,$ i.e. $F'(r_0)r_0-2F(r_0)=0.$ Then we find 
$$r_+\sim r_0+ \frac{a}{\tilde{\o}}H+{\mathcal O}(a^2),\qq H=\frac{\tilde{k}4F^\frac12(r_0)}{r_0 (F'(r)r-2F(r))'(r_0)}=\frac{\tilde{k} 4F^\frac12(r_0)r_0^2}{8Q^2-6Mr_0} $$
and
$$r_+:=r(x_+)=r_0+\frac{a}{\tilde{\o}}\frac{\tilde{k} 4F^\frac12(r_0)r_0^2}{8Q^2-6Mr_0}+{\mathcal O}(a^2(|\tilde{k}|^2+|\tilde{\o}|^2 )).$$

Considering order $h^0$ in (\ref{star}) we get for small $a$ $$\cF^{r,+}_0(\tilde{\l},\tilde{k})=\frac{\tilde{\l}}{\ga(x_0)}+{\left(\cF^{r,+}_0(\tilde{\l},\tilde{k})\right)'_a}_{|_{a=0}}a+{\mathcal O}(a^2). $$ Then we expand  
$$r_+^2=r_0^2+2r_0H\frac{a}{\tilde{\o}}+{\mathcal O}(a^2),\qq F^\frac12(r_+)=F^\frac12(r_0)+ \frac12 F^{-\frac12}(r_0)F'(r_0)H\frac{a}{\tilde{\o}}+{\mathcal O}(a^2).$$ 

Now, consider the following equation in $\tilde{\o}:$
$ c_0(x_+)-\tilde{\l}+ \tilde{\o}\ga(x_+)=0.$ As $c_0= \frac{a E\tilde{k}}{r^2+a^2}=\frac{a\tilde{k}}{r_0^2}+{\mathcal O}(a^2)$ and
$$\ga(x_+)=\frac{F^\frac12(r_+)}{r_+}=\frac{F^\frac12(r_0)}{r_0}+\frac{a}{\tilde{\o}} H\left(\frac12\frac{F^{-\frac12}(r_0)F'(r_0)}{r_0}-\frac{F^\frac12(r_0) }{r_0^2}\right)+{\mathcal O}(a^2), $$ we get equation in  $\tilde{\o}:$
$$\frac{a\tilde{k}}{r_0^2}-\tilde{\l}+\tilde{\o}\left[\frac{F^\frac12(r_0)}{r_0}+\frac{a}{\tilde{\o}} H\left(\frac12\frac{F^{-\frac12}(r_0)F'(r_0)}{r_0}-\frac{F^\frac12(r_0) }{r_0^2}\right)\right]={\mathcal O}(a^2). $$

We get 
\[\lb{viktig-rad}\cF^{r,+}_0(\tilde{\l},\tilde{k})=\tilde{\l}r_0F^{-\frac12}(r_0)+\frac{a}{r_0}\left[H\left(1-r_0\frac12 F^{-1}(r_0)F'(r_0)\right)-F^{-\frac12}(r_0)\tilde{k} \right]\]

The radial quantization symbol $\cF^{r,-}(m,\tilde{\l},\tilde{\nu},\tilde{k};h)$ can be obtained as the solution $\tilde{\o}+ih\tilde{\mu}$ to the equation
$$\beta_-(\tilde{\l},\tilde{\nu},\tilde{\o},\tilde{\mu},\tilde{k};h)=-ihm,\qq m\in\Z,\qq 0\leq m\leq C_m.$$ In the leading order    near $x=x_-$ we have 
\[\lb{star_minus_bis}(c_0(x_-)-\tilde{\l}- \tilde{\o}\ga (x_-))\sqrt{|\tilde{\o}|\ga(x_-)}=-ihm\sqrt{ (c_0-|\tilde{\o}|\ga)''(x_-)},\qq c_0=\frac{a E\tilde{k}}{r^2+a^2}.\] As for the case $x_+,$ we  solve (\ref{star_minus_bis})  \wrt{} $\tilde{\o}.$ Then we get 
\[\lb{viktig-rad_minus_bis}\cF^{r,-}_0(\tilde{\l},\tilde{k})=-\tilde{\l}r_0F^{-\frac12}(r_0)-\frac{a}{r_0}\left[H\left(1-r_0\frac12 F^{-1}(r_0)F'(r_0)\right)-F^{-\frac12}(r_0)\tilde{k} \right]\] and
  \[\lb{symF_bis}\cF^{r,-}_0(\tilde{\l},\tilde{k})=-\cF^{r,+}_0(\tilde{\l},\tilde{k})
.\]

We compare with  \cite{Dyatlov2012}. Put $Q=0.$ Then $$H=-\frac{\tilde{k} 4F^\frac12(r_0)r_0}{6M},\qq r_0=3M,\qq F(r)=1-\frac{2M}{r}-\frac{\Lambda}{3}r^2,\qq F'(r)=\frac{2M}{r^2}-\frac23\L r$$ and
$$F(r_0)=\frac13- 3M^2\Lambda,\qq F'(r_0)=\frac29 \frac{1}{M}-2\L M,\qq H= -\tilde{k} \frac{2}{\sqrt{3}}\sqrt{1- 9M^2\Lambda}.$$  
We get the leading term at $a=0$ $$\tilde{\l}r_0F^{-\frac12}(r_0)=\frac{3M}{\sqrt{\frac13- 9M^2\Lambda}}\tilde{\l}=\frac{3\sqrt{3}M}{\sqrt{1- 9M^2\Lambda}}\tilde{\l}. $$
Order $a^1:$
\begin{align*}&F^{-\frac12}(r_0)\left[H\left(F^\frac12(r_0)r_0^{-1}-\frac12 F^{-\frac12}(r_0)F'(r_0)\right)-\tilde{k}r_0^{-1} \right]\\
&=H\left(r_0^{-1}-\frac12 F^{-1}(r_0)F'(r_0)\right)-\tilde{k}r_0^{-1} F^{-\frac12}(r_0)\\
&=-\tilde{k} \frac{2}{\sqrt{3}}\sqrt{1- 9M^2\Lambda}\left(\frac{1}{3M}-\frac{\frac19 \frac{1}{M}-\L M}{\frac13- 3M^2\Lambda}\right)-\frac{\tilde{k}}{3M}\frac{1}{\sqrt{\frac13- 3M^2\Lambda}}\\&=-\frac{\tilde{k}}{M\sqrt{3}\sqrt{1- 9M^2\Lambda}}.
\end{align*} Then
$$\cF_0^{r,\pm}(\tilde{\l},\tilde{k})^2=\frac{27M^2}{1- 9M^2\Lambda}\tilde{\l}^2-\frac{6a\tilde{\l}\tilde{k}}{1- 9M^2\Lambda}+{\mathcal O}(a^2)$$
coincides with expression obtained 
in \cite{Dyatlov2012}, Proposition 1.5.

\section{Combination of both quantization conditions.}\lb{s-together}

In this section we calculate the poles of the resolvent $R(\l,k)$ by combining the poles of the angular and radial resolvents (see Remark \ref{Remark-poles-together}) 
Proposition \ref{Prop1.4Dy}. We define $\cF^{\l,\pm}(m, \tilde{l},\tilde{k};h)$ to be solution $\tilde{\l}+ih\tilde{\nu}$ to the equation \[\lb{(1.17}\left(\cF^{r,\pm}(m,\tilde{\l},\tilde{\nu},\tilde{k};h)\right)^2=\cF^\vt(\tilde{l},\tilde{\l},\tilde{\nu},\tilde{k};h).\] 
Note that such solutions are not unique.

The proof of (\ref{1.8}) repeats the arguments from  the proof of Proposition 1.2 in \cite{Dyatlov2012} (with few straightforward modifications)  and we do not need to reproduce it here.  Note that in our case we need to choose the contour $\gamma=\gamma_-\cup\gamma_+$ in Fig. \ref{admcontour} containing two parts: with $\gamma_+$ as in  Figure 1 in  \cite{Dyatlov2012} and $\gamma_-$ the mirror image with respect to the imaginary axis of the contour $\gamma_+.$
We prove 
  (\ref{main-formula_bis}) .

As   $\cF^{r,-}_0(\tilde{\l},\tilde{k})=-\cF^{r,+}_0(\tilde{\l},\tilde{k})$ (see (\ref{symF_bis}) it is enough to consider $\cF^{r,+}_0(\tilde{\l},\tilde{k}).$
We use   (\ref{viktig-ang}) from Proposition \ref{Prop1.6}: $$\cF_0^{\vt}(\pm \tilde{k},\tilde{\l},\tilde{k})=E^2(\tilde{k}-\frac{a}{E}\tilde{\l})^2,\qq \partial_{\tilde{l}}\cF_0^{\vt}(\pm \tilde{k},\tilde{\l},\tilde{k})=\pm\tilde{k}+{\mathcal O}(a^2),$$ implying $\partial_{\tilde{k}}\cF^\vt_0(\pm \tilde{k},\tilde{\l},\tilde{k})=-2a\tilde{\l}+{\mathcal O}(a^2).$ Then we get \[\lb{derang}\partial_{\tilde{k}}\sqrt{\cF^\vt_0(\pm \tilde{k},\tilde{\l},\tilde{k})}=-\frac{a\tilde{\l}}{E\tilde{k}}+{\mathcal O}(a^2).\]

We use  (\ref{viktig-rad}).
Suppose ${\tilde \l}$ is solution to the equation
$$\cF_0^{r}(\tilde{\l},\tilde{k})=\pm\left(\cF_0^{\vt}(\tilde{l},\tilde{\l},\tilde{k})\right)^\frac12.$$

In $
\cF_0^{\vt}(\pm \tilde{k},\tilde{\l},\tilde{k})=E^2(\tilde{k}-\frac{a}{E}\tilde{\l})^2,$   $\tilde{\l}$ can be exchanged with $$\tilde{\l}_{|a=0, \tilde{l}=\pm\tilde{k}}=\pm Ez_0 \tilde{k},\qq z_0=\left(\frac{M}{r_0^3}-\frac{Q^2}{r_0^4}-\frac{\L}{3}\right)^\frac12.$$ Then $\tilde{k}-\frac{a}{E}\tilde{\l}=\tilde{k}(1\mp a z_0).$ 


We have $$\partial_{\tilde{k}}\cF_0^r(\tilde{\l},\tilde{k})=r_0F^{-\frac12}(r_0)\partial_{\tilde{k}}\tilde{\l}+\frac{a}{r_0}\left[\frac{4F^\frac12(r_0)r_0^2}{8Q^2-6Mr_0}\left(1-r_0\frac12 F^{-1}(r_0)F'(r_0)\right)-F^{-\frac12}(r_0)\right].$$ Comparing with (\ref{derang}) we get that
$$r_0F^{-\frac12}(r_0)\partial_{\tilde{k}}\tilde{\l}+\frac{a}{r_0}\left[\frac{4F^\frac12(r_0)r_0^2}{8Q^2-6Mr_0}\left(1-r_0\frac12 F^{-1}(r_0)F'(r_0)\right)-F^{-\frac12}(r_0)\right]=\mp\frac{a\tilde{\l}_{|a=0}}{E\tilde{k}}.$$ Here $$\tilde{\l}_{|a=0}=Ez_0 \tilde{k},\qq z_0=\left(\frac{M}{r_0^3}-\frac{Q^2}{r_0^4}-\frac{\L}{3}\right)^\frac12.$$
We get 
  (\ref{main-formula_bis}) 
$$\partial_{\tilde{k}}\tilde{\l}=-\frac{a}{r_0^2}\left[\frac{4F(r_0)r_0^2}{8Q^2-6Mr_0}\left(1-r_0\frac12 F^{-1}(r_0)F'(r_0)\right)-1\right]\mp az_0r_0^{-1}F^{\frac12}(r_0) .$$

If $Q=0,$ we get $1-r_0\frac12 F^{-1}(r_0)F'(r_0)=0$ and we get formula (0.4) in \cite{Dyatlov2012} if we choose  sign $-$.


\newpage \no{\LARGE {\em Extention to the massive Dirac fields.}}

\section{Preliminaries.}\lb{s-Prelmass}

\subsection{Evolution equation and separation of variables for real $\l$.}
We consider the  charged Dirac fields with mass $\gm$ represented by 4-spinors $\psi$ belonging to the Hilbert space $\cH=L^2\left(\R\times\S^2,dxd\vt d\vp;\,\C^4 \right)$ and satisfying the evolution equation
$$i\partial_t\psi=\cD\psi,\qq \cD=J^{-1}\cD_0,\qq \cD_0=\G^1D_x+\gb(x)\G^0+c(x,D_\varphi)+\ga(x)\cD_{\S^2}.$$

Here $\cD_{\S^2}$ is an angular Dirac operator on 2-sphere $\S^2$ given in (\ref{3.10mass})
$$\cD_{\S^2}=\sqrt{\D_\vt}\left[\G^2\left(D_\vt+\frac{i\L a^2\sin(2\vt)}{12\Delta_\vt}\right)+\G^3\frac{E}{\Delta_\vt\sin\vt}D_\vp\right]-a\gm\cos\vt \G^5.$$

 Now, note the following identities \begin{align*}&(\cD-\l)\psi=\phi\qq\Leftrightarrow\qq J^{-1} \left(\G^1D_x+\gb(x)\G^0+c(x,D_\varphi)+\ga(x)\cD_{\S^2}-\l J\right)\psi=\phi\\
& \qq\Leftrightarrow\qq\left[\G^1D_x+\gb(x)\G^0+c(x,D_\varphi)+\ga(x)\cD_{\S^2}-\l \left(I_4+\ga(x)b(\vt)\G^3\right)\right]\psi= J\phi\\
& \qq\Leftrightarrow\qq\left[\G^1D_x+\gb(x)\G^0+c(x,D_\varphi)-\l+\ga(x)\left\{\cD_{\S^2}-\l b(\vt)\G^3\right\}\right]\psi= J\phi.
\end{align*}

Then the stationary Dirac equation $\cD\psi=\l\psi$ can be re-written as $\cD(\l)\psi=0,$ where
\[\lb{3.26mass}\cD(\l)=\G^1D_x+\gb(x)\G^0+c(x,D_\varphi)-\l+\ga(x)\left\{\cD_{\S^2}-\l b(\vt)\G^3\right\},
\]  and $(\cD-\l)\psi=f$ is equivalent to $\cD(\l)\psi=Jf$ and $\psi=(\cD-\l)^{-1}f=[\cD(\l)]^{-1}Jf.$
Therefore, we have for the resolvent of $\cD$
$$  (\cD-\l)^{-1}=[\cD(\l)]^{-1}J.$$
Let \[\lb{3.27mass} A_{\S^2}(\l)=\cD_{\S^2}-\l b(\vt)\G^3,\qq \cH_{\S^2}=L^2(\S^2,\,d\vt d\vp;\,\C^4).\] 

For real $\l,$ we decompose $\cH_{\S^2}$  onto the angular modes $\{e^{ik\vp}\}_{ k\in \frac12+\Z}$ that are eigenfunctions for $D_\vp$ with anti-periodic boundary conditions (see \cite{BelgiornoCacciatori2009}). Then
\[\lb{3.12mass}
\cH_{\S^2}=\bigoplus_{ k\in \frac12+\Z}\cH_{\S^2}^k,\qq \cH_{\S^2}^k=L^2( (0,\pi),\,d\vt;\,\C^4).
 \] The reduced subspaces $\cH_{\S^2}^k$ remain invariant under the action of $A_{\S^2}(\l)$
and we denote $A_k(\l)=A_{\S^2}(\l)_{|\cH_{\S^2}^k}.$ We have explicitly  
\[\lb{3.30mass}A_k(\l)= \sqrt{\D_\vt}\left[\G^2\left(D_\vt+\frac{i\L a^2\sin(2\vt)}{12\Delta_\vt}\right)+\G^3\left(\frac{kE}{\Delta_\vt\sin\vt} -\l\frac{a\sin\vt}{\Delta_\vt}\right)\right]-a\gm\cos\vt \G^5.\]

For each $k\in 1/2 +\Z,$ operator  $A_k(\l)$ is self-adjoint and has discrete simple spectrum $\sigma(A_k(\l))=\{\mu_{kl}(\l)\}_{l\in\Z^*}$ with associated set of eigenfunctions $\{u_{kl}^\l\}_{l\in\Z^*},$ 
$$A_k(\l) u_{kl}^\l(\vt)=\mu_{k,l}(\l)u_{kl}^\l(\vt).$$ Here $\Z^*=\Z\setminus\{0\}.$ Since $\sigma(A_k(\l))$ is discrete, it has no accumulation point and thus
$$\forall k\in 1/2 +\Z,\qq |\mu_{k,l}(\l)| \rightarrow \infty \qq\mbox{as}\qq l\rightarrow \pm\infty.$$
Eigenvalues $\mu_{k,l}(\l)$ of $A_k(\l)$ are also  the  eigenvalues  of $A_{\S^2}(\l)$ with eigenfunctions\\
$Y_{k,l}(\l):=Y_{k,l}(\l,\vt,\vp)=u_{kl}^\l(\vt) e^{ik\vp}.$


The analogue of Lemma \ref{L-DNangular} is still valid in the massive case. 
By using the cylindrical symmetry we decompose the Hilbert space $\cH$ onto the angular modes $\{e^{ik\vp}\}_{ k\in \frac12+\Z},$ 
\[\lb{3.12mass_bis}
\cH=\bigoplus_{ k\in \frac12+\Z}\cH_k,\qq \cH_k=L^2(\R\times (0,\pi),\,dxd\vt;\,\C^4)=L^2(\R;\C^4)\otimes L^2((0,\pi),\,d\vt;\,\C^4).
 \]
 

We choose half-integers $k$  as we want  the anti-periodic conditions in variable $\vp:$ the spinors change the sign after a complete rotation (see \cite{BelgiornoCacciatori2009}). Note that $$\mu_{k,-l}(\l)=-\mu_{kl}(\l),\qq Y_{k,-l}(\l)=\G_1 Y_{kl}(\l).$$

Using these results we have the decomposition 
$$\cH=\bigoplus_{(k,l)\in I}\cH_{kl}(\l),\qq I=\left(\frac12+\Z\right)\times\N^*,\qq \cH_{kl}(\l)=L^2(\R;\C^4)\otimes Y_{kl}(\l).$$ We choose $I$ instead of $K$  in order to have subspaces $\cH_{kl}$ remain invariant under the action of $\cD(\l)$ (see Section 3.2 in \cite{DaudeNicoleau2015} for details).

Let ${\displaystyle \cD_k(\l):=\cD(\l)_{|\cH_k}=\G^1D_x+\gb(x)\G^0+c(x,D_\varphi)-\l+\ga(x)A_k(\l)}$ be restriction of $\cD(\l)$ to $\cH_k.$

Radial operator $\G^1D_x+\gb(x)\G^0+c(x,k)-\l$ lets invariant  $\cH_{kl}$ and its action on $\psi=\psi_{kl}\otimes Y_{kl}(\l)\in \cH_{kl}$ is given by
$$ \left[\G_1D_x+\gb(x)\G^0+c(x,k)-\l\right]\psi= \left(\left[\G^1D_x+\gb(x)\G^0+c(x,k)-\l\right]\psi_{kl}\right)\otimes Y_{kl}(\l) .$$

Angular operator $A_{\S^2}(\l)$ lets invariant $\cH_{kl}$ and its action on $\psi=\psi_{kl}\otimes Y_{kl}(\l)\in \cH_{kl}$ is given by
\[\lb{3.39mass} A_{\S^2}(\l)\psi=(\mu_{kl}(\l)\,\G^2\,\psi_{kl})\otimes Y_{kl}(\l) . \]

\subsection{ Properties of eigenvalues $\mu_{kl}$ of angular operator for real $\l$.}\lb{ss-eigangmass} 
The angular operator
$$A_{\S^2}(\l)=\cD_{\S^2}+\l \frac{a\sin\vt}{\sqrt{\Delta_\vartheta}}\G^3,$$ where $\cD_{\S^2}$ is given in (\ref{3.10mass}) has similar properties as its massless version considered in Section \ref{ss-eigang}

 For $\l\in\R$ operator $A_{\S^2}(\l)$ is self-adjoint  on
$\cH_{\S^2}=L^2(\S^2;\C^4)$ and has positive discrete spectrum $\s(A_{\S^2}(\l))=\{\mu_{kl}(\l)\}_{(k,l)\in I},$ 
ordered in such a way that for  each $k\in\frac12+\Z$ and $l\in\N^*$ it follows  $0<\mu_{kl}(\l)< \mu_{k(l+1)}(\l).$ Put $A_k(\l)=A_{\S^2}(\l)_{|\cH_{\S^2}^k}.$ 
\vspace{0.5cm}

  Let $\zeta=\frac{a^2\L}{3},$ $\xi=a\l,$ $\nu=a\gm$ and consider operator $A_k(\l)$ as operator-valued function $A_k(\z,\x,\nu)$ of  complex parameters $\z,\x,\nu.$ 
Put
\[\lb{Akmass}A_k(\z,\x,\nu)=A(\z) {\Bbb D}_{\S^2}^k+B(\z,\x,\nu), \]
with $$A(\z)=\sqrt{1+\z\cos^2\vt},\qq B=i\G^2\frac{\z\sin(2 \vt)}{4\sqrt{1+\z\cos^2\vt}}+\G^3\frac{(\z k-\xi)\sin\vt}{\sqrt{1+\z\cos^2\vt}} -\nu\cos\vt\, \G^5 .$$

 Operator $$ A_k(0,0,0)\equiv{\Bbb D}_{\S^2}^k=\G^2D_\vt+\G^3\frac{k}{\sin\vt}$$ is the restriction of the standard Dirac operator on $\S^2$ onto the angular mode $\{e^{ik\vp}\},$ $k\in 1/2+\Z.$
 The domain of $A_k(0,0,0)$ is given by $$\gD=\{u\in \cH_{\S^2}^k,\,\, u\,\,\mbox{is absolutely continuous, }\,\,{\Bbb D}_{\S^2}^ku\in\cH_{\S^2}^k,\,\, u(\pi)=-u(0)\}.$$ 
The spectrum of   $A_k(0,0,0)$ is simple discrete given by
\[\lb{A.3}\mu_{k,l}(0,0,0)=\sgn(l)\left(|k|-\frac12+|l|\right),\qq l\in\Z^*\] and
 $A_k(0,0,0)$ has compact resolvent.

According to ({\ref{2.5-2.6}) we have $\z\in [0,7-4\sqrt{3}]\subset [0,\frac{1}{13.8}]$ and $\x,\nu\in\R$ respectively. Now,
we allow parameters $\z,\x,\nu$ to be complex  $(\z,\x,\nu)\in B(0,\frac{1}{13})\times S^2,$ where
$B(0,r)=\{z\in\C;\,\,|z|<r\}$   and  $S$ is a narrow strip containing the real axis.  

The operators $A(\z),$ $B(\z,\x,\nu)$ are bounded matrix-valued multiplications operators analytic in the variables $(\z,\x,\nu)\in B(0,\frac{1}{13})\times S^2.$ Since the operator $A(\z)$ is also invertible, the operators domain of $A_k(\z,\x,\nu)$ is independent on  $(\z,\x,\nu)\in B(0,\frac{1}{13})\times S^2.$ 

Moreover, since for all $u\in\gD,$ $A_k(\z,\x,\nu)u$ is a vector-valued analytic function in $(\z,\x,\nu),$ and since for $(\z,\x,\nu)\in [0,\frac{1}{13}]\times\R^2$ is self-adjoint on $\cH_{\S^2}^k=L^2( (0,\pi),\,d\vt;\,\C^2),$ then $A_k(\z,\x,\nu)$ forms  a self-adjoint holomorphic family of type (A) in variable
$(\z,\x,\nu)\in B(0,\frac{1}{13})\times S^2$ according to Kato's classification.

Then, using the analytic perturbation theory by Kato \cite{Kato1966},    $A_k(\z,\x,\nu)$ has compact resolvent for all  $(\z,\x,\nu)\in B(0,\frac{1}{13})\times S^2,$ and 
 for a fixed $k\in1/2+\Z,$ the eigenvalues $$\mu_{kl}(\z,\x,\nu),\qq k\in1/2+\Z,\qq l\in\Z^*,$$ of $A_k(\z,\x,\nu)$ are simple and depend holomorphically on $(\z,\x,\nu)$ in a complex \neigh{} of $[0,\frac{1}{13}]\times\R^2.$

Moreover, 
for all $\z\in [0,\frac{1}{13}],$ $k\in\frac12+\Z$ and $l\in\N^*,$
\[\lb{mu0-k}|\mu_{kl}(\z,0)-(|k|-\frac12 +l)|\leq\left(e^\frac{1}{26}-1\right)(|k|-\frac12 +l)+2\left(e^\frac{1}{26}-1\right)\left(1+\frac{1}{26}\right)(|k|+\frac14 )+|a|\gm,\]
\begin{align*}\mu_{kl}(\z,0)&\geq \left(2-e^\frac{1}{26}\right)(|k|-\frac12 +l)-2\left(e^\frac{1}{26}-1\right)\left(1+\frac{1}{26}\right)(|k|+\frac14 )-|a|\gm \\
&\geq \left(4+\frac{1}{13}-(3+\frac{1}{13})e^\frac{1}{26}\right)\left(|k|+\frac14 \right)-|a|\gm >0.8 \left(|k|+\frac14 \right)-|a|\gm>0.7 \left(|k|+\frac14 \right)\end{align*} for $|a|\gm<0.1 \left(|k|+\frac14\right). $
For $l<0$ we use $\mu_{k,-l}=-\mu_{kl}$ and get $\mu_{kl}(\z,0) <-0.7 \left(|k|+\frac14 \right).$

Note that $$B(\z,\x,\nu)-B(\z,0,\nu)=-\G^3\frac{\xi\sin\vt}{\sqrt{1+\z\cos^2\vt}} $$ and the same argument works as in Section 5.2 which implies an analogue of  (\ref{3.5-3.5}) in the massive case: $$ \| u\|_{L^2}\leq \frac{\| A_k(\z,\x) -\o)u\|_{L^2}}{d(\o, \R\setminus (0.7\left(|k|+\frac14\right))(-1,1) )-|a \l|}, $$ provided that the denominator is positive and $|a|\gm$ small enough. 

Note that from (\ref{mu0-k}) it follows and an analogue of (\ref{mubounds}):\\
for all $\l\in\R,$ for all $k\in\frac12 +\Z$ and for all $l\in\N^*,$ the 
eigenvalues $\mu_{kl}(\l)$   for some  constants $C_1$ and $C_2$ independent of $k,l,$  satisfy bound  
\begin{align}&\left(2-e^\frac{1}{26}\right)\left(|k|-\frac12 +l\right)-C_1|k|-C_2-|a|(|\l|+\gm )\leq\mu_{kl}(\l)\lb{muboundsmass}\\&\leq e^\frac{1}{26}\left( |k|-\frac12+l\right)+C_1|k| +C_2+|a|(|\l|+\gm ).\nonumber\end{align}

\section{Resolvent.}\lb{s-resolventmass} 

\subsection{Decomposition in radial and angular parts.}
Now, we apply the method from \cite{Dyatlov2011}, Section 3.
Take $k\in 1/2+\Z$ and an arbitrary $\delta>0.$ Let the angle of admissible contours (see Section \ref{s-resolvent}). Then it follows that the
  resolvents  $$R_r(\l,\omega,k)=(H_k(\l)+\o)^{-1}:\,\,L^2_{\rm comp}(\R,\,dx;\C^4)\,\,\mapsto\,\, H^1_{\rm loc}(\R;\C^4),$$ $$ R_\vt(\l,\omega,k)=(A_k(\l)-\o
)^{-1}:\,\,L^2( (0,\pi),\,d\vt;\,\C^4)\,\,\mapsto\,\, H^1( (0,\pi),\,d\vt;\,\C^4) $$ are {\em meromorphic families} of operators in the sense of Definition 2.1 in \cite{Dyatlov2011}. In particular,  for a fixed values of $\l,$ these families are meromorphic in $\o$ with poles of {\em finite rank} (see Definition 2.2 in \cite{Dyatlov2011}).

We have the following property which assures that admissible contour exists at every regular point.\\
{\em For any compact $K_\l\subset \Omega\subset\C$ there exist constants $C$ and $R$ such that for $\l\in K_\l$ and $|\o|\geq R$,\\
\no (i) for $|\arg\o|\leq\psi$ and $|\pi-\arg\o|\leq\psi$ we have $(\l,\o)\not\in Z_r$ and $\| R_r(\l,\o,k)\|\leq 1/|\o|;$\\
\no (ii) for $\psi\leq |\arg\o|\leq\pi-\psi,$ we have $(\l,\o)\not\in Z_\vt$ and $\|R_\vt(\l,\o,k)\|\leq 1/|\o|.$}

 Then we can construct restriction of the resolvent $$\cR(\l)=(\cD-\l)^{-1}=[\cD(\l)]^{-1}J
=[\tilde{\cD}(\l)]^{-1}[\ga(x)]^{-1}J$$ to $\cH_k$ at any regular point $\l$ as a contour integral
\[\lb{(2.1)}
R(\l,k)=\tilde{R}(\l,k)[\ga(x)]^{-1},\qq \tilde{R}(\l,k)=\frac{1}{2\pi i}\int_\gamma R_r(\l,\omega,k)\otimes R_\vt(\l,\omega,k) d\o\] for some admissible $\g.$ The orientation of $\gamma$ is chosen so that $\G_r$ always stays on the left.

We represent the Hamiltonian as tensor product $H_k(\l)\otimes I_2 +I_2\otimes A_k(\l)$  acting in $L^2(\R,\, dx;\C^4)\otimes \cH_{\S^2}^k,$ where $\cH_{\S^2}^k=L^2( (0,\pi),\,d\vt;\,\C^4).$ Here  $H_k,$ $A_k$ are the operators specified below. 

In the exterior region  $\ga\neq 0$ and we introduce operator
 $$\tilde{\cD}(\l)=[\ga(x)]^{-1}\cD(\l)=[\ga(x)]^{-1}\left(\G^1D_x+\gb(x)\G^0+c(x,k)-\l\right)+ (\cD_{\S^2}-\l b(\vt)\G^3).$$ Denote $\tilde{\cD}_k(\l)=\tilde{\cD}(\l)_{|\cH_k}$ its restriction to $\cH_k.$
For real $\l,$ its radial $\tilde{\cD}_k^r(\l)$ and angular   $\tilde{\cD}_k^\vt(\l)$    parts  let invariant  $\cH_{kl}$ and the action  $\tilde{\cD}_k(\l)=\tilde{\cD}_k^r(\l)+\tilde{\cD}_k^\vt(\l)$  on $\psi=\psi_{kl}\otimes Y_{kl}(\l)\in \cH_{kl}$ is given by
$$\tilde{\cD}_k(\l)\psi= \left(\tilde{\cD}_k^r(\l)\psi_{kl}\right)\otimes Y_{kl}(\l)+(\G^2\psi_{kl})\otimes (A_{\S^2}(\l)Y_{kl}(\l)) .$$

Let $\phi_{kl}=\G^2\psi_{kl}\in L^2(\R,\, dx;\C^4)$ and 
 $$H_k(\l):=[\ga(x)]^{-1}(\G^1 D_x+\gb(x)\G^0+c(x,k)-\l)\G^2,\qq  A_k(\l)=A_{\S^2}(\l)_{|\cH_{\S^2}^k}$$  acting in $L^2(\R,\, dx;\C^4),$  $\cH_{\S^2}^k={\rm Span}_{l\in Z^*}\,  (Y_{kl}(\l))$ 
respectively.
Then for any $\l\in\C,$ $$\tilde{\cD}_k(\l)=H_k(\l)\otimes I_2 +I_2\otimes A_k(\l)$$  acts in  $L^2(\R,\, dx;\C^4)\otimes\cH_{\S^2}^k.$ 

Let $R_r(\l,\omega,k)=(H_k(\l)+\o
)^{-1},$
$R_\vt(\l,\omega,k)=(A_k(\l)-\o
)^{-1}.$

\subsection{Square of the angular operator.}
Here, we extend results of Section \ref{s-red-to-Schr} for the angular operator to the massive case.
 
We consider an $h$-dependent version of the operator $A_{\S^2}(\l)$ $$A_{\S^2,h}(\l)= \sqrt{\D_\vt}\left[\G^2\left(hD_\vt+hq_1(\vt)\right)+\G^3\left(q_2(\vt)hD_\varphi -\l q_3(\vt)\right)-hq_4(\vt)\G^5\right],$$
$$q_1(\vt)= \frac{i\L a^2\sin(2\vt)}{12\Delta_\vt},\qq q_2(\vt)=\frac{E}{\Delta_\vt\sin\vt},\qq q_3(\vt)=\frac{a\sin\vt}{\Delta_\vartheta},\qq q_4(\vt)=\frac{a\gm\cos\vt}{\sqrt{\D_\vt}}. $$

Let $\sigma_j,$ $i=1,2,3,$ be the Pauli matrices (\ref{Paulim}) and $\s_0=I_2.$
We choose the following representation of Dirac matrices satisfying   $\G^i\G^j+\G^j\G^i=2\delta_{ij}I_4,$ $j=0,1,2,3,4,5:$   \begin{align*}&\G^1=\left(
         \begin{array}{cc}
        \s_0 & 0 \\
             0 & -\s_0 
         \end{array}
       \right)=\left(
         \begin{array}{cccc}
        1 & 0& 0&0 \\
             0 & 1& 0&0\\
0 & 0& -1&0\\
0 & 0& 0&-1\\
         \end{array}
       \right),\,\, \G^2=\left(
         \begin{array}{cc}
        0 & -i\s_2 \\
             i\s_2 & 0 
         \end{array}
       \right)=\left(
         \begin{array}{cccc}
        0 & 0& 0&1 \\
             0 & 0& -1&0\\
0 & -1& 0&0\\
1 & 0& 0&0\\
         \end{array}
       \right),\\
 &\G^0=\left(
         \begin{array}{cc}
        0 & -i\s_3 \\
             i\s_3 & 0 
         \end{array}
       \right)=\left(
         \begin{array}{cccc}
        0 & 0& -i& 0 \\
             0 & 0& 0&i\\
i & 0& 0&0\\
0 & -i& 0&0\\
         \end{array}
       \right),\,\, \G^3=\left(
         \begin{array}{cc}
        0 & -i\s_1 \\
             i\s_1 & 0 
         \end{array}
       \right)=\left(
         \begin{array}{cccc}
        0 & 0& 0& -i \\
             0 & 0& -i&0\\
0 & i& 0&0\\
i & 0& 0&0\\
         \end{array}
       \right),
\\
&\G^5=\G^0\G^1\G^2\G^3.\end{align*}

Using that
$$\G^2\G^3=\ma -i\s_3 & 0\\ 0 &-i\s_3\am,\qq \G^2\G^5=\G^0\G^1\G^3=\ma -i\s_2 & 0\\ 0 &i\s_2\am, $$ we get

$$\left[A_{\S^2,h}(\l)\right]^2=  \ma \cP_+& 0_2\\0_2 &\cP_- \am,$$

\begin{align*}\left[\D_\vt\right]^{-1}\cP_\pm=&\left[hD_\vt+hq_1(\vt)\right]^2+\left[q_2(\vt)hD_\vp-\l q_3(\vt)\right]^2-h\s_3\left[q_2'(\vt)hD_\vp -\l q_3'(\vt)\right]\pm h\s_2q_4'(\vt)\\ \\&\hspace{2cm}-h^2q_4^2(\vt). \end{align*}

The leading term of order $h^0$ of $\cP_\pm$ is given by
\begin{align*}& \cP_{+,0}(\l)= \cP_{-,0}(\l)=I_2\left(\D_\vt h^2D_\vt^2+\left[q_2(\vt)hD_\vp-\l q_3(\vt)\right]^2\right)\\
&=I_2\left(\Delta_\vt h^2D_\vt^2+\frac{\left(1+\frac{\L a^2}{3}\right)^2}{\sin^2\vt \D_\vt}\left(hD_\vp -\l a\sin^2\vt\frac{1}{1+\frac{\L a^2}{3}}\right)^2\right) ,\end{align*} where the operator on the diagonal is exactly the same as in the massless case.

\subsection{Radial operator.}

Here we extend the content of Section \ref{s-proofP3.2} to the massive case. 

We consider the radial resolvent
$R_r(\l,\o,k)=(H_k(\l)+\omega)^{-1},$
where
$$H_k(\l):=[\ga(x)]^{-1}(\G^1D_x+\gb(x)\G^0+c(x,k)-\l)\G^2.$$
Note that $R_r(\l,\o,k)$ satisfies
\begin{align*}&\left([\ga(x)]^{-1}(\G^1D_x+\gb(x)\G^0+c(x,k)-\l) \G^2+\o\right)R_r(\l,\o,k)f=f,\qq f\in C_0^\infty (\R,\,dx:\C^2),\\
&\Leftrightarrow\qq \left(\underline{\G^1 D_x+\gb(x)\G^0+c(x,k)-\l+\o \ga(x) \G^2}\right) \G^2 R_r(\l,\o,k)f=\ga(x)f\\
&\Leftrightarrow\qq R_r(\l,\o,k)=\G^2\left[\G^1D_x+\gb(x)\G^0+c(x,k)-\l+\o \ga(x) \G^2\right]^{-1} \ga(x).
\end{align*}

Let  $R(\l,\o,k)=[ \cD_k^r(\l,\o)]^{-1}$ be resolvent  of  $${\displaystyle \cD_k^r(\l,\o)=\G^1D_x+\gb(x)\G^0+c(x,k)-\l+\o \ga(x) \G^2}.$$ 

Here $$c(x,k)=\frac{aEk+\gq Qr}{r^2+a^2},\qq  \ga(x)=\frac{\sqrt{\Delta_r}}{r^2+a^2},\qq \gb(x)=\gm\frac{r\sqrt{\Delta_r}}{r^2+a^2}.$$ satisfying

\[\lb{3.15mass} \ga(x)=a_\pm e^{\k_\pm x}+{\mathcal O}\left(e^{3\k_\pm x}\right),\qq \gb(x)=b_\pm e^{\k_\pm x}+{\mathcal O}\left(e^{3\k_\pm x}\right)\qq x\rightarrow\pm\infty,\]
\[\lb{3.17mass}c(x,k)=\Omega_\pm(k)+c_\pm e^{2\k_\pm x}+{\mathcal O}\left(e^{4\k_\pm x}\right),\qq x\rightarrow\pm\infty,\]
$$\Omega_\pm=\frac{aEk+\gq Qr_\pm}{r_\pm^2+a^2},\qq \O_- >\O_+.$$ 

Outgoing solutions $f^\pm$ can be defined similarly to the massless case, using that as $x\rightarrow\pm\infty$ they satisfy $ {\displaystyle \left[\G^1D_x+\Omega_\pm(k)-\l\right]}f=0,$ 
$$f(x)= \sum_{i=1}^4 f_i(x)e_i,  \ \ \ e_1= \ma 1\\0\\0\\0\am,
   \   e_2= \ma 0\\1\\0 \\0\am,\ldots$$
 or
$$
 \left\{\begin{array}{c}
           -if_1'+\Omega_\pm(k)f_1=\l f_1 \\
-if_2'+\Omega_\pm(k)f_2=\l f_2 \\
           if_3'+\Omega_\pm(k)f_3=\l f_3\\
  if_4'+\Omega_\pm(k)f_4=\l f_4
         \end{array}\right. ,\qq \l \in \C.
$$ The outgoing solutions $f^\pm$ are then defined by their asymptotics
 $$f^-\sim f^{0,-}(x;\l,k)= \ma 0_2\\e^{-i(\l-\Omega_-(k)) x}I_2\am,\qq x\rightarrow-\infty;$$
$$ f^+\sim f^{0,+}(x;\l,k)=  \ma e^{i(\l-\Omega_+(k)) x} I_2\\ 0_2\am,\qq x\rightarrow+\infty.$$
Then the properties of the radial resolvent follows as in the massless case, Section \ref{s-proofP3.2}.

\section{Semiclassical leading part and square of the radial operator.}
Here, we extend results of Section \ref{s-red-to-Schr} for the radial operator to the massive case.

We consider semiclassical version of the Dirac operator $\G^1D_x+\gb(x)\G^0+c(x,k)-\l+\o \ga(x) \G^2$ of the form  $$\G^1hD_x+h\gb(x)\G^0+c(x,\tilde{k})-(\tilde{\l}+ih\tilde{\nu})+(\tilde{\o} +ih\tilde{\mu}) \ga(x) \G^2$$ with the leading part $\cD_{h}=\cD_{h,+}=\G^1hD_x+c_0(x,\tilde{k})-\tilde{\l}+\tilde{\o}  \ga(x) \G^2.$ Note that   the leading part  is independent of the field mass $\gm.$ Put $c=c_0,$ $q=\tilde{\o}  \ga(x).$

In $\cH=L^2(\R )^4$ we consider  two Dirac operators $$\cD_{h,\pm}(\l):=-\cD_{h,0}\pm (c(x)-\l)I_4,\qq \cD_{h,0}:=-ih\G^1\partial_x+q(x)\G^2.$$ 

The product of the operators is given by
\[\lb{produktmass}\cD_{h,+}\cD_{h,-}=  \cD_{h,0}^2-(c-\l)^2 I_4-\left[\cD_{h,0},(c-\l) I_4\right].\]Here
$$ \cD_{h,0}^2=(-h^2\partial_x^2+q^2(x))I_4-\ma 0&\s_2\\ \s_2 &0\am hq'(x)$$  is the  matrix Schr{\"o}dinger operator. The commutator is given by
$$\left[\cD_{h,0},(c-\l) I_4\right]= -ih\G^1 c'(x).$$
The operator $\cD_{h,0}^2$ is self-adjoint in  $\cH=L^2(\R )^4$ and unitary equivalent to
\[\lb{diagonmass}U\cD_{h,0}^2U^{-1}=\ma \cP_{h,0}^-&0\\
0&\cP_{h,0}^+\am,\qq \cP_{h,0}^\pm=\left(-h^2\partial_x^2 +q^2\pm h q'\right) I_2.\]
Here $U$ is constant matrix. 

 We get also
\[\lb{diagon_bismass}U\left(\cD_{h,0}^2-(c-\l)^2 I_4\right)U^{-1}=\ma \cP_h^-(\l)&0\\
0&\cP_h^+(\l)\am,\]
$$ \cP_h^\pm(\l)=\left(-h^2\partial_x^2 +q^2-(c-\l)^2\pm h q'\right) I_2$$

Now, using (\ref{produktmass}) we get

\[\lb{produkt-bismass} \cD_{h,+}\cD_{h,-}= U^{-1}\ma \cP_{h}^-(\l)&0\\
0&\cP_{h}^+(\l)\am U +ih c'(x)\G^1.
\] If $\cP_h^\pm(\l)$ are invertible, denote $$\cR(\l)=U^{-1}\ma [\cP_{h}^-(\l)]^{-1}&0\\
0&[\cP_{h}^+(\l)]^{-1}\am U$$ and write
$$\cD_{h,+}\cD_{h,-}\cR(\l)=  I_4 +ih c'(x)\G_1\cR(\l)$$
$$\cD_{h,-}\cR(\l)=  [\cD_{h,+}]^{-1} + [\cD_{h,+}]^{-1}ih c'(x)\G^1\cR(\l)=[\cD_{h,+}]^{-1}\left(I_4 + ih c'(x)\G^1\cR(\l)\right) $$
or
\begin{align*}&[\cD_{h,+}(\l)]^{-1}\left(I_4+ih\G^1 c'(x)\left(\cD_{h,0}^2-(c-\l)^2 I_4\right)^{-1}\right)=\cD_{h,-}(\l)\left(\cD_{h,0}^2-(c-\l)^2 I_4\right)^{-1}\end{align*}
\begin{align*}&[\cD_{h,+}(\l)]^{-1}=\cD_{h,-}(\l)\left(\cD_{h,0}^2-(c-\l)^2 I_4\right)^{-1}\left(I_4+ih\G^1 c'(x)\left(\cD_{h,0}^2-(c-\l)^2 I_4\right)^{-1}\right)^{-1}\end{align*}
which leads to
\begin{align}&[\cD_{h,+}(\l)]^{-1}=\cD_{h,-}(\l)U^{-1}\ma [\cP_h^-(\l)]^{-1}&0\\
0&[\cP_h^+(\l)]^{-1}\am U\lb{rad-res-id}\\&\cdot\left(I_4+ih\G^1 c'(x)U^{-1}\ma [\cP_h^-(\l)]^{-1}&0\nonumber\\
0&[\cP_h^+(\l)]^{-1}\am U\right)^{-1}.\end{align}

As $U$ is constant matrix, identity (\ref{rad-res-id}) can be extended to $q=\o\ga$ complex with $\o\in\C.$ 

\section{Normal form of the radial operator.}\lb{s-normal formmass}
In this section we extend Proposition \ref{Prop4.3} to the massive case which together with the properties of the angular operator and quantization conditions as in the massless case concludes the proof of Theorem \ref{Th1}.
 
Let $\tilde{\cD}_k^r=\G^1hD_x+h\gb(x)\G^0+c(x,\tilde{k})-(\tilde{\l}+ih\tilde{\nu})+(\tilde{\o} +ih\tilde{\mu}) \ga(x) \G^2$ with
the leading part $$\cD_{h,0}=\G^1hD_x+c_0(x,\tilde{k})-\tilde{\l}+\tilde{\o}  \ga(x) \G^2=$$ 
$$\left(
         \begin{array}{cccc}
        hD_x+c_0(x,\tilde{k})-\tilde{\l} & 0& 0&\tilde{\o}  \ga(x) \\
             0 & hD_x+c_0(x,\tilde{k})-\tilde{\l}& -\tilde{\o}  \ga(x)&0\\
0 & -\tilde{\o}  \ga(x)& -hD_x+c_0(x,\tilde{k})-\tilde{\l}&0\\
\tilde{\o}  \ga(x) & 0& 0&-hD_x+c_0(x,\tilde{k})-\tilde{\l}\\
         \end{array}
       \right).$$ The determinant of its principal symbol
$$\det\left(
         \begin{array}{cccc}
        \xi+c_0(x,\tilde{k})-\tilde{\l} & 0& 0&\tilde{\o}  \ga(x) \\
             0 & \xi+c_0(x,\tilde{k})-\tilde{\l}& -\tilde{\o}  \ga(x)&0\\
0 & -\tilde{\o}  \ga(x)& -\xi+c_0(x,\tilde{k})-\tilde{\l}&0\\
\tilde{\o}  \ga(x) & 0& 0&-\xi+c_0(x,\tilde{k})-\tilde{\l}\\
         \end{array}
       \right)$$ is equal to
$(\xi^2+(\tilde{\o}  \ga(x))^2-(c_0-\tilde{\l})^2)^2.$ Thus eigenvalues of the principal symbol are of multiplicity two and equal to $$\mu_\pm=\mu_\pm(x,\xi; \tilde{\l},\tilde{\o},\tilde{k})=c_0(x,\tilde{k})-\tilde{\l}\pm\sqrt{\xi^2+\tilde{\o}^2\ga^2(x)}, $$ which are the same as in the massless case (where the multiplicity was one).

Now, by iteration as in Section 3.1 of \cite{HelfferSjostrand1990} and Section 4.1 of \cite{BruneauRobert1999} we can reduce operator $\tilde{\cD}_k^r$ to the block-diagonal form up to the infinite order of $h.$

\begin{proposition}[Decoupling]\lb{Prop-decouplingmass} There exist unitary $U=U(x,hD_x,h)$  such that 
\[\lb{decouplingmass} U^*\tilde{\cD}_k^r(\l,\o) U= \left(\begin{array}{cc}
   \mu_{+}(x,hD,h)&0 \\
   0  &   \mu_{-}(x, hD,h)\\
         \end{array}
       \right)+{\mathcal O}(h^\infty)\] in $\cL(L^2,L^2).$ Here
 $\mu_{\pm}(x,hD,h)$ is \pseudor{} with classical symbol $\mu_{\pm}(x,\xi,h)=\sum_{j=0}^{\infty}h^j\mu_{j,\pm}(x,\xi)$ and $\mu_{0,\pm}(x,\xi)=\mu_{\pm}(x,\xi)I_2.$
\end{proposition}
 \begin{proposition}\lb{Prop4.3mass} The microlocal normal form of $\tilde{\cD}_k^r(\l,\o)$ is given by
 $$\left(\begin{array}{cc}
 hxI_2D_x-\beta_+&0 \\
   0  &  -hxI_2D_x+\beta_-\\
         \end{array}
       \right) ,$$ where $\beta_\pm=\beta_\pm(\tilde{\l},\tilde{\nu},\tilde{\o},\tilde{\mu},\tilde{k};h)$ is classical symbol. Moreover, the principal part $\beta_{0,\pm}$ of $\beta_\pm$ is real-valued, independent of $\tilde{\nu},$ $\tilde{\mu}$ and mass $\gm,$ and vanishes if and only if $\tilde{\o}=\tilde{\o}_0(\tilde{\l},\tilde{k}),$  where $\tilde{\o}_0^{-2}$ is defined by   (\ref{valf}). Moreover,
$$\beta_{0,+}=-\frac{(c_0(x_+)-\tilde{\l}+ \tilde{\o}\ga (x_+))\sqrt{|\tilde{\o}|\ga(x_+)}}{\sqrt{ \left|(c_0+|\tilde{\o}|\ga)''(x_+)\right|}}I_2+{\mathcal O}\left( c_0(x_+)-\tilde{\l}+ \tilde{\o}\ga(x_+)\right)^2 $$
and  $\beta_{0,-}(\tilde{\l},\tilde{\nu},\tilde{\o},\tilde{k})=\beta_{0,+}(-\tilde{\l},\tilde{\nu},\tilde{\o},-\tilde{k}).$
\end{proposition}
Now, using that the leading parts of the angular and radial operators are independent of mass $\gm,$ the angular and radial quantization conditions are independent are independent of mass in the leading order and Theorem (\ref{Th1}) follows. 

\newpage


\end{document}